\documentclass[leqno,final]{siamart171218}

\usepackage{amsmath,dsfont}
\usepackage{graphicx,color}
\usepackage{amssymb}

\usepackage{algorithm}
\usepackage{algorithmic}
\usepackage{float}
\usepackage{enumitem}   

\def\E{\mathbb{E}}

\newcommand{\Div}{\mbox{{\rm div}\,}}
\newcommand{\bu}{\mathbf{u}}

\newcommand{\bv}{\mathbf{v}}

\newcommand{\bF}{\mathbf{F}}
\newcommand{\bG}{\mathbf{G}}
\newcommand{\bW}{\mathbf{W}}
\newcommand{\bw}{\mathbf{w}}

\newcommand{\bL}{\mathbf{L}}

\newcommand{\bH}{\mathbf{H}}

\newcommand{\bUh}{\mathbb{U}_h}
\newcommand{\bVh}{\mathbb{V}_h}
\newcommand{\bphi}{\boldsymbol{\phi}}
\newcommand{\bpsi}{\boldsymbol{\psi}}
\newcommand{\bPhi}{\boldsymbol{\Phi}}

\newcommand{\td}{\text{\rm d}}

\newcommand{\Rh}{\cal{R}_h}
\newcommand{\Ph}{\cal{P}_h}
\newcommand{\beu}{{\bf e}_{\bu}}
\newcommand{\bbeu}{{\bf \overline{e}}_{\bu}}
\newcommand{\bev}{{\bf e}_{\bv}}
\newcommand{\tJ}{\widetilde{\mathcal{J}}}

\def\mbf#1{\mathbf{#1}}
\def\bb#1{\mathbb{#1}}
\def\cal#1{\mathcal{#1}}
\def\exc#1{\bb{E}\left[#1\right]}

\def\normbf#1{\left\lVert{#1}\right\rVert^{2}_{\mbf{L}^2}}
\def\normbfhalf#1{\left\lVert{#1}\right\rVert_{\mbf{L}^2}}
\def\normbfp#1{\left\lVert{#1}\right\rVert^{2p}_{\mbf{L}^2}}
\def\normbfh#1{\left\lVert{#1}\right\rVert^{2}_{\mbf{H}^1}}

\def\half{\frac{1}{2}}
\def\veps{\varepsilon}

\def\normdel#1{\lambda\normbf{\Div(#1)}+\mu\normbf{\veps(#1)}}

\def\normdelbig#1{\lambda\normbf{\Div\big(#1\big)}+\mu\normbf{\veps\big(#1\big)}}

\def\normdelpbig#1{\lambda\normbfp{\Div\big(#1\big)}+\mu\normbfp{\veps\big(#1\big)}}

\def\innerprd#1#2{\lambda\big(\Div({#1}),\Div({#2})\big)+\mu\big(\veps({#1}),\veps({#2})\big)}

\def\bferr#1#2{\mbf{e}_\mbf{{#1}}^{#2}}

\def\bfErr#1#2{\mbf{E}_\mbf{{#1}}^{#2}}

\newcommand{\opL}{\cal{L}}
\newcommand{\opLh}{\cal{L}_{h}}

\newcommand\commentone[1]{\textcolor{black}{#1}}

\newtheorem{remark}{Remark}
\newcommand{\mE}{\mathbb{E}}

\title{Optimal Order Space-Time Discretization Methods for the Nonlinear Stochastic Elastic Wave Equations with Multiplicative Noise}

\author{
	Xiaobing Feng\thanks{Department of Mathematics, The University of Tennessee, Knoxville, TN 37996, U.S.A. This author was partially supported by the NSF grant DMS-2309626. ({\tt xfeng@utk.edu}).}  
	\and
	Yukun Li\thanks{Department of Mathematics, University of Central Florida, Orlando, FL, 32816, U.S.A. This author was partially supported by the NSF grant DMS-2110728. ({\tt yukun.li@ucf.edu}).}
	\and
	Liet Vo\thanks{School of Mathematical and Statistical Sciences, The University of Texas Rio Grande Valley, Edinburg, TX 78539, U.S.A.  ({\tt liet.vo@utrgv.edu}).} 
}

\begin{document}
	\maketitle
	
	\begin{abstract}
		This paper develops and analyzes an optimal-order semi-discrete scheme and its fully discrete finite element approximation for nonlinear stochastic elastic wave equations with multiplicative noise. 
		A non-standard time-stepping scheme is introduced for time discretization, it is showed that the scheme converges with rates $O(\tau)$ and $O(\tau^{\frac32})$ respectively in the energy- and  $L^2$-norm,  which are optimal with respect to the time regularity of the PDE solution. 
		 For spatial discretization, the standard finite element method is employed. It is proven that the 
		 fully discrete method  converges with optimal  rates $O(\tau + h)$ and $O(\tau^{\frac{3}{2}} + h^2)$  respectively in the energy- and $L^2$-norm.
		 The cruxes of the analysis are to establish some high-moment stability results and utilize a refined error estimate for the trapezoidal quadrature rule to control the nonlinearities from the drift term and the multiplicative noise. Numerical experiments are also provided to validate the theoretical results.
		
	\end{abstract}
	
	

	\section{Introduction}\label{sec-1}
	This paper is concerned with both semi-discrete and fully discrete approximations for the following stochastic elastic wave equations with multiplicative noise:
	\begin{subequations} \label{eq1.1}
	\begin{alignat}{2}\label{eqn:elastic-wave-equation}
		\bu_{tt}-\Div(\sigma(\bu))&=\bF[\bu]+\bG[\bu]\xi\quad &\text{in}&\  \cal{D}_{T}:=(0,T)\times\cal{D},\\
		\bu(0,\cdot)&=\bu_{0},\quad\bu_{t}(0,\cdot)=\bv_{0}\quad &\text{in}&\ \cal{D},\label{20220526_1}\\
		\bu(t,\cdot)&=0&\text{on}&\ \partial\cal{D}_{T}:=(0,T)\times\partial\cal{D},\label{20220526_2} 
	\end{alignat}
	where 	$\bu_t=\frac{\td\bu}{\td t}$, $\xi = \dot{\bW} = \frac{\td \bW}{\td t}$ is the white noise, $\cal{D}\subset\bb{R}^{d}\ (d=2,3)$ is a bounded domain, $(\bu_{0},\bv_{0})$ is a pair of $\bf{H}_{0}^{1}\times\bf{L}^{2}$-valued random variables, and
	\begin{align}
		\sigma(\bu)&=(\lambda \Div\bu)\bf{I}+\mu\veps(\bu),\label{20220526_3}\\
		\veps(\bu)&=\frac{1}{2}\bigl(\nabla\bu+(\nabla\bu)^{\bf{T}}\bigr),\label{20220526_4}\\
		\bF[\bu]&=\bF(\bu,\nabla \bu),\label{20220526_5}\\
		\bG[\bu]&=\bG(\bu,\nabla \bu).\label{20220526_6}
	\end{align}
	\end{subequations}

	Here $\bf{I}$ stands for the identity matrix, $\bF[\bu]$ and $\bG[\bu]$ are two d-dimensional nonlinear mappings (see Section \ref{sec-2} for the precise definitions), and $W$ is a $\bb{R}$-valued Wiener process defined on the filtered probability space $(\Omega,\cal{F},\{\cal{F}_{t}\}_{0\leq t\leq T},\bb{P})$. We note that this simpler setting is considered to avoid technicalities but the results of this paper can be extended to the more general cases when $\bW$ is a $\bb{R}^d$-valued Wiener process and $\bG[\bu]$ is a scalar nonlinear mapping and when $\bW$ is a $\bb{R}^{l}$-valued  Wiener process and $\bG[\bu]$ is a $d\times l$ matrix.
	
	Wave propagation as one of primary means for energy transmission arises in broad applications from many scientific and engineering fields including biology, physics, earth sciences, medical science, telecommunications, and defense industry. For the mathematical analysis of various wave equations, we refer readers to \cite{chow2002stochastic, chow2006asymptotics, chow2009nonlinear, chow2015stochastic, dalang2009stochastic, dalang1998stochastic, millet2000stochastic, millet1999stochastic} and the references therein. The numerical analysis of wave propagation problems, in both deterministic and stochastic contexts, has also been extensively studied.
	In the deterministic setting, finite element methods have been well developed with rigorous error estimates, as demonstrated in \cite{baker1976error, dupont19732, falk1999explicit}. Similarly, discontinuous Galerkin methods have been studied in \cite{adjerid2011discontinuous, baccouch2012local, chou2014optimal, chung2006optimal, chung2009optimal, grote2006discontinuous, grote2009optimal, monk2005discontinuous, riviere2003discontinuous, safjan1993high, sun2021, xing2013energy, zhong2011numerical}, not even mention those developments in finite difference, spectral and other types of methods. In the stochastic setting, substantial progress has also been achieved.  We refer to \cite{cohen2013trigonometric, cui2019strong, gubinelli2018renormalization, hausenblas2010weak, kovacs2010finite} for studies related to additive noise, and \cite{Anton2016full, cohen2018numerical, cohen2015fully, feng2022higher, feng2024optimal, hong2021energy, li2022finite, kovacs2010finite, quer2006space} for results with multiplicative noise.
	
	While the above listed references were mostly concerned with approximating acoustic wave equations, numerical methods for deterministic elastic wave equations have also been studied in parallel in \cite{ciarlet2002finite, gauthier1986two, igel1995anisotropic, marfurt1984accuracy, reddy1973convergence, saenger2000modeling, virieux1984sh}, and for stochastic elastic wave equations, including those involving random field coefficients and sources, in \cite{babuvska2014stochastic, feng2022fully, motamed2015analysis, motamed2013stochastic}. 
	Relating to this paper, we note that several semi-discrete in-time schemes for stochastic acoustic wave equations were proposed and analyzed in \cite{feng2022higher}, rates of convergence $O(\tau^{\frac12})$ in the energy norm and $O(\tau^{\frac32})$ in the $L^2$-norm were established. Two additional semi-discrete schemes for nonlinear stochastic acoustic wave equations were proposed in \cite{feng2024optimal} and proved to have rates of convergence $O(\tau)$ in the energy norm and $O(\tau^{\frac32})$ in the $L^2$-norm. However, both papers did not analyze the fully discrete schemes. 
	Later,  both semi-discrete in-time  and fully discrete finite element methods were proposed and analyzed 
	for the stochastic elastic wave equations with multiplicative noise in  \cite{feng2022fully}.  Rates of convergence $O(\tau^{\frac12} + h)$ in the energy norm and $O(\tau + h^2)$ in the $L^2$-norm were proved for the fully discrete methods. However,  these rates of convergence in time  is lower than the optimal rate obtained in \cite{feng2022higher, feng2024optimal}. Moreover, the error estimates for the fully discrete method were established only for linear stochastic elastic wave equations under the assumption that the multiplicative noise depends on the PDE solution linearly. 
	
	The primary goal of this paper is to develop and analyze an optimal-order semi-discrete scheme in time and its fully discrete finite element approximation for the nonlinear stochastic elastic wave problem \eqref{eqn:elastic-wave-equation}--\eqref{20220526_6}. We extend the optimal error estimate result of \cite{feng2024optimal} for the acoustic waves to the nonlinear stochastic elastic wave equations with multiplicative noise. For spatial discretization, we only consider the standard finite element method in this paper although other spatial discretization methods can be handled similarly. We are able to establish space-time error estimates for the proposed methods with optimal rate of convergence $O(\tau + h)$ in the energy norm and $O(\tau^{\frac{3}{2}} + h^2)$ in the $L^2$-norm. 
	To achieve these results, the cruxes of the analysis are to establish some high-moment stability results and utilize a refined error estimate for the trapezoidal quadrature rule to control the nonlinearities from the drift term and the multiplicative noise.  We also note that our numerical solutions retain several important properties of the PDE solution, such as energy stability and high-moment bounds. This is because our numerical methods are designed with the help of the structure of the stochastic elastic wave equations, where high-order terms arising from noise are incorporated into the methods, rather than being discarded as done previously in the literature. This approach ensures that the numerical solutions approximate the PDE solution with higher accuracy.
	
	The remainder of the paper is organized as follows.  In Section \ref{sec-2}, we introduce notation, assumptions on the nonlinearities, and quote some facts about the Gaussian noise increments. In Section \ref{sec-3}, we formulate our nonstandard semi-discrete time-stepping scheme and prove its stability and optimal order error estimates. In Section \ref{sec-4}, we formulate the finite element spatial discretization and establish its stability and optimal order error estimates. Section \ref{sec-5} states our fully discrete error estimates, which are obtained by combining the temporal and spatial error estimates from Sections \ref{sec-3} and \ref{sec-4}. Finally, Section \ref{sec-6} contains some numerical experiment results which validate our theoretical error estimates.
	
	\section{Preliminaries}\label{sec-2}
	\subsection{Notations}
	Standard notations for functions and spaces are adopted in this paper. For example, $\bL^p$ denotes $(L^p(\cal{D}))^d$ for $1\leq p\leq \infty$, $(\cdot,\cdot)$ denotes the standard $\bL^2(\cal{D})$-inner product and $\bH^m(\cal{D})$ denotes the Sobolev space of order $m$. Throughout this paper, $C$ will denote a generic positive constant independent of the mesh parameters $\tau$ and $h$. 
	
	Moreover, let $(\Omega,\mathcal{F}, \{\mathcal{F}_t\},\mathbb{P})$ be a filtered probability space with the probability measure $\mathbb{P}$, the 
	$\sigma$-algebra $\mathcal{F}$ and the continuous filtration $\{\mathcal{F}_t\} \subset \mathcal{F}$. For a random variable $v$ 
	defined on $(\Omega,\mathcal{F}, \{\mathcal{F}_t\},\mathbb{P})$,
	${\mathbb E}[v]$ denotes the expected value of $v$. 
	For a vector space $X$ with norm $\|\cdot\|_{X}$,  and $1 \leq p < \infty$, we define the B\"ochner space
	$\bigl(L^p(\Omega;X); \|v\|_{L^p(\Omega;X)} \bigr)$, where
	$\|v\|_{L^p(\Omega;X)}:=\bigl({\mathbb E} [ \Vert v \Vert_X^p]\bigr)^{\frac1p}$.
	Let $\mathcal{K}, \mathcal{H}$ be two separable Hilbert spaces, we use 
	$\mathcal{L}(\mathcal{K}^m, \mathcal{H})$ to denote the space of all multi-linear maps from $\mathcal{K} \times \cdots \times \mathcal{K}$ ($m$-times) to $\mathcal{H}$ for $m \geq 1$.  For a mapping $\Psi: \bH_0^1  \rightarrow {\bL}^2$, we introduce the notation $D_{\bu} \Psi(\bu) \in \mathcal{L}\left({\bH}_0^1, {\bL}^2\right)$ for the G\^ateaux derivative with respect to  $\bu$, whose action is seen as
	$$
	\pmb{\xi} \mapsto D_{\bu} \Psi(\bu)(\pmb{\xi}) \qquad \text { for } \pmb{\xi} \in {\bH}_0^1(D) .
	$$

	Next, 	let $N>>1$ be a positive integer and $\tau:=T/N$. Let $\{t_n\}_{n=0}^N$ be a uniform partition of  the interval $[0,T]$ with mesh size $ \tau$. 	We define the useful notation $\overline{\Delta}W_{n}$ and $\widetilde{\Delta}W_n$ as follows:
	\begin{align*}
		\overline{\Delta}W_{n} &:= W(t_{n+1}) - W(t_n), \\ 
		\widetilde{\Delta}W_n &:= \int_{t_n}^{t_{n+1}}[W(t_{n+1})-W(s)]ds\\
		&= \tau W(t_{n+1}) - \int_{t_n}^{t_{n+1}} W(s)\, ds. 
	\end{align*}
	In addition, following  \cite{feng2022higher},  we approximate the last integral with higher accuracy as follows:
	\begin{align}\label{approx_integral}
		\int_{t_n}^{t_{n+1}} W(s)\, ds \approx \tau^3 \sum_{\ell=1}^{[\tau^{-2}]} W(t_{n,\ell}),
	\end{align}
	where $\{W({t_{n,\ell}})\}_{\ell=1}^{[\tau^{-2}]}$ is the piecewise affine approximation of $W$ on  $[t_n, t_{n+1}]$ over an equidistant mesh $\{t_{n,\ell}\}_{\ell=1}^{[\tau^{-2}]}$ with mesh size $\tau^3= t_{n,\ell + 1} - t_{n,\ell}$.
	
	Thus, we obtain the following high order approximation of $\widetilde{\Delta}W_n$ (see Lemma \ref{increment}):
	\begin{align}
		\widehat{\Delta} W_n &:= \tau W(t_{n+1}) - \tau^3 \sum_{\ell = 1}^{[\tau^{-2}]} W(t_{n,\ell}).
	\end{align}
	We will frequently write $	\widetilde{\Delta} W_n \approx \widehat{\Delta} W_n$ in the subsequent sections. 
	
	\subsection{Assumptions}\label{ssec:assumptions}
	The following structural conditions will be imposed on the mappings $\bF[\cdot]$ and $\bG[\cdot]$:
	\begin{subequations}\label{Assumptions}
	\begin{gather}
		\|\bF[0]\|_{\bL^2}+\|\bG[0]\|_{\bL^2}\leq C_A,\label{assump:F(0)}\\
		\|D_{\bu}\bF[\cdot]\|_{\bL^\infty}+\|D_{\bu} \bG[\cdot]\|_{\bL^\infty}\leq C_A,\label{assump:GradientF}\\
		\|\bF_{\bu_i\bu_j}[\cdot]\|_{\bL^\infty} +\|\bG_{\bu_i\bu_j}[\cdot]\|_{\bL^\infty}  \leq C_A,\quad 1\le i,j\le d,\label{assump:Fuu}\\
		\normbfhalf{\bF[\bv]-\bF[\bw]}\leq C_B\Bigl(\lambda\|\mathrm{div}(\bv-\bw)\|_{{\bf L}^2}^2+\mu\|\epsilon(\bv-\bw)\|_{{\bf L}^2}^2+\normbfhalf{\bv-\bw}^2\Bigr)^\half,\label{assump:LipsF}\\
		\normbfhalf{\bG[\bv]-\bG[\bw]}\leq C_B\Bigl(\lambda\|\mathrm{div}(\bv-\bw)\|_{{\bf L}^2}^2+\mu\|\epsilon(\bv-\bw)\|_{{\bf L}^2}^2+\normbfhalf{\bv-\bw}^2\Bigr)^\half,\label{assump:LipsG}
	\end{gather}
\end{subequations}
	where $\bF_{\bu_i\bu_j}[\cdot]$ denotes the second derivative of $\bF$ with respect to $\bu_i, \bu_j$ , and $C_A$ and $C_B$ are two positive constants.  We note that, by the Korn's inequality,  \eqref{assump:LipsF} and \eqref{assump:LipsG} imply $\bF$ and $\bG$ are Lipschitz continuous from $\bH^1$ to $\bL^2$. 
	
	\subsection{Useful facts} 
	In this subsection, we collect a few useful facts that will be crucially used in the proofs of our convergence results in the subsequent sections.  The first one is the following refined estimate for the trapezoidal quadrature rule \cite[Theorem 2]{Dragomir}.
	
	\begin{lemma}
		Let $\phi \in C^{1,\alpha} ([0, T]; \mathbb{R})$ for some $\alpha \in (0, 1]$. Then, there exists a constant $C_e >0$ such that
		\begin{align*}
			\biggl|\frac{\phi(0)  + \phi(T)}{2} - \frac{1}{T} \int_0^T \phi(\xi)\, d\xi\biggr| \leq \frac{C_e}{(\alpha + 2)(\alpha + 3)} \, T^{1+\alpha},
		\end{align*}
		where the constant $C_e$ satisfies
		\begin{align*}
			|\phi'(t) - \phi'(s)| \leq C_e |t-s|^{\alpha}.
		\end{align*}
	
	\end{lemma}
	
	Next, we state a few properties of $\overline{\Delta}W_{n}$, $\widetilde{\Delta}W_n$, and $ \widehat{\Delta} W_n$ which are defined in the previous subsection.  Their proofs can be found in \cite[Remarks 1 and 2]{feng2022higher}.
	
		\begin{lemma}\label{increment} 
			There hold the following estimates for the increments $\widetilde{\Delta}W_n$ and $ \widehat{\Delta} W_n$:
			\begin{enumerate}
				\item[\rm (i)] $\mE\bigl[|\widetilde{\Delta}W_n - \widehat{\Delta}W_n|^2\bigr] \leq C\tau^5$.
				\item[\rm (ii)] $\mE\bigl[|\widetilde{\Delta}W_n|^{2p}\bigr] + \mE\bigl[|\widehat{\Delta}W_n|^{2p}\bigr] \leq C\tau^{3p}$ for $1\leq p <\infty$.
			\end{enumerate}
		\end{lemma}

		\subsection{Variational weak formulation and properties of weak solutions}\label{ssec:weakform}
		In this subsection, we first give the definition of variational weak formulation and weak solutions for problem \eqref{eqn:elastic-wave-equation}--\eqref{20220526_6}. We then establish several technical lemmas that will be used in the subsequent sections.
		
		Equations \eqref{eqn:elastic-wave-equation}--\eqref{20220526_2} can be written as
		\begin{subequations}\label{eq2.4}
			\begin{align}\label{eqn:elastic-wave-equation-uv}
				\td\bu&=\bv\td t,\\
				\td\bv&=(\opL\bu+\bF[\bu])\commentone{\td t}+\bG[\bu]\td W(t),\quad \opL\bu := \Div \sigma(\bu), 
				\label{20220526_8}\\
				\bu(0)&=\bu_{0},\quad\bv(0)=\bv_{0},\label{20220526_9}\\
				\bu(t,\cdot)&=0.\label{20220526_10}
			\end{align}
		\end{subequations}

		\begin{definition}\label{def:weakform-solution}
			The weak formulation for problem \eqref{eqn:elastic-wave-equation-uv}--\eqref{20220526_10} 
			is defined as seeking $(\bu, \bv)\in\bL^2\bigl(\Omega;\mbf{C}([0, T]; \commentone{\bL^2})\cap\bL^2((0, T),\mbf{H}^1_0)\bigr)\times\bL^2\bigl(\Omega;\mbf{C}([0, T], \bL^2)\bigr)$ such that 
			\begin{subequations}
			\begin{alignat}{2}
				\big(\bu(t), \bphi\big)=&\int_{0}^{t}\big(\bv(s), \bphi\big)\td s + (\bu_0, \bphi)&&\forall \bphi \in \bL^2,\label{eqn:weakform-1}\\
				\big(\bv(t), \bpsi\big)=&-\int_{0}^{t}\lambda\big(\Div(\bu(s)), \Div(\bpsi)\big)\td s-\int_{0}^{t}\mu\big(\veps(\bu(s)), \veps(\bpsi)\big)\td s\label{eqn:weakform-2}\\
				+\int_{0}^{t}&\big(\bF[\bu(s)], \bpsi\big)\td s+\int_{0}^{t}\big(\bG[\bu(s)]\td W(s), \bpsi\big) + (\bv_0, \bpsi)&&\forall \bpsi \in \bH^1_0 \nonumber
			\end{alignat}
			\end{subequations}
			for all $(\phi, \psi)\in\bf{L}^2\times\bf{H}^1_0$. Such a pair $(\bu,\bv)$, if it exists,  is called
			a (variational) weak solution to problem  \eqref{eqn:elastic-wave-equation-uv}--\eqref{20220526_10}. Moreover, if a weak solution $(\bu, \bv)$ belongs to $\bL^2\big(\Omega; \mbf{C}([0, T];$ $\bH^2\cap\bH^1_0)\big)\times\bL^2\big(\Omega; \mbf{C}([0, T];\bH^1_0)\big)$, then $(\bu, \bv)$ is called a strong solution to problem  \eqref{eqn:elastic-wave-equation-uv}--\eqref{20220526_10}.
		\end{definition}
		
		We recall the following stability estimates from \cite[Lemma 2.3]{feng2022fully}.
		\begin{lemma}\label{lem:stability}
			Let $(\bu, \bv)$ be a strong solution to equations \eqref{eqn:elastic-wave-equation-uv}--\eqref{20220526_10}. Under the assumptions \eqref{assump:F(0)}--\eqref{assump:LipsF},
			there hold
			\begin{subequations}
			\begin{align}
				\sup\limits_{0\leq t\leq T}\exc{\normbf{\bv}}+\sup\limits_{0\leq t\leq T}\exc{\normdel{\bu}}\leq& C_{s1},\label{eqn:stability-uh1}\\
				\sup\limits_{0\leq t\leq T}\exc{\normdel{\bv}}+\sup\limits_{0\leq t\leq T}\exc{\normbf{\opL\bu}}\leq& C_{s2},\label{eqn:stability-vh1}\\
				\sup\limits_{0\leq t\leq T}\exc{\normbf{\partial_{x_j}\bv}}+\sup\limits_{0\leq t\leq T}\exc{\normdel{\partial_{x_j}\bu}}\leq& C_{s3}\label{20220624_7}
			\end{align}
			\end{subequations}
			for $1\leq j\leq d$, and
			\begin{align*}
				C_{s1}&=\Big(\exc{\normbf{\bv_0}}+\exc{\normdel{\bu_0}}+4C_A^2\Big)e^{CC_B^2},\\
				C_{s2}&=\Big(\exc{\normdel{\bv_0}}+\exc{\normbf{\opL\bu_0}}+CC_A^2C_{s1}\Big)e^T,\\
				C_{s3}&=\Big(\exc{\normbf{\partial_{x_j}\bv_0}}+\exc{\normdel{\partial_{x_j}\bu_0}}\Big)e^{CC_A^2}.
			\end{align*}
		\end{lemma}
		
		Next, we state the following high moment stability result for $(\bu,\bv)$.  Since its proof is similar to \cite[Lemma 3.2]{feng2022higher} with slight modification, we omit it to save space.
		\begin{lemma}\label{lem:highmomentstability}
			Let $(\bu, \bv)$ be a strong solution to equations \eqref{eqn:elastic-wave-equation-uv}--\eqref{20220526_10}. Assume that $(\bu_0,\bv_0) \in L^p(\Omega; (\bH^2\cap\bH_0^1)\times \bH_0^1)$. For any integer $1 \leq p <\infty$, under the assumptions \eqref{assump:F(0)}--\eqref{assump:LipsF},
			there holds
			\begin{align}
				\mE\left[\sup_{0\leq t \leq T}\left(\|\bu(t)\|^{2p}_{\bH^2} + \|\bv(t)\|^{2p}_{\bH^1}\right)\right] \leq \tilde{K}_{1,p}
			\end{align}
			where the constant $\tilde{K}_{1,p} =\tilde{K}(\bu_0,\bv_0, p,C_A,C_B)$.
		\end{lemma}

		The next lemmas establish some time H\"older continuity results for the solution $(\bu,\bv)$ with respect to the $\bf{L}^2$-norm, $\bf{H}^1$-seminorm and $\bf{H}^2$-seminorm. Their proofs can be found in \cite[Lemma 2.6, Lemma 2.7]{feng2022fully}.
		
		\begin{lemma}\label{lem:holder-cont-u}
			{\itshape
				Let $(\bu, \bv)$ be a strong solution to problem \eqref{eqn:elastic-wave-equation-uv}--\eqref{20220526_10}. For any integer $1\leq p <\infty$, under the assumptions \eqref{assump:F(0)}--\eqref{assump:LipsF}, for any $s, t\in[0, T]$, we have
				\begin{subequations}
				\begin{align}
					\exc{\normbfp{\bu(t)-\bu(s)}}\leq& C_{s1}|t-s|^{2p},\label{eqn:holder-cont-u-1}\\
					\exc{\normdelpbig{\bu(t)-\bu(s)}}\leq& C_{s2}|t-s|^{2p},\label{eqn:holder-cont-u-2}\\
					\exc{\normbf{\opL\big(\bu(t)-\bu(s)\big)}}\leq& C_{s4}|t-s|^2.\label{eqn:holder-cont-u-3}
				\end{align}	
				\end{subequations}
			}
		\end{lemma}
		
		\begin{lemma}\label{lem:holder-cont-v}
			{\itshape
				Let $(\bu, \bv)$ be a strong solution to problem \eqref{eqn:elastic-wave-equation-uv}--\eqref{20220526_10}. Under the assumptions \eqref{assump:F(0)}--\eqref{assump:LipsF}, for any $s, t\in[0, T]$, we have
				\begin{subequations}
				\begin{align}
					\exc{\normbf{\bv(t)-\bv(s)}}\leq& C_{s5}|t-s|,\label{eqn:holder-cont-v-1}\\
					\exc{\normdelbig{\bv(t)-\bv(s)}}\leq& C_{s6}|t-s|,\label{eqn:holder-cont-v-2}\\
					\exc{\normbf{\opL\big(\bv(t)-\bv(s)\big)}}\leq& C_{s7}|t-s|,\label{eqn:holder-cont-v-3}
				\end{align}
				\end{subequations}
				where
				\begin{align*}
					C_{s5}&=e^T(CC_B^2C_{s1}+4C_A^2),\\
					C_{s6}&=e^T(CC_A^2C_{s1}+CC_B^2C_{s1}+2C_A^2),\\
					C_{s7}&=CC_A^2(C_{s1}+C_{s2}+C_{s4}).
				\end{align*}
			}
		\end{lemma}
		
		\section{Time-stepping scheme}\label{sec-3} 
		
		\subsection{Formulation of the time-stepping scheme}
		In this section, we consider the time discretization of  \eqref{eqn:weakform-1}--\eqref{eqn:weakform-2} in the following algorithm. A related time-stepping scheme considered in \cite{feng2024optimal} for the stochastic acoustic wave equations.
		\smallskip
		
		{\bf Algorithm 1.}\, Let $N>>1$ be a positive integer and $\tau:=T/N$. Let $\{t_n\}_{n=0}^N$ be a uniform partition of  the interval $[0,T]$ with mesh size $ \tau$. Find $\mathcal{F}_{t_n}$ adapted process $\{(\bu^n,\bv^n)\}_{n=0}^N$ such that there hold $\mathbb{P}$-almost surely
		\begin{subequations} \label{theta_scheme}
			\begin{align}
				(\bu^{n+1}-\bu^n,\bphi) &= \tau (\bv^{n+1},\bphi) -\bigl(\bG[\bu^n]\widehat{\Delta}W_n,\bphi \bigr) \qquad\qquad\quad  \forall \bphi \in \bH_0^1(D),\label{eq20230703_30}\\
				(\bv^{n+1} - \bv^n, \bpsi) &+ \tau\lambda \bigl(\Div \bu^{n,\frac12}, \Div \bpsi \bigr) + \tau\mu \bigl(\veps(\bu^{n,\frac12}),\veps(\bpsi) \bigr) \label{eq20230703_31} \\
				&= \tau\bigl(\bF^{n,\frac12}, \bpsi\bigr) 
				+ \bigl(\bG[\bu^n]\overline{\Delta} W_n, \bpsi \bigr)  \nonumber\\
				&\quad +  \bigl(D_{\bu}\bG[\bu^n]\bv^n\widehat{\Delta}W_n, \bpsi \bigr) \qquad\qquad\qquad\quad  \forall \, \bpsi \in \bH^1_0(D), \nonumber
			\end{align}
		\end{subequations}
		where 
		\begin{align*}
			\bu^{n,\frac12} :=  \frac12(\bu^{n+1} + \bu^{n-1}),\qquad  \bF^{n,\frac12} := \frac12( \bF[u^{n+1}] +  \bF[u^{n-1}]).
		\end{align*}

		\begin{remark}
			The equation \eqref{eq20230703_31} in Algorithm 1 iterates from $n = 1$ instead of $n=0$. Therefore, $\bu^1$ is also an initial value and must be specified.  Following \cite{feng2024optimal}  we choose $\bu^1$ as follows:
			\begin{align}\label{eq3.3}
				\bu^1 = \tau \bv^0 + \bu^0 -\frac{\tau^2}{2}\opL \bu^0 - \frac{\tau^2}{2}\bF[\bu^0] - \bG[\bu^0]\widehat{\Delta} W_0 + \tau\bG[\bu^0]\overline{\Delta} W_0,
			\end{align}
			where $(\bu^0,\bv^0) = (\bu_0,\bv_0)$.  With this choice of $\bu^1$,  it is easy to verify that
			\begin{align}\label{eq3.4}
				\mE\bigl[\|\opL(\bu(t_1) - \bu^1)\|^2_{\bL^2}\bigr] \leq C\tau,
			\end{align}
			where $C>0$ is a constant. This estimate will be used in the proof of Theorem \ref{theorem3.3} later.
		\end{remark}

		Next, we state and prove a crucial stability estimate for $\{(\bu^n, \bv^n)\}$ generated by Algorithm 1 in the next lemma. To state the result, we introduce 	the forward difference operator $\delta_t \bu^n:= (\bu^{n+1} - \bu^n)/\tau$ and define the following energy functional:
		\begin{align*}
			\mathcal{J}(\bphi,\bpsi) = \|\bpsi\|^2_{L^2}+\frac{\lambda}{2}\|\Div \bphi\|^2_{L^2} + \frac{\mu}{2}\|\varepsilon(\bphi)\|^2_{L^2}.
		\end{align*}
		
		\begin{lemma}\label{lemma3.1}
			Let $\{(\bu^n, \bv^n)\}$ denote the numerical solutions generated by Algorithm 1. Suppose that $(\bu_{0},\bv_0) \in L^2(\Omega; \bH_0^1 \times \bL^2)$. Under the assumptions \eqref{assump:F(0)}--\eqref{assump:Fuu}, there exists a constant $C_1>0$ such that
			\begin{align*}
				\max_{1 \leq n \leq N}\mE\left[\mathcal{J}(\bu^n,\bv^n)\right] \leq C_{1},
			\end{align*}
			where $C_1=Ce^{C(1+C_A)}\mE\bigl[\mathcal{J}(\bu^1,\bv^1) + \mathcal{J}(\bu^0,\bv^0)\bigr]$.
		\end{lemma}
		\begin{proof}\label{prf:time-semiform-stab}
			First, we rewrite \eqref{eq20230703_30} as
			\begin{align}\label{eq3.2}
				\Bigl(\frac{\bu^{n+1} - \bu^{n-1}}{\tau},\bphi\Bigr) &=\Bigl( \frac{\bu^{n+1} - \bu^n}{\tau} + \frac{\bu^{n} - \bu^{n-1}}{\tau},\bphi\Bigr)\\\nonumber
				&= \bigl(\bv^{n+1} + \bv^n,\bphi\bigr) -\frac{1}{\tau} \bigl(\bG[\bu^n]\widehat{\Delta}W_n, \bphi\bigr) \\\nonumber
				&\qquad- \frac{1}{\tau}\bigl(\bG[\bu^{n-1}]\widehat{\Delta}W_{n-1},\bphi\bigr).
			\end{align}
			Taking $\bpsi = (\bu^{n+1} - \bu^{n-1})/\tau$ in \eqref{eq20230703_31} and combining with \eqref{eq3.2}, we obtain
			\begin{align}\label{eq3.11}
				&\bigl(\bv^{n+1} - \bv^n, \bv^{n+1} + \bv^n\bigr)	+ \frac{\lambda}{2} \bigl(\Div(\bu^{n+1} + \bu^{n-1}), \Div(\bu^{n+1}  - \bu^{n-1})\bigr) \\\nonumber
				&\qquad\qquad+ \frac{\mu}{2} \bigl(\varepsilon(\bu^{n+1} + \bu^{n-1}), \varepsilon(\bu^{n+1}  - \bu^{n-1})\bigr) \\\nonumber
				&= \frac{1}{\tau}\bigl(\bG[\bu^n]\widehat{\Delta}W_n + \bG[\bu^{n-1}]\widehat{\Delta}W_{n-1}, \bv^{n+1} - \bv^n\bigr) \\\nonumber
				&\qquad+ \Bigl(\bG[\bu^n]\overline{\Delta}W_n, \frac{\bu^{n+1}- \bu^{n-1}}{\tau}  \Bigr) \\\nonumber
				&\qquad+ \Bigl(D_{\bu}\bG[\bu^n]\bv^n \widehat{\Delta}W_n , \frac{\bu^{n+1} - \bu^{n-1}}{\tau}\Bigr)\\\nonumber
				&\qquad + \frac{\tau}{2}\Bigl(\bF[\bu^{n+1}] + \bF[\bu^{n-1}], \frac{\bu^{n+1} - \bu^{n-1}}{\tau}\Bigr)\\\nonumber
				&:= \mathcal{I}_1 + \mathcal{I}_2 + \mathcal{I}_3 + \mathcal{I}_4.
			\end{align}
			
			The left-hand side of \eqref{eq3.11} can be written as
			\begin{align*}
				&\bigl(\bv^{n+1} - \bv^n, \bv^{n+1} + \bv^n\bigr)	+ \frac{\lambda}{2} \bigl(\Div(\bu^{n+1} + \bu^{n-1}), \Div(\bu^{n+1}  - \bu^{n-1})\bigr) \\\nonumber
				&\qquad+ \frac{\mu}{2} \bigl(\varepsilon(\bu^{n+1} + \bu^{n-1}), \varepsilon(\bu^{n+1}  - \bu^{n-1})\bigr) \\\nonumber
				&=\bigl[\|\bv^{n+1}\|^2_{\bL^2} - \|\bv^n\|^2_{\bL^2}\bigr] + \frac{\lambda}{2} \bigl[\|\Div \bu^{n+1}\|^2_{\bL^2} - \|\Div \bu^{n-1}\|^2_{\bL^2}\bigr] \\\nonumber
				&\qquad+ \frac{\mu}{2} \bigl[\|\varepsilon(\bu^{n+1})\|^2_{\bL^2} - \|\varepsilon(\bu^{n-1})\|^2_{\bL^2}\bigr].
			\end{align*}
			
			Next,  we bound the right-hand side of \eqref{eq3.11} as follows.  First,  choosing $\bpsi = \bG[\bu^n]\widehat{\Delta}W_n + \bG[\bu^{n-1}]\widehat{\Delta}W_{n-1}$ in \eqref{eq20230703_31} and then substituting it into $\mathcal{I}_1$ we get
			\begin{align*}
				\mE[\mathcal{I}_1] &= \frac{1}{\tau}\mE\bigl[\bigl(\bG[\bu^n]\widehat{\Delta}W_n + \bG[\bu^{n-1}]\widehat{\Delta}W_{n-1}, \bv^{n+1} -\bv^n\bigr)\bigr] \\\nonumber
				& = -\frac{\lambda}{2}\mE\Bigl[\bigl(\Div(\bu^{n+1} + \bu^{n-1}), \Div\bigl(\bG[\bu^n]\widehat{\Delta}W_n + \bG[\bu^{n-1}]\widehat{\Delta}W_{n-1}\bigr)\bigr)\Bigr]\\\nonumber
				&\qquad -\frac{\mu}{2}\mE\Bigl[\bigl(\varepsilon(\bu^{n+1} + \bu^{n-1}), \varepsilon\bigl(\bG[\bu^n]\widehat{\Delta}W_n + \bG[\bu^{n-1}]\widehat{\Delta}W_{n-1}\bigr)\bigr)\Bigr]\\\nonumber
				&\qquad+ \frac12\mE\Bigl[\bigl(\bF[\bu^{n+1}] + \bF[\bu^{n-1}], \bG[\bu^n]\widehat{\Delta}W_n + \bG[\bu^{n-1}]\widehat{\Delta}W_{n-1}\bigr)\Bigr]\\\nonumber
				&\qquad +\frac{1}{\tau}\mE\Bigl[\bigl(\bG[\bu^n]\overline{\Delta}W_{n}, \bG[\bu^n]\widehat{\Delta}W_n + \bG(\bu^{n-1})\widehat{\Delta}W_{n-1}\bigr)\Bigr] \\\nonumber
				&\qquad+\frac{1}{\tau}\mE\Bigl[\bigl(D_{\bu}\bG[\bu^n]\bv^n\widehat{\Delta}W_n, \bG[\bu^n]\widehat{\Delta}W_n + \bG[\bu^{n-1}]\widehat{\Delta}W_{n-1}\bigr) \Bigr]\\\nonumber
				&:= \mathcal{I}_{1,1} +  \mathcal{I}_{1,2} +  \mathcal{I}_{1,3} +  \mathcal{I}_{1,4} + \mathcal{I}_{1,5}.   
			\end{align*}
			
			It follows from Lemma $\ref{increment}(\rm ii)_2$, \eqref{assump:GradientF} and the Korn's inequality that
			\begin{align*}
				\mathcal{I}_{1,1} + \mathcal{I}_{1,2}&\leq C\tau \lambda\mE\bigl[\|\Div \bu^{n+1}\|^2_{\bL^2} + \|\Div \bu^{n-1}\|^2_{\bL^2}\bigr]\\\nonumber
				&\qquad + \frac{1}{\tau} \lambda\mE\bigl[\|\Div\bigl(\bG[\bu^n]\widehat{\Delta}W_n + \bG[\bu^{n-1}]\widehat{\Delta}W_{n-1}\bigr)\|^2_{\bL^2}\bigr]\\\nonumber
				&\qquad+C\tau \mu\mE\bigl[\|\varepsilon(\bu^{n+1})\|^2_{\bL^2} + \|\varepsilon( \bu^{n-1})\|^2_{\bL^2}\bigr]\\\nonumber
				&\qquad + \frac{1}{\tau} \mu\mE\bigl[\|\varepsilon\bigl(\bG[\bu^n]\widehat{\Delta}W_n + \bG[\bu^{n-1}]\widehat{\Delta}W_{n-1}\bigr)\|^2_{\bL^2}\bigr]\\\nonumber
				&\leq C\tau \lambda \mE\bigl[\|\Div \bu^{n+1}\|^2_{\bL^2} + \|\Div \bu^{n-1}\|^2_{\bL^2}\bigr] + CC_A\tau^2 \lambda \mE\bigl[\|\Div \bu^{n}\|^2_{\bL^2} \\ \nonumber
				&\qquad + \|\Div \bu^{n-1}\|^2_{\bL^2}\bigr] 
				+C\tau \mu \mE\bigl[\|\varepsilon(\bu^{n+1})\|^2_{\bL^2} + \|\varepsilon( \bu^{n-1})\|^2_{\bL^2}\bigr] \\ \nonumber
				&\qquad + CC_A\tau^2 \mu \mE\bigl[\|\varepsilon(\bu^{n})\|^2_{\bL^2} + \|\varepsilon (\bu^{n-1})\|^2_{\bL^2}\bigr]\\\nonumber
				&\leq CC_A\tau \lambda \mE\bigl[\|\Div \bu^{n+1}\|^2_{\bL^2} +\|\Div \bu^n\|^2_{\bL^2}+ \|\Div \bu^{n-1}\|^2_{\bL^2}\bigr]\\\nonumber
				&\qquad+CC_A\tau \mu \mE\bigl[\|\varepsilon(\bu^{n+1})\|^2_{\bL^2} +\|\varepsilon (\bu^n)\|^2_{\bL^2}+ \|\varepsilon(\bu^{n-1})\|^2_{\bL^2}\bigr].
			\end{align*}

			Using Lemma $\ref{increment}(\rm ii)_2$, \eqref{assump:F(0)}, \eqref{assump:GradientF}, Poincare's inequality and Korn's inequality we obtain
			\begin{align*}
				\mathcal{I}_{1,3} &\leq CC_A\tau \mE\left[\lambda\|\Div \bu^{n+1}\|^2_{\bL^2} +\mu \|\varepsilon(\bu^{n+1})\|^2_{\bL^2} + \|\bu^{n+1}\|^2_{\bL^2}\right] \\\nonumber
				&\qquad+ CC_A\tau \mE\left[\lambda\|\Div \bu^{n}\|^2_{\bL^2} +\mu \|\varepsilon(\bu^{n})\|^2_{\bL^2} + \|\bu^{n}\|^2_{\bL^2}\right] \\\nonumber
				&\qquad+ CC_A\tau \mE\left[\lambda\|\Div \bu^{n-1}\|^2_{\bL^2} +\mu \|\varepsilon(\bu^{n-1})\|^2_{\bL^2} + \|\bu^{n-1}\|^2_{\bL^2}\right] + CC_A\tau\\\nonumber
				&\leq CC_A\tau \mE\left[\lambda\|\Div \bu^{n+1}\|^2_{\bL^2} +\mu \|\varepsilon(\bu^{n+1})\|^2_{\bL^2}\right] \\\nonumber
				&\qquad+ CC_A\tau \mE\left[\lambda\|\Div \bu^{n}\|^2_{\bL^2} +\mu \|\varepsilon(\bu^{n})\|^2_{\bL^2}\right] \\\nonumber
				&\qquad+ CC_A\tau \mE\left[\lambda\|\Div \bu^{n-1}\|^2_{\bL^2} +\mu \|\varepsilon(\bu^{n-1})\|^2_{\bL^2} \right] + CC_A\tau.
			\end{align*}
			
			Similarly to $\mathcal{I}_{1,3}$, we also can show
			\begin{align*}
				\mathcal{I}_{1,4} &\leq C\mE\left[\|\bG[\bu^n]\overline{\Delta} W_n\|^2_{\bL^2}\right]  + \frac{C}{\tau^2}\mE\left[\|\bG[\bu^n]\widehat{\Delta}W_n + \bG[\bu^{n-1}]\widehat{\Delta}W_{n-1}\|^2_{\bL^2}\right]\\\nonumber
				&\leq  CC_A\tau \mE\left[\lambda\|\Div \bu^{n}\|^2_{\bL^2} +\mu \|\varepsilon(\bu^{n})\|^2_{\bL^2}\right] \\\nonumber
				&\qquad+ CC_A\tau \mE\left[\lambda\|\Div \bu^{n-1}\|^2_{\bL^2} +\mu \|\varepsilon(\bu^{n-1})\|^2_{\bL^2}\right] + CC_A\tau.
			\end{align*}	
			
			Next, using \eqref{assump:GradientF}, Remark $\ref{increment}(\rm ii)_2$,  \eqref{assump:Fuu}, Poincare's inequality and Korn's inequality we obtain
			\begin{align*}
				\mathcal{I}_{1,5} &\leq CC_A\tau \mE\bigl[\|\bv^n\|^2_{\bL^2} +  CC_B\tau \mE\left[\lambda\|\Div \bu^{n}\|^2_{\bL^2} +\mu \|\varepsilon(\bu^{n})\|^2_{\bL^2} + \|\bu^{n}\|^2_{\bL^2}\right] \\\nonumber
				&\qquad+ CC_A\tau \mE\left[\lambda\|\Div \bu^{n-1}\|^2_{\bL^2} +\mu \|\varepsilon(\bu^{n-1})\|^2_{\bL^2} + \|\bu^{n-1}\|^2_{\bL^2}\right] + CC_A\tau\\\nonumber
				&\leq CC_A\tau \mE\bigl[\|\bv^n\|^2_{\bL^2} +  CC_A\tau \mE\left[\lambda\|\Div \bu^{n}\|^2_{\bL^2} +\mu \|\varepsilon(\bu^{n})\|^2_{\bL^2} \right] \\\nonumber
				&\qquad+ CC_A\tau \mE\left[\lambda\|\Div \bu^{n-1}\|^2_{\bL^2} +\mu \|\varepsilon(\bu^{n-1})\|^2_{\bL^2} \right] + CC_A\tau.
			\end{align*}
			In summary, we have 
			\begin{align*}
				\mE\bigl[\mathcal{I}_1\bigr] &\leq CC_A\tau \mE\left[\mathcal{J}(\bu^{n+1},\bv^{n+1})\right] +CC_A\tau \mE\left[\mathcal{J}(\bu^{n},\bv^{n})\right] \\\nonumber
				&\qquad+ CC_A\tau \mE\left[\mathcal{J}(\bu^{n-1},\bv^{n-1})\right]+ CC_A\tau.
			\end{align*}
			
			Now, we continue to bound $\mathcal{I}_2$.  By using the martingale property of the increments $\overline{\Delta}W_n$ and \eqref{eq20230703_31} with $\bpsi = \bG[\bu^n]\overline{\Delta}W_{n}$, we get 
			\begin{align*}
				\mE\bigl[\mathcal{I}_2\bigr] &= \frac{1}{\tau}\mE\bigl[\bigl(\bG[\bu^n]\overline{\Delta}W_n, \bu^{n+1} - \bu^{n-1}\bigr)\bigr] \\\nonumber
				&= \mE\bigl[\bigl(\bG[\bu^n]\overline{\Delta}W_n, \bv^{n+1} + \bv^n\bigr)\bigr]\\\nonumber
				&= \mE\bigl[\bigl(\bG[\bu^n]\overline{\Delta}W_n, \bv^{n+1}\bigr)\bigr]\\\nonumber
				&=  \mE\bigl[\bigl(\bG[\bu^n]\overline{\Delta}W_n, \bv^{n+1} - \bv^n\bigr)\bigr]\\\nonumber
				&=  -\frac{\tau}{2}\lambda \mE\bigl[\bigl(\Div\bG[\bu^n] \overline{\Delta}W_n, \Div(\bu^{n+1} + \bu^{n-1})\bigr)\bigr] \\\nonumber
				&\qquad-\frac{\tau}{2}\mu \mE\bigl[\bigl(\varepsilon(\bG[\bu^n]\overline{\Delta}W_n), \varepsilon(\bu^{n+1} + \bu^{n-1})\bigr)\bigr] \\\nonumber
				&
				\qquad+  \frac{\tau}{2}\mE\bigl[\bigl(\bG[\bu^n]\overline{\Delta}W_n, \bF[\bu^{n+1}] + \bF[\bu^{n-1}]\bigr)\bigr] + \mE\bigl[\|\bG[\bu^n]\overline{\Delta}W_n\|^2_{L^2}\bigr]\\\nonumber
				&\qquad+  \mE\bigl[\bigl(\bG[\bu^n]\overline{\Delta}W_n, D_u\bG[\bu^n]\bv^n\widehat{\Delta}W_n\bigr)\bigr]\\\nonumber
				&=  \mathcal{I}_{2,1}  +  \mathcal{I}_{2,2} +  \mathcal{I}_{2,3} +  \mathcal{I}_{2,4}.   
			\end{align*}
			Each of  $\mathcal{I}_{2,1},  \mathcal{I}_{2,2}, \mathcal{I}_{2,3},$ and $\mathcal{I}_{2,4}$ can be controlled using the same techniques for bounding  $ \mathcal{I}_{1,1},  \mathcal{I}_{1,2}, \mathcal{I}_{1,3}, $ and $\mathcal{I}_{1,4}$ above (we skip the repetitions to save space). In summary, we have
			\begin{align*}
				\mE\bigl[\mathcal{I}_2\bigr] &\leq CC_A\tau \mE\left[\mathcal{J}(\bu^{n+1},\bv^{n+1})\right] +CC_A\tau \mE\left[\mathcal{J}(\bu^{n},\bv^{n})\right] \\\nonumber
				&\qquad+ CC_A\tau \mE\left[\mathcal{J}(\bu^{n-1},\bv^{n-1})\right]+ CC_A\tau.
			\end{align*}
			
			To bound $\mathcal{I}_3$, using the martingale property of the increments, \eqref{assump:GradientF}, \eqref{assump:Fuu}, and \eqref{eq20230703_30} with $\bphi = D_u\bG[\bu^n]\bv^n\widehat{\Delta}W_n$, we obtain
			\begin{align*}
				\mE\bigl[\mathcal{I}_3\bigr] &= \frac{1}{\tau}\mE\bigl[\bigl(D_{\bu}\bG[\bu^n]\bv^n\widehat{\Delta}W_n, \bu^{n+1} - \bu^n\bigr)\bigr]\\\nonumber
				&= \mE\bigl[\bigl(D_{\bu}\bG[\bu^n]\bv^n\widehat{\Delta}W_n, \bv^{n+1}\bigr)\bigr] - \frac{1}{\tau}\mE\bigl[\bigl(D_{\bu}\bG[\bu^n]\bv^n\widehat{\Delta}W_n, \bG[\bu^n]\widehat{\Delta}W_n\bigr)\bigr]\\\nonumber
				&\leq C\tau\mE\bigl[\|\bv^{n+1}\|^2_{\bL^2}\bigr] + CC_A\tau\mE\bigl[\|\bv^n\|^2_{\bL^2}\bigr] \\\nonumber
				&\qquad+ CC_A\tau\mE\bigl[\| \bv^n\|^2_{\bL^2}\bigr] + C\mE\bigl[\|\bG[\bu^n]\widehat{\Delta}W_n\|^2_{\bL^2}\bigr]\\\nonumber
				&\leq C\tau\mE\bigl[\mathcal{J}(\bu^{n+1},\bv^{n+1})\bigr] + CC_A\tau\mE\bigl[\mathcal{J}(\bu^n,\bv^n)\bigr] + CC_A\tau.
			\end{align*}
			
			Finally, we bound $\mathcal{I}_4$ as follows.  Using \eqref{eq20230703_30} with test function  $\bphi = \bF[\bu^{n+1}] + \bF[\bu^{n-1}]$ and applying \eqref{assump:F(0)}, \eqref{assump:Fuu} yield
			\begin{align}
				\mE\bigl[\mathcal{I}_4\bigr] &= \frac{\tau}{2}\mE\bigl[\bigl(\bF[\bu^{n+1}] + \bF[\bu^{n-1}], \bv^{n+1} + \bv^n\bigr)\bigr] \\\nonumber
				&\qquad- \frac12\mE\bigl[\bigl(\bF[\bu^{n+1}] + \bF[\bu^{n-1}], \bG[\bu^n]\widehat{\Delta}W_n\bigr)\bigr] \\\nonumber
				&\qquad- \frac12\mE\bigl[\bigl(\bF[\bu^{n+1}] + \bF[\bu^{n-1}], \bG[\bu^{n-1}]\widehat{\Delta}W_{n-1}\bigr)\bigr]\\\nonumber
				&\leq C\tau\mE\bigl[\|\bv^{n+1}\|^2_{\bL^2} + \|\bv^n\|^2_{\bL^2}\bigr] + C\tau\mE\bigl[\|\bF[\bu^{n+1}]\|^2_{\bL^2} + \|\bF[\bu^{n-1}]\|^2_{\bL^2}\bigr]\\\nonumber
				&\qquad+ \frac{C}{\tau}\mE\bigl[\|\bG[\bu^n]\widehat{\Delta}W_n\|^2_{\bL^2} + \|\bG[\bu^{n-1}]\widehat{\Delta}W_{n-1}\|^2_{\bL^2}\bigr] \\\nonumber
				&\leq CC_A\mE\bigl[\mathcal{J}(\bu^{n+1},\bv^{n+1}) + \mathcal{J}(\bu^n,\bv^n) + \mathcal{J}(\bu^{n-1},\bv^{n-1})\bigr] + CC_A\tau.
			\end{align} 
			
			Collecting all the estimates for $\mathcal{I}_1, \mathcal{I}_2, \mathcal{I}_3, \mathcal{I}_4$ and substituting them into \eqref{eq3.11} followed  by taking summation from $n = 1$ to $n = \ell$ for $1 \leq \ell <N$, we get
			\begin{align*}
				&	\mE\bigl[\|\bv^{\ell +1}\|^2_{\bL^2}\bigr] + \frac{\lambda}{2} \mE\bigl[\|\Div\bu^{\ell+1}\|^2_{\bL^2}\bigr] + \frac{\mu}{2} \mE\bigl[\|\varepsilon(\bu^{\ell+1})\|^2_{\bL^2}\bigr]   \\\nonumber
				&\leq C\mE\bigl[\|\bv^1\|^2_{\bL^2}\bigr] + \frac{\lambda}{2}\mE\bigl[\|\Div\bu^1\|^2_{\bL^2} + \|\Div\bu^0\|^2_{\bL^2}\bigr] + \frac{\mu}{2}\mE\bigl[\|\varepsilon(\bu^1)\|^2_{\bL^2} + \|\varepsilon(\bu^0)\|^2_{\bL^2}\bigr]  \\\nonumber
				&\,\,+C(1+C_A)\tau\sum_{n=1}^{\ell}\mE\bigl[\mathcal{J}(\bu^{n+1},\bv^{n+1})+\mathcal{J}(\bu^{n},\bv^{n})+\mathcal{J}(\bu^{n-1},\bv^{n-1})\bigr] + CC_A,
			\end{align*}
			which implies that
			\begin{align}\label{eq3.13}
				\mE\bigl[\mathcal{J}(\bu^{\ell+1},\bv^{\ell+1})\bigr]
				&\leq C\mE\bigl[\mathcal{J}(\bu^1,\bv^1) + \mathcal{J}(\bu^0,\bv^0)\bigr] \\\nonumber
				&\qquad+C(1+C_A)\tau\sum_{n=1}^{\ell}\mE\bigl[\mathcal{J}(\bu^{n+1},\bv^{n+1})\\\nonumber
				&\qquad+\mathcal{J}(\bu^{n},\bv^{n})+\mathcal{J}(\bu^{n-1},\bv^{n-1})\bigr] + CC_A.
			\end{align}
			
			The proof is complete by applying the discrete Gronwall's inequality to \eqref{eq3.13}.
		\end{proof}
		
	 To state the next lemma concerning high moment stability estimates, we introduce the following discrete energy functional:
		\begin{align}\label{Q}
			\mathcal{Q}(\bu^{n+1},\overline{\bu}^{n+1}):= \frac12 \|\bu^{n+1}\|^2_{\bL^2} + \frac{\tau^2}{4}\lambda \|\Div \overline{\bu}^{n+1}\|^2_{\bL^2} + \frac{\tau^2}{2}\mu\|\varepsilon(\overline{\bu}^{n+1})\|^2_{\bL^2}.
		\end{align}
		
		
		\begin{lemma}\label{lemma3.1.2}
			Let $\{(\bu^n, \bv^n)\}$ denote the numerical solutions generated by Algorithm 1. Suppose that $(\bu_0,\bv_0) \in L^p(\Omega; \bH_0^1\times \bL^2)$. For any $1\leq p <\infty$, under the assumptions \eqref{assump:F(0)}--\eqref{assump:Fuu}, there exists a constant $C_{2,p}>0$ such that
			\begin{align*}
				\max_{1 \leq n \leq N}	\mE\left[(\mathcal{Q}(\bu^n,\overline{\bu}^{n}))^p\right] \leq C_{2,p},
			\end{align*}
			where 
			\begin{align*}
				C_{2,p} &= \biggl\{2^{2p}\mE\Bigl[\Bigl\{\mathcal{Q}(\bu^1,\overline{\bu}^1) + \frac{\tau^2}{4}\lambda\|\Div\overline{\bu}^{0}\|^2_{\bL^2} + \frac{\tau^2}{4}\mu\|\varepsilon(\overline{\bu}^{0})\|^2_{\bL^2}\Bigr\}^p\Bigr]\\\nonumber
				&\quad + C_p \mE\left[\|\bv_{0} - \frac{\tau}{2}\mathcal{L}\bu^0 + \bG[\bu^0]\overline{\Delta}W_0 + \tau\mathcal{L}\bu^{\frac12}\|^{2p}_{\bL^2}\right] \\\nonumber 
				&\quad+C_pC_A^p T^2 +C_pC_A^{4p}\tau^{3p}T^2+C_pC_A^{4p}\tau^{4p}T^2\biggr\}e^{C_p+C_pC_A^p(2C_A^{p}+1 + 2C_A^{3p})T}.
			\end{align*}
		\end{lemma}
		
		\begin{proof}
		Substituting \eqref{eq20230703_30} into \eqref{eq20230703_31} we obtain
			\begin{align}\label{eq3.233}
				&\bigl((\bu^{m+1} - \bu^m) - (\bu^m - \bu^{m-1}), \bphi\bigr) + \bigl(\bG[\bu^m] \widehat{\Delta} W_m - \bG[\bu^{m-1}]\widehat{\Delta} W_{m-1}, \bphi\bigr)\\\nonumber
				&\qquad \qquad +\tau^2 \lambda\bigl(\Div \bu^{m,\frac12},\Div \bphi\bigr) +\tau^2\mu \bigl(\veps (\bu^{m,\frac12}),\veps(\bphi)\bigr) \\\nonumber
				&\quad = \tau^2 \bigl(\bF^{m,\frac12},\bphi\bigr) + \tau \bigl(\bG[\bu^m]\overline{\Delta} W_m, \bphi\bigr) + \tau\bigl(D_{\bu}\bG[\bu^m]\bv^m\widehat{\Delta}W_m, \bphi\bigr).
			\end{align}
			Applying the summation operator $\sum_{m=1}^n$ ($1\leq n <N$) to \eqref{eq3.233} and denoting $\overline{\bu}^{n} := \sum_{m=1}^n \bu^{m}$ and $\bu^{\frac12} = \bu^{0+\frac12} = \frac12\bigl(\bu^1 + \bu^0\bigr)$, we get 
			\begin{align}\label{eq3.244}
				&	\bigl(\bu^{n+1} - \bu^n, \bphi\bigr)  + \tau^2\lambda\bigl(\Div\overline{\bu}^{n,\frac12},\Div\bphi\bigr) + \tau^2\mu\bigl(\veps(\overline{\bu}^{n,\frac12}),\veps(\bphi)\bigr)\\\nonumber
				&= \bigl(\bu^1 - \bu^0,\bphi\bigr) + \bigl(\bG[\bu^0]\widehat{\Delta} W_0 + \tau^2\opL \bu^{\frac12},\bphi\bigr)   \\\nonumber
				&\qquad- \bigl(\bG[\bu^n]\widehat{\Delta} W_n, \bphi\bigr)+ \tau^2\sum_{m=1}^n \bigl(\bF^{m,\frac12}, \bphi\bigr)\\\nonumber 
				&\qquad + \tau\sum_{m=1}^n \bigl(\bG[\bu^m],\bphi\bigr)\overline{\Delta} W_m +\tau\sum_{m=1}^n\bigl(D_{\bu}\bG[\bu^m]\bv^m\widehat{\Delta}W_m, \bphi\bigr).
			\end{align}
			
			Without loss of generality, we proceed by assuming $\bF \equiv 0$ to avoid technicalities and save space.
			Taking $\bphi = \bu^{n+\frac12}$ in \eqref{eq3.244}, we obtain 
			\begin{align}\label{eq3.10}
				\frac12[\|\bu^{n+1}\|^2_{\bL^2} &- \|\bu^n\|^2_{\bL^2}] + \frac{\tau^2}{4}\lambda\left[\|\Div \overline{\bu}^{n+1}\|^2_{\bL^2} - \|\Div\overline{\bu}^{n-1}\|^2_{\bL^2}\right]\\\nonumber
				&+\frac{\tau^2}{4}\mu\left[\|\varepsilon(\overline{\bu}^{n+1})\|^2_{\bL^2} - \|\varepsilon(\overline{\bu}^{n-1})\|^2_{\bL^2}\right]\\\nonumber
				&= \left(\bu^1-\bu^0 + \bG[\bu^0]\widehat{\Delta}W_0 + \tau^2\mathcal{L}\bu^{\frac12},\bu^{n+\frac12}\right)\\\nonumber
				&\qquad- \bigl(\bG[\bu^n]\widehat{\Delta} W_n, \bu^{n+\frac12}\bigr)+ \tau\left(\sum_{m=1}^n\bG[\bu^m]\overline{\Delta}W_m, \bu^{n+\frac12}\right) \\\nonumber
				&\qquad + \tau\left(\sum_{m=1}^n D_{\bu}\bG[\bu^m]\bv^m\widehat{\Delta}W_m, \bu^{n+\frac12}\right).
			\end{align}
			Using the definition of $\mathcal{Q}$ and \eqref{eq3.10} we arrive at
			\begin{align}\label{eq3.111}
				\mathcal{Q}(\bu^{n+1},\overline{\bu}^{n+1}) - \mathcal{Q}(\bu^n,\overline{\bu}^n) &= \frac{\tau^2}{4}\lambda\left[\|\Div\overline{\bu}^{n-1}\|^2_{\bL^2} - \|\Div\overline{\bu}^{n}\|^2_{\bL^2}\right] \\\nonumber
				&\qquad+ \frac{\tau^2}{4}\mu\left[\|\varepsilon(\overline{\bu}^{n-1})\|^2_{\bL^2} - \|\varepsilon(\overline{\bu}^{n}\|^2_{\bL^2})\right]\\\nonumber
				&\qquad+\left(\bu^1-\bu^0 + \bG[\bu^0]\widehat{\Delta}W_0 + \tau^2\mathcal{L}\bu^{\frac12},\bu^{n+\frac12}\right)\\\nonumber
				&\qquad- \bigl(\bG[\bu^n]\widehat{\Delta} W_n, \bu^{n+\frac12}\bigr)\\\nonumber
				&\qquad+ \tau\left(\sum_{m=1}^n\bG[\bu^m]\overline{\Delta}W_m, \bu^{n+\frac12}\right) \\\nonumber
				&\qquad+ \tau\left(\sum_{m=1}^n D_{\bu}\bG[\bu^m]\bv^m\widehat{\Delta}W_m, \bu^{n+\frac12}\right).
			\end{align}
			Applying the summation operator $\sum_{n=1}^{\ell}$ for any $0<\ell\leq N-1$ to \eqref{eq3.111}, we get
			\begin{align}\label{eq3.12}
				\mathcal{Q}(\bu^{\ell+1},\overline{\bu}^{\ell+1}) &\leq \Bigl\{\mathcal{Q}(\bu^1,\overline{\bu}^1) + \frac{\tau^2}{4}\lambda\|\Div\overline{\bu}^{0}\|^2_{\bL^2} + \frac{\tau^2}{4}\mu\|\varepsilon(\overline{\bu}^{0})\|^2_{\bL^2} \Bigr\}\\\nonumber
				&\qquad+\left(\bu^1-\bu^0 + \bG[\bu^0]\widehat{\Delta}W_0 + \tau^2\mathcal{L}\bu^{\frac12},\sum_{n=1}^{\ell}\bu^{n+\frac12}\right)\\\nonumber
				&\qquad+ \tau\sum_{n=1}^{\ell}\left(\sum_{m=1}^n\bG[\bu^m]\overline{\Delta}W_m, \bu^{n+\frac12}\right) \\\nonumber
				&\qquad+ \tau\sum_{n=1}^{\ell}\left(\sum_{m=1}^n D_{\bu}\bG[\bu^m]\bv^m\widehat{\Delta}W_m, \bu^{n+\frac12}\right)\\\nonumber
				&\qquad- \sum_{n=1}^{\ell}\bigl(\bG[\bu^n]\widehat{\Delta} W_n, \bu^{n+\frac12}\bigr).
			\end{align}
			Taking the $p$th-power to \eqref{eq3.12} for any $p\geq 1$ and using the inequality $(a+b) \leq 2^p(a^p + b^p)$,  we obtain 
			\begin{align}\label{eq3.133}
				\mathcal{Q}^p(\bu^{\ell+1},\overline{\bu}^{\ell+1}) &\leq 2^{2p}\Bigl\{\mathcal{Q}(\bu^1,\overline{\bu}^1) + \frac{\tau^2}{4}\lambda\|\Div\overline{\bu}^{0}\|^2_{\bL^2} + \frac{\tau^2}{4}\mu\|\varepsilon(\overline{\bu}^{0})\|^2_{\bL^2}\Bigr\}^p\\\nonumber
				&\qquad+2^{2p}\left\{\left(\bu^1-\bu^0 + \bG[\bu^0]\widehat{\Delta}W_0 + \tau^2\mathcal{L}\bu^{\frac12},\sum_{n=1}^{\ell}\bu^{n+\frac12}\right)\right\}^p\\\nonumber
				&\qquad+ 2^{2p}\left\{\tau\sum_{n=1}^{\ell}\left(\sum_{m=1}^n\bG[\bu^m]\overline{\Delta}W_m, \bu^{n+\frac12}\right)\right\}^p \\\nonumber
				&\qquad+2^{3p}\left\{ \tau\sum_{n=1}^{\ell}\left(\sum_{m=1}^n D_{\bu}\bG[\bu^m]\bv^m\widehat{\Delta}W_m, \bu^{n+\frac12}\right)\right\}^p\\\nonumber
				&\qquad+2^{3p} \left\{\sum_{n=1}^{\ell}\|\bG[\bu^n]\widehat{\Delta} W_n\|_{\bL^2} \|\bu^{n+\frac12}\|_{\bL^2}\right\}^p\\\nonumber
				&:=K_1 + K_2 + K_3 + K_4+ K_5.
			\end{align}

			We notice that choosing the initial values $(\bu_0,\bv_0)$ and computing $\bu^1$ by using \eqref{eq3.3} help to control $\mE[K_1]$. Next, to bound $K_2$, we use \eqref{eq3.3} to obtain
			\begin{align*}
				\mE[K_2] &= 2^{2p}\mE\left[\left\{\left(\bu^1-\bu^0 + \bG[\bu^0]\widehat{\Delta}W_0 + \tau^2\mathcal{L}\bu^{\frac12},\sum_{n=1}^{\ell}\bu^{n+\frac12}\right)\right\}^p\right]\\\nonumber
				&=2^{2p}\mE\left[\left\{\left(\tau \bv^0 - \frac{\tau^2}{2}\mathcal{L}\bu^0 + \tau\bG[\bu^0]\overline{\Delta}W_0 + \tau^2\mathcal{L}\bu^{\frac12},\sum_{n=1}^{\ell}\bu^{n+\frac12}\right)\right\}^p\right]\\\nonumber
				&\leq C_p \mE\left[\|\bv_{0} - \frac{\tau}{2}\mathcal{L}\bu^0 + \bG[\bu^0]\overline{\Delta}W_0 + \tau\mathcal{L}\bu^{\frac12}\|^{2p}_{\bL^2}\right] + C_p\mE\left[\left\|\tau\sum_{n=1}^{\ell}\bu^{n+\frac12}\right\|^{2p}_{\bL^2}\right]\\\nonumber
				&\leq  C_p \mE\left[\|\bv_{0} - \frac{\tau}{2}\mathcal{L}\bu^0 + \bG[\bu^0]\overline{\Delta}W_0 + \tau\mathcal{L}\bu^{\frac12}\|^{2p}_{\bL^2}\right] + C_p\mE\left[\tau\sum_{n=1}^{\ell}\|\bu^{n+\frac12}\|^{2p}_{\bL^2}\right]\\\nonumber
				&\leq  C_p \mE\left[\|\bv_{0} - \frac{\tau}{2}\mathcal{L}\bu^0 + \bG[\bu^0]\overline{\Delta}W_0 + \tau\mathcal{L}\bu^{\frac12}\|^{2p}_{\bL^2}\right] \\\nonumber&\qquad+ C_p\mE\left[\tau\sum_{n=1}^{\ell}(\mathcal{Q}^p(\bu^{n+1},\overline{\bu}^{n+1}) + \mathcal{Q}^p(\bu^n,\overline{\bu}^n))\right],
			\end{align*}
			where the second inequality  is obtained by using the discrete H\"older's inequality, while the last inequality follows from the fact that  $\|\bu^{n+1}\|_{\bL^2}^2 \leq 2\mathcal{Q}(\bu^{n+1},\overline{\bu}^{n+1})$.
			
			Next, we estimate $K_3$ by using the discrete H\"older's inequality and then the Cauchy-Schwarz inequality as follows:
			\begin{align*}
				\mE[K_3] &\leq C_p \mE\left[\tau\sum_{n=1}^{\ell} \left|\left(\sum_{m=1}^n\bG[\bu^m]\overline{\Delta}W_m, \bu^{n+\frac12}\right)\right|^p\right]\\\nonumber
				&\leq C_p \mE\left[\tau\sum_{n=1}^{\ell} \left\|\sum_{m=1}^n\bG[\bu^m]\overline{\Delta}W_m\right\|^{p}_{\bL^2} \|\bu^{n+\frac12}\|^p_{\bL^2}\right]\\\nonumber
				&\leq C_p \mE\left[\tau\sum_{n=1}^{\ell} \left\|\sum_{m=1}^n\bG[\bu^m]\overline{\Delta}W_m\right\|^{2p}_{\bL^2}\right] +  C_p\mE\left[\tau\sum_{n=1}^{\ell}\|\bu^{n+\frac12}\|^{2p}_{\bL^2}\right]\\\nonumber
				&\leq C_p \mE\left[\tau\sum_{n=1}^{\ell} \tau\sum_{m=1}^n\|\bG[\bu^m]\|_{\bL^2}^{2p}\right] +  C_p\mE\left[\tau\sum_{n=1}^{\ell}\|\bu^{n+\frac12}\|^{2p}_{\bL^2}\right]\\\nonumber
				&\leq C_pC_A^p \mE\left[\tau\sum_{n=1}^{\ell} \tau\sum_{m=1}^n\mathcal{Q}^p(\bu^m,\overline{\bu}^m)\right] + C_pC_A^p T^2\\\nonumber
				&\qquad+ C_p\mE\left[\tau\sum_{n=1}^{\ell}\mathcal{Q}^p(\bu^{n+1},\overline{\bu}^{n+1}) + \mathcal{Q}^p(\bu^n,\overline{\bu}^n)\right],
			\end{align*}
			where third inequality  is obtained by using the inequality from \cite[Lemma 2.3]{feng2024NS}, while the last inequality is obtained by using the assumptions \eqref{assump:F(0)} and \eqref{assump:GradientF}.
			
			Similarly, we can control $K_4$. Using the discrete the H\"older's and Cauchy-Schwarz inequalities, we obtain
			\begin{align}\label{eq3.14}
				\mE[K_4] &\leq C_p\mE\left[\tau\sum_{n=1}^{\ell} \left\|\sum_{m=1}^nD_{\bu}\bG[\bu^m]\bv^m\widehat{\Delta}W_m\right\|^{p}_{\bL^2}\|\bu^{n+\frac12}\|^{p}_{\bL^2}\right]\\\nonumber
				&\leq C_p\mE\left[\tau\sum_{n=1}^{\ell} \left\|\sum_{m=1}^nD_{\bu}\bG[\bu^m]\bv^m\widehat{\Delta}W_m\right\|^{2p}_{\bL^2}\right] + C_p\mE\left[\tau\sum_{n=1}^{\ell}\|\bu^{n+\frac12}\|^{2p}_{\bL^2}\right].
			\end{align}
			It follows from the definition of $\widehat{\Delta}W_m$ that 
			\begin{align}\label{eq20250227}
				\widehat{\Delta}W_m = \tau W(t_{m+1}) - \tau^3\sum_{k = 1}^{\tau^{-2}}W(t_{m,k}) 
				= \tau \overline{\Delta}W_{m} - \tau^3\sum_{k = 1}^{\tau^{-2}}\left[W(t_{m}) - W(t_{m,k})\right],
			\end{align}
			which  and \eqref{eq20250227} imply 
			\begin{align*}
				\mE[K_4] &\leq C_p\mE\left[\tau\sum_{n=1}^{\ell} \left\|\sum_{m=1}^nD_{\bu}\bG[\bu^m]\tau\bv^m\overline{\Delta}W_m\right\|^{2p}_{\bL^2}\right] \\\nonumber
				&\qquad+ C_p\mE\left[\tau\sum_{n=1}^{\ell} \left\|\tau\sum_{m=1}^nD_{\bu}\bG[\bu^m]\bv^m\tau^2\sum_{k = 1}^{\tau^{-2}}\left[W(t_{m}) - W(t_{m,k})\right]\right\|^{2p}_{\bL^2}\right] \\\nonumber
				&\qquad+ C_p\mE\left[\tau\sum_{n=1}^{\ell}\|\bu^{n+\frac12}\|^{2p}_{\bL^2}\right]\\\nonumber
				&:= K_{4,1} + K_{4,2} + K_{4,3}.
			\end{align*}
			
			Now, we  control $K_{4,1}$ by using \cite[Lemma 2.3]{feng2024NS}, \eqref{assump:F(0)}, \eqref{assump:GradientF}, and \eqref{eq20230703_30} as follows:
			\begin{align*}
				K_{4,1} &\leq C_p\mE\left[\tau\sum_{n=1}^{\ell} \tau\sum_{m=1}^n\|D_{\bu}\bG[\bu^m]\tau\bv^m\|^{2p}_{\bL^2}\right] \\\nonumber
				&\leq C_pC_A^{2p}\mE\left[\tau\sum_{n=1}^{\ell} \tau\sum_{m=1}^n\|\bu^{m} - \bu^{m-1}\|^{2p}_{\bL^2}\right]  \\\nonumber
				&\qquad+C_pC_A^{2p}\mE\left[\tau\sum_{n=1}^{\ell} \tau\sum_{m=1}^n\|\bG[\bu^{m}]\widehat{\Delta}W_m\|^{2p}_{\bL^2}\right]  \\\nonumber
				&\leq C_pC_A^{2p}\mE\left[\tau\sum_{n=1}^{\ell} \tau\sum_{m=1}^n\|\bu^{m} - \bu^{m-1}\|^{2p}_{\bL^2}\right]  \\\nonumber
				&\qquad+C_pC_A^{4p}\mE\left[\tau\sum_{n=1}^{\ell} \tau^{1+3p}\sum_{m=1}^n\|\bu^{m}\|^{2p}_{\bL^2}\right]  + C_pC_A^{4p}\tau^{3p}\\\nonumber
				&\leq C_pC_A^{2p}\mE\left[\tau\sum_{n=1}^{\ell} \tau\sum_{m=1}^n\left(\mathcal{Q}^p(\bu^m,\overline{\bu}^m)+\mathcal{Q}^p(\bu^{m-1},\overline{\bu}^{m-1})\right)\right]  \\\nonumber
				&\qquad+C_pC_A^{4p}\mE\left[\tau\sum_{n=1}^{\ell} \tau^{1+3p}\sum_{m=1}^n\mathcal{Q}^p(\bu^m,\overline{\bu}^m)\right]  + C_pC_A^{4p}\tau^{3p}T^2,
			\end{align*}
			where the third inequality is obtained by using Lemma \ref{increment} (ii).
			
			To control $K_{4,2}$, we  use the discrete H\"older's inequality, the assumptions   \eqref{assump:F(0)} and \eqref{assump:GradientF}, the fact that $\mE[|W(t) - W(s)|^{2p}] \leq C|t - s|^{p}$ from \cite[Corollary 1.1]{Ichikawa}, and  \eqref{eq20230703_30} to get 
			\begin{align*}
				K_{4,2} &\leq C_p\mE\left[\tau\sum_{n=1}^{\ell} \tau\sum_{m=1}^n\|D_{\bu}\bG[\bu^m]\bv^m\tau^2\sum_{k = 1}^{\tau^{-2}}\left[W(t_{m}) - W(t_{m,k})\right]\|^{2p}_{\bL^2}\right] \\\nonumber
				&\leq C_pC_A^{2p}\mE\left[\tau\sum_{n=1}^{\ell} \tau\sum_{m=1}^n\|\bv^m\|^{2p}_{\bL^2}\left|\tau^2\sum_{k = 1}^{\tau^{-2}}\left[W(t_{m}) - W(t_{m,k})\right]\right|^{2p}\right] \\\nonumber
				&\leq C_pC_A^{2p}\mE\left[\tau\sum_{n=1}^{\ell} \tau\sum_{m=1}^n\|\bv^m\|^{2p}_{\bL^2}\tau^2\sum_{k = 1}^{\tau^{-2}}\left|W(t_{m}) - W(t_{m,k})\right|^{2p}\right] \\\nonumber
				&\leq C_pC_A^{2p}\mE\left[\tau\sum_{n=1}^{\ell} \tau\sum_{m=1}^n\|\bv^m\|^{2p}_{\bL^2}\tau^2\sum_{k = 1}^{\tau^{-2}}\left|t_m - t_{m,k}\right|^{p}\right] \\\nonumber
				&\leq C_pC_A^{2p}\mE\left[\tau\sum_{n=1}^{\ell} \tau\sum_{m=1}^n\|\bv^m\|^{2p}_{\bL^2} \tau^{3p}\right] \\\nonumber
				&= C_pC_A^{2p}\mE\left[\tau\sum_{n=1}^{\ell} \tau\sum_{m=1}^n\|\bu^m - \bu^{m-1}\|^{2p}_{\bL^2} \tau^{p}\right] \\\nonumber
				&\qquad+ C_pC_A^{2p}\mE\left[\tau\sum_{n=1}^{\ell} \tau\sum_{m=1}^n\|\bG[\bu^m]\widehat{\Delta}W_m\|^{2p}_{\bL^2} \tau^{p}\right] \\\nonumber
				&\leq C_pC_A^{2p}\mE\left[\tau\sum_{n=1}^{\ell} \tau\sum_{m=1}^n\left(\mathcal{Q}^p(\bu^m,\overline{\bu}^{m}) + \mathcal{Q}^p(\bu^{m-1},\overline{\bu}^{m-1})\right) \tau^{p}\right] \\\nonumber
				&\qquad+ C_pC_A^{4p}\mE\left[\tau\sum_{n=1}^{\ell} \tau\sum_{m=1}^n\mathcal{Q}^p(\bu^m,\overline{\bu}^m) \tau^{4p}\right] + C_pC_A^{4p}\tau^{4p}T^2,
			\end{align*}
			where the last inequality holds because of  \eqref{eq20230703_30}, \eqref{assump:F(0)}, and \eqref{assump:GradientF}. 
			
			Next, using the fact that $\|\bu^{n}\|^2_{\bL^2} \leq 2\mathcal{Q}(\bu^n,\overline{\bu}^n)$, we get
			\begin{align*}
				K_{4,3} \leq C_p\mE\left[\tau\sum_{n=1}^{\ell}(\mathcal{Q}^p(\bu^{n+1},\overline{\bu}^{n+1}) + \mathcal{Q}^p(\bu^n,\overline{\bu}^n))\right].
			\end{align*}
			
			Finally, using the discrete H\"older's inequality and Lemma \ref{increment}, we obtain
			\begin{align*}
				\mE[K_5] &\leq C_p\mE\left[\tau^{1-p}\sum_{n=1}^{\ell}\|\bG[\bu^n]\widehat{\Delta}W_n\|^{p}_{\bL^2}\|\bu^{n+\frac12}\|^{p}_{\bL^2}\right]\\\nonumber
				&\leq C_p\tau^{1-p}\sum_{n=1}^{\ell} \left(\mE\left[\|\bG[\bu^n]\widehat{\Delta}W_n\|^{2p}_{\bL^2}\right]\right)^{\frac12}\left(\mE\left[\|\bu^{n+\frac12}\|^{2p}_{\bL^2}\right]\right)^{\frac12}\\\nonumber
				&\leq C_p\tau^{1+\frac{p}{2}}\sum_{n=1}^{\ell} \left(\mE\left[\|\bG[\bu^n]\|^{2p}_{\bL^2}\right]\right)^{\frac12}\left(\mE\left[\|\bu^{n+\frac12}\|^{2p}_{\bL^2}\right]\right)^{\frac12}\\\nonumber
				&\leq C_p\tau^{1+\frac{p}{2}}\sum_{n=1}^{\ell} \left(\mE\left[\|\bG[\bu^n]\|^{2p}_{\bL^2}\right] + \mE\left[\|\bu^{n+\frac12}\|^{2p}_{\bL^2}\right]\right)
				\\\nonumber
				&\leq C_pC_A^{2p}\tau^{1+\frac{p}{2}}\sum_{n=1}^{\ell} \mE\left[\mathcal{Q}^p(\bu^n,\overline{\bu}^n)\right] + C_pC_A^{2p}T\tau^{p/2}\\\nonumber
				&\qquad+ C_p\tau^{1+p/2}\sum_{n=1}^{\ell} \mE\left[\mathcal{Q}^p(\bu^{n+1},\overline{\bu}^{n+1}) + \mathcal{Q}^p(\bu^{n-1}, \overline{\bu}^{n-1})\right].
			\end{align*}
			
			Now, combining all the estimates for $K_1, K_2, K_3, K_4$ and $K_5$,  and substituting them into \eqref{eq3.133} we obtain
			\begin{align}\label{eq3.16}
				\mE[\mathcal{Q}^p(\bu^{\ell+1},\overline{\bu}^{\ell+1})] &\leq 2^{2p}\mE\Bigl[\Bigl\{\mathcal{Q}(\bu^1,\overline{\bu}^1) + \frac{\tau^2}{4}\lambda\|\Div\overline{\bu}^{0}\|^2_{\bL^2} + \frac{\tau^2}{4}\mu\|\varepsilon(\overline{\bu}^{0})\|^2_{\bL^2}\Bigr\}^p\Bigr]\\\nonumber
				&\quad + C_p \mE\left[\|\bv_{0} - \frac{\tau}{2}\mathcal{L}\bu^0 + \bG[\bu^0]\overline{\Delta}W_0 + \tau\mathcal{L}\bu^{\frac12}\|^{2p}_{\bL^2}\right] \\\nonumber 
				&\quad+C_pC_A^p T^2 +C_pC_A^{4p}\tau^{3p}T^2+C_pC_A^{4p}\tau^{4p}T^2\\\nonumber
				&\quad+ 4C_p\mE\left[\tau\sum_{n=1}^{\ell}(\mathcal{Q}^p(\bu^{n+1},\overline{\bu}^{n+1}) + \mathcal{Q}^p(\bu^n,\overline{\bu}^n))\right]\\\nonumber
				&\quad+\tilde{C}_p\mE\left[\tau\sum_{n=1}^{\ell} \tau\sum_{m=1}^n\left(\mathcal{Q}^p(\bu^m,\overline{\bu}^m)+\mathcal{Q}^p(\bu^{m-1},\overline{\bu}^{m-1})\right)\right],
			\end{align}
			where $\tilde{C}_p=C_pC_A^p(2C_A^{p}+1 + 2C_A^{3p})$.
			
			The proof is complete after applying the discrete Gronwall's inequality.
		\end{proof}
		
		\medskip
		
		It turns out that the above stability results are not sufficient to establish error estimates for fully
		discrete finite element methods,  which require stronger stability results to be given below.  To the end, we introduce  the following new functional:
		\begin{align*}
			\widetilde{\mathcal{J}}(\bphi, \bpsi):= \left[\|\mathcal{L}\bphi\|^2_{\bL^2} + \lambda\|\Div \bpsi\|^2_{\bL^2} + \mu\|\varepsilon(\bpsi)\|^2_{\bL^2}\right].
		\end{align*}
		
		\begin{lemma}\label{lemma:time-semiform-stab-uh2-vh1}
			Let $(\bu^0, \bv^0)=(\bu_0, \bv_0)\in (\bH^2\cap\bH^1_0)\times\bH^1_0$. Under the assumptions \eqref{assump:F(0)}--\eqref{assump:Fuu}, the solution $\{(\bu^n, \bv^n);1\leq n\leq N\}$ of Algorithm 1 satisfies
			\begin{align}\label{eqn:time-semiform-stab-uh2-vh1}
				\max_{1\leq n\leq N}\mE\left[\widetilde{\mathcal{J}}(\bu^n,\bv^n)\right]
				\leq C_{3},
			\end{align}
			where $C_{3} = Ce^{CC_AT}\mE\bigl[\tJ(\bu^1,\bv^1) + \tJ(\bu^0,\bv^0)\bigr].$
		\end{lemma}		
		\begin{proof}
			Taking $\bpsi = -\frac{1}{\tau}\opL(\bu^{n+1} - \bu^{n-1})$ in \eqref{eq20230703_31}, we obtain
			\begin{align}\label{eq3.7}
				&\lambda\left(\Div\left(\frac{\bu^{n+1} - \bu^{n-1}}{\tau}\right), \Div(\bv^{n+1} - \bv^n)\right)\\\nonumber
				&\qquad +\mu\left(\veps\left(\frac{\bu^{n+1} - \bu^{n-1}}{\tau}\right), \veps(\bv^{n+1} - \bv^n)\right)  + \left(\opL \bu^{n,\frac12}, \opL(\bu^{n+1} - \bu^{n-1})\right)\\\nonumber
				&=\lambda\left(\Div \bG[\bu^n]\widehat{\Delta}W_n, \frac{1}{\tau}\Div(\bu^{n+1} - \bu^{n-1})\right) \\\nonumber
				&\qquad+ \mu\left(\veps(\bG[\bu^n]\widehat{\Delta}W_n), \frac{1}{\tau}\veps(\bu^{n+1} - \bu^{n-1})\right)\\\nonumber
				&\qquad+\lambda\left(\Div (D_{\bu}\bG[\bu^n]\bv^n)\widehat{\Delta}W_n, \frac{1}{\tau}\Div(\bu^{n+1} - \bu^{n-1})\right)\\\nonumber&\qquad + \mu\left(\veps(D_{\bu}\bG[\bu^n]\bv^n\widehat{\Delta}W_n), \frac{1}{\tau}\veps(\bu^{n+1} - \bu^{n-1})\right)\\\nonumber
				&\qquad +\lambda\left(\Div(\bF^{n,\frac12}), \Div(\bu^{n+1} - \bu^{n-1})\right) + \lambda\left(\veps(\bF^{n,\frac12}), \veps(\bu^{n+1} - \bu^{n-1})\right),
			\end{align}
			and taking $\bphi = - \opL(\bv^{n+1} - \bv^n)$ in \eqref{eq3.2}, we have 
			\begin{align}\label{eq3.8}
				&\lambda\left(\Div\left(\frac{\bu^{n+1} - \bu^{n-1}}{\tau}\right), \Div(\bv^{n+1} - \bv^n)\right)\\\nonumber
				&+\mu\left(\veps\left(\frac{\bu^{n+1} - \bu^{n-1}}{\tau}\right), \veps(\bv^{n+1} - \bv^n)\right)\\\nonumber
				&= \lambda\left(\Div(\bv^{n+1} + \bv^n),\Div(\bv^{n+1} - \bv^n)\right) + \mu\left(\veps(\bv^{n+1} + \bv^n),\veps(\bv^{n+1} - \bv^n)\right)\\\nonumber
				&\qquad-\frac{\lambda}{\tau}\left(\Div(\bG[\bu^n] \widehat{\Delta}W_n), \Div(\bv^{n+1} - \bv^n)\right)  \\\nonumber
				&\qquad-\frac{\mu}{\tau}\left(\veps(\bG[\bu^n] \widehat{\Delta}W_n), \veps(\bv^{n+1} - \bv^n)\right)\\\nonumber
				&\qquad-\frac{\lambda}{\tau}\left(\Div(\bG[\bu^{n-1}] \widehat{\Delta}W_{n-1}), \Div(\bv^{n+1} - \bv^n)\right)  \\\nonumber
				&\qquad-\frac{\mu}{\tau}\left(\veps(\bG[\bu^{n-1}] \widehat{\Delta}W_{n-1}), \veps(\bv^{n+1} - \bv^n)\right).
			\end{align}
			Substituting \eqref{eq3.8} to the left-hand side of \eqref{eq3.7} yields
			\begin{align}\label{20241218_3}
				&\lambda\left[\|\Div\bv^{n+1}\|^2_{\bL^2} - \|\Div\bv^n\|^2_{\bL^2}\right] +\mu\left[\|\veps(\bv^{n+1})\|^2_{\bL^2} - \|\veps(\bv^n)\|^2_{\bL^2}\right]\\\nonumber
				&\quad +\frac12\left[\|\opL \bu^{n+1}\|^2_{\bL^2} - \|\opL\bu^{n-1}\|^2_{\bL^2}\right]\\\nonumber
				&=\frac{\lambda}{\tau}\left(\Div(\bG[\bu^n]\widehat{\Delta}W_n+\bG[\bu^{n-1}]\widehat{\Delta}W_{n-1}),\Div(\bv^{n+1} -\bv^n)\right)\\\nonumber
				&\qquad+\frac{\mu}{\tau}\left(\veps(\bG[\bu^n]\widehat{\Delta}W_n+\bG(\bu^{n-1})\widehat{\Delta}W_{n-1}),\veps(\bv^{n+1} -\bv^n)\right)\\\nonumber
				&\qquad+\lambda\left(\Div \bG[\bu^n] \widehat{\Delta}W_n, \frac{1}{\tau}\Div(\bu^{n+1} - \bu^{n-1})\right) \\\nonumber
				&\qquad+ \mu\left(\veps(\bG[\bu^n] \widehat{\Delta}W_n), \frac{1}{\tau}\veps(\bu^{n+1} - \bu^{n-1})\right)
			\end{align}
			\begin{align}\nonumber
				&\qquad+\lambda\left(\Div D_{\bu}\bG[\bu^n] \bv^n\widehat{\Delta}W_n, \frac{1}{\tau}\Div(\bu^{n+1} - \bu^{n-1})\right)\\\nonumber&\qquad + \mu\left(\veps(D_{\bu}\bG[\bu^n] \bv^n\widehat{\Delta}W_n), \frac{1}{\tau}\veps(\bu^{n+1} - \bu^{n-1})\right)\\\nonumber
				&\qquad +\lambda\left(\Div(\bF^{n,\frac12}), \Div(\bu^{n+1} - \bu^{n-1})\right) + \lambda\left(\veps(\bF^{n,\frac12}), \veps(\bu^{n+1} - \bu^{n-1})\right)\\\nonumber
				&:=\mathcal{X}_1 + \mathcal{X}_2 +\cdots + \mathcal{X}_8.
			\end{align}
			
		On noting that  \eqref{eq20230703_31} can be rewritten as
			\begin{align}\label{eq20241218_1}
				\bv^{n+1} - \bv^n - \tau\lambda \opL\bu^{n,\frac12} 
				= \tau\bF^{n,\frac12}+\bG[\bu^n]\overline{\Delta} W_n  
				+D_{\bu}\bG[\bu^n] \bv^n\widehat{\Delta}W_n.
			\end{align}
			Then we obtain
			\begin{align}\label{eq20241218_2}
				&(\Div(\bv^{n+1} - \bv^n), \Div\bpsi) + \tau\lambda\bigl(\Div\opL \bu^{n,\frac12}, \Div \bpsi \bigr)\\
				&= \tau\bigl(\Div\bF^{n,\frac12}, \Div\bpsi\bigr) 
				+ \bigl(\Div\bG[\bu^n]\overline{\Delta} W_n, \Div\bpsi \bigr)\notag\\  &\quad+\bigl(\Div(D_{\bu}\bG[\bu^n]\bv^n)\widehat{\Delta}W_n, \Div\bpsi \bigr) \qquad\qquad \forall \, \bpsi \in \bH^2(D),\notag\\
				&(\veps(\bv^{n+1} - \bv^n), \veps(\bpsi)) + \tau\lambda \bigl(\veps\opL \bu^{n,\frac12}, \veps(\bpsi)\bigr)\label{eq20241218_5}\\
				&= \tau\bigl(\veps\bF^{n,\frac12}, \veps(\bpsi)\bigr) 
				+ \bigl(\veps\bG[\bu^n] \overline{\Delta} W_n, \veps(\bpsi)\bigr)\notag\\  &\quad+\bigl(\veps(D_{\bu}\bG[\bu^n] \bv^n)\widehat{\Delta}W_n, \veps(\bpsi)\bigr) \qquad\qquad \forall \, \bpsi \in \bH^2(D).\notag
			\end{align}
			
			Setting $\bpsi = \bG[\bu^n]\widehat{\Delta}W_n + \bG[\bu^{n-1}]\widehat{\Delta}W_{n-1}$ in \eqref{eq20241218_2} and \eqref{eq20241218_5} followed by substituting them into $\mathcal{X}_1$ and $\mathcal{X}_2$, respectively, we get
			\begin{align*}
				\mE[\mathcal{X}_1+\mathcal{X}_2] &= \frac{\lambda}{\tau}\mE\bigl[\bigl(\Div(\bG[\bu^n]\widehat{\Delta}W_n + \bG[\bu^{n-1}]\widehat{\Delta}W_{n-1}), \Div(\bv^{n+1} -\bv^n)\bigr)\bigr] \\\nonumber
				&\qquad+\frac{\mu}{\tau}\left(\veps(\bG[\bu^n] \widehat{\Delta}W_n+\bG[\bu^{n-1}]\widehat{\Delta}W_{n-1}),\veps(\bv^{n+1} -\bv^n)\right)\\\nonumber
				& = -\frac{\lambda}{2}\mE\Bigl[\bigl(\opL(\bu^{n+1} + \bu^{n-1}), \opL\bigl(\bG[\bu^n]\widehat{\Delta}W_n + \bG[\bu^{n-1}]\widehat{\Delta}W_{n-1}\bigr)\bigr)\Bigr]\\\nonumber
				&\qquad+ \frac{\lambda}{2}\mE\Bigl[\bigl(\Div(\bF[\bu^{n+1}] + \bF[\bu^{n-1}]), \Div(\bG[\bu^n]\widehat{\Delta}W_n + \bG[\bu^{n-1}] \widehat{\Delta}W_{n-1})\bigr)\Bigr]\\\nonumber
				&\qquad +\frac{\lambda}{\tau}\mE\Bigl[\bigl(\Div\bG[\bu^n] \overline{\Delta}W_{n}, \Div(\bG[\bu^n] \widehat{\Delta}W_n + \bG[\bu^{n-1}]\widehat{\Delta}W_{n-1})\bigr)\Bigr] \\\nonumber
				&\qquad+\frac{\lambda}{\tau}\mE\Bigl[\bigl(\Div D_{\bu}\bG[\bu^n]\bv^n\widehat{\Delta}W_n, \Div(\bG[\bu^n]\widehat{\Delta}W_n + \bG[\bu^{n-1}]\widehat{\Delta}W_{n-1})\bigr) \Bigr]\\\nonumber
				&\qquad+ \frac{\mu}{2}\mE\Bigl[\bigl(\veps(\bF[\bu^{n+1}] + \bF[\bu^{n-1}]), \veps(\bG[\bu^n]\widehat{\Delta}W_n + \bG[\bu^{n-1}]\widehat{\Delta}W_{n-1})\bigr)\Bigr]\\\nonumber
				&\qquad +\frac{\lambda}{\tau}\mE\Bigl[\bigl(\veps(\bG[\bu^n]\overline{\Delta}W_{n}), \veps(\bG[\bu^n]\widehat{\Delta}W_n + \bG[\bu^{n-1}]\widehat{\Delta}W_{n-1})\bigr)\Bigr] \\\nonumber
				&\qquad+\frac{\lambda}{\tau}\mE\Bigl[\bigl(\veps( D_{\bu}\bG[\bu^n]\bv^n\widehat{\Delta}W_n), \veps(\bG[\bu^n]\widehat{\Delta}W_n + \bG[\bu^{n-1}]\widehat{\Delta}W_{n-1})\bigr) \Bigr]\\\nonumber
				&:= \mathcal{X}_{1,1} +  \mathcal{X}_{1,2} +  \mathcal{X}_{1,3} +  \mathcal{X}_{1,4} +  \mathcal{X}_{1,5} +  \mathcal{X}_{1,6} +  \mathcal{X}_{1,7}.   
			\end{align*}
			
			It follows from Lemma $\ref{increment}(\rm ii)_2$, \eqref{assump:Fuu} and the Korn's inequality that
			\begin{align*}
				\mathcal{X}_{1,1}&\leq C\tau \lambda\mE\bigl[\|\opL \bu^{n+1}\|^2_{\bL^2} + \|\opL \bu^{n-1}\|^2_{\bL^2}\bigr]\\\nonumber
				&\qquad + \frac{1}{\tau} \lambda\mE\bigl[\|\opL\bigl(\bG[\bu^n] \widehat{\Delta}W_n + \bG[\bu^{n-1}] \widehat{\Delta}W_{n-1}\bigr)\|^2_{\bL^2}\bigr]\\\nonumber
				&\leq C\tau \lambda \mE\bigl[\|\opL \bu^{n+1}\|^2_{\bL^2} + \|\opL \bu^{n-1}\|^2_{\bL^2}\bigr] + CC_A^2\tau^2 \lambda \mE\bigl[\|\opL \bu^{n}\|^2_{\bL^2} + \|\opL \bu^{n-1}\|^2_{\bL^2}\bigr]\\\nonumber
				&\leq CC_A^2\tau \mE\bigl[\|\opL \bu^{n+1}\|^2_{\bL^2} +\|\opL \bu^n\|^2_{\bL^2}+ \|\opL \bu^{n-1}\|^2_{\bL^2}\bigr].\\
			\end{align*}
			Using Lemma $\ref{increment}(\rm ii)_2$, \eqref{assump:GradientF} and the Korn's inequality, we obtain
			\begin{align*}
				\mathcal{X}_{1,2}+\mathcal{X}_{1,5} &\leq CC_A\tau^{\frac32} \mE\left[\|\Div \bu^{n+1}\|^2_{\bL^2} + \|\varepsilon(\bu^{n+1})\|^2_{\bL^2}\right] \\
				&\qquad +CC_A\tau^{\frac32} \mE\left[\|\Div \bu^{n-1}\|^2_{\bL^2} + \|\varepsilon(\bu^{n-1})\|^2_{\bL^2}\right] \\
				&\qquad +CC_A\tau^{\frac32} \mE\left[\|\Div \bu^{n}\|^2_{\bL^2} + \|\varepsilon(\bu^{n})\|^2_{\bL^2}\right].
			\end{align*}
			
			Similarly, 
			we also can show
			\begin{align*}
				\mathcal{X}_{1,3}+\mathcal{X}_{1,6} \leq & CC_A\tau^2 \mE\left[\|\Div \bu^{n-1}\|^2_{\bL^2} + \|\varepsilon(\bu^{n-1})\|^2_{\bL^2}\right] \\
				&\qquad +CC_A\tau^2 \mE\left[\|\Div \bu^{n}\|^2_{\bL^2} + \|\varepsilon(\bu^{n})\|^2_{\bL^2}\right].
			\end{align*}	
			
			It follows from  \eqref{assump:GradientF}, \eqref{assump:Fuu}, Remark $\ref{increment}(\rm ii)_2$, the Korn's inequality,  and Lemma \ref{lemma3.1} that
			\begin{align*}
				\mathcal{X}_{1,4}+\mathcal{X}_{1,7} &\leq CC_A\tau \mE\bigl[\|\bv^n\|^2_{\bL^2} +CC_A\tau \mE\bigl[\|\Div\bv^n\|^2_{\bL^2}+CC_A\tau \mE\bigl[\|\veps(\bv^n)\|^2_{\bL^2}\\
				&\qquad +CC_A\tau \mE\left[\|\Div \bu^{n}\|^2_{\bL^2} + \|\varepsilon(\bu^{n})\|^2_{\bL^2}\right] + CC_A\tau \mE\left[\|\Div \bu^{n-1}\|^2_{\bL^2} \right. \\
				&\qquad \left. + \|\varepsilon(\bu^{n-1})\|^2_{\bL^2}\right]\\
				&\le CC_A\tau \mE\bigl[\|\Div\bv^n\|^2_{\bL^2}+CC_A\tau \mE\bigl[\|\veps(\bv^n)\|^2_{\bL^2}+CC_A\tau.
			\end{align*}
			
			In summary, we obtain 
			\begin{align*}
				\mE\bigl[\mathcal{X}_1\bigr] &\leq CC_A\tau \mE\left[\widetilde{\mathcal{J}}(\bu^{n+1},\bv^{n+1})\right] +CC_A\tau \mE\left[\widetilde{\mathcal{J}}(\bu^{n},\bv^{n})\right] \\
				&\qquad+ CC_A\tau \mE\left[\widetilde{\mathcal{J}}(\bu^{n-1},\bv^{n-1})\right]+ CC_A\tau.
			\end{align*}
			
			Next, we bound $\mathcal{X}_3$ and $\mathcal{X}_4$.  By using the martingale property of the increments $\overline{\Delta}W_n$, and \eqref{eq20241218_2} and \eqref{eq20241218_5} with $\bpsi = \bG[\bu^n]\widehat{\Delta}W_{n}$, we get
			\begin{align*}
				\mE\bigl[\mathcal{X}_3+\mathcal{X}_4\bigr] &=\lambda\left(\Div \bG[\bu^n]\widehat{\Delta}W_n, \frac{1}{\tau}\Div(\bu^{n+1} - \bu^{n-1})\right) \\
				&\qquad+ \mu\left(\veps(\bG[\bu^n] \widehat{\Delta}W_n), \frac{1}{\tau}\veps(\bu^{n+1} - \bu^{n-1})\right)\\
				&=\mE\bigl[\bigl(\Div\bG([\bu^n]\widehat{\Delta}W_n, \Div(\bv^{n+1} - \bv^n)\bigr)\bigr]+\mE\bigl[\bigl(\veps\bG[\bu^n]\widehat{\Delta}W_n, \veps(\bv^{n+1} - \bv^n)\bigr)\bigr]\\
				&=-\frac{\tau}{2} \mE\bigl[\bigl(\opL\bG[\bu^n]\widehat{\Delta}W_n, \opL(\bu^{n+1} + \bu^{n-1})\bigr)\bigr] \\\nonumber
				&\qquad+  \frac{\tau}{2}\mE\bigl[\bigl(\Div\bG[\bu^n]\widehat{\Delta}W_n, \Div(\bF(\bu^{n+1}) + \bF(u^{n-1}))\bigr)\bigr]\\
				&\qquad+\mE\bigl[\bigl(\Div\bG(\bu^n)\overline{\Delta}W_n,\Div\bG[\bu^n]\widehat{\Delta}W_n\bigr)\bigr)]\\
				&\qquad+\mE\bigl[\bigl(\Div\bG[\bu^n]\widehat{\Delta}W_n,\Div( D_{\bu}\bG[\bu^n]\bv^n)\widehat{\Delta}W_n\bigr)\bigr]\\\nonumber
				&=  \mathcal{X}_{2,1}  +  \mathcal{X}_{2,2} +  \mathcal{X}_{2,3} +  \mathcal{X}_{2,4}.
			\end{align*}
			
			 Using the same techniques for bounding  $ \mathcal{X}_{1,1},  \mathcal{X}_{1,2}, \mathcal{X}_{1,3}$, and $\mathcal{X}_{1,4}$, we can control each of $ \mathcal{X}_{2,1},  \mathcal{X}_{2,2}, \mathcal{X}_{2,3}$, and $\mathcal{X}_{2,4}$, hence, we get    
			\begin{align*}
				\mE\bigl[\mathcal{X}_2\bigr] &\leq CC_A\tau \mE\left[\widetilde{\mathcal{J}}(\bu^{n+1},\bv^{n+1})\right] +CC_A\tau \mE\left[\widetilde{\mathcal{J}}(\bu^{n},\bv^{n})\right] \\\nonumber
				&\qquad+ CC_A\tau \mE\left[\widetilde{\mathcal{J}}(\bu^{n-1},\bv^{n-1})\right]+ CC_A\tau.
			\end{align*}
			
			To bound $\mathcal{X}_5$ and $\mathcal{X}_6$, from  \eqref{eq20230703_30}, we get
			\begin{align}\label{eq20241220_1}
				(\Div(\bu^{n+1}-\bu^n),\Div\bphi) &= \tau (\Div\bv^{n+1},\Div\bphi) \\
				&\quad -\bigl(\Div\bG[\bu^n]\widehat{\Delta}W_n,\Div\bphi \bigr) \qquad \forall \bphi \in \bH^2(D),\notag\\
				(\veps(\bu^{n+1}-\bu^n),\veps(\bphi)) &= \tau (\veps(\bv^{n+1}),\veps(\bphi))\label{eq20241220_2}\\
				&\quad-\bigl(\veps(\bG[\bu^n]\widehat{\Delta}W_n),\veps(\bphi) \bigr) \qquad \forall \bphi \in \bH^2(D).\notag
			\end{align}            
			It follows from using the martingale property of the increments, \eqref{assump:GradientF}, and \eqref{eq20241220_1} and \eqref{eq20241220_2} with $\bphi = D_{\bu}\bG[\bu^n]\bv^n\widehat{\Delta}W_n$ that 
			\begin{align*}
				\mE\bigl[\mathcal{X}_5\bigr]+\mE\bigl[\mathcal{X}_6\bigr]=& \frac{\lambda}{\tau}\mE\bigl[\bigl(\Div (D_{\bu}\bG[\bu^n]\bv^n\widehat{\Delta}W_n), \Div(\bu^{n+1} - \bu^n)\bigr)\bigr]\\
				&\quad +\frac{\mu}{\tau}\left(\veps(D_{\bu}\bG[\bu^n]\bv^n\widehat{\Delta}W_n), \veps(\bu^{n+1} - \bu^{n-1})\right)\\\nonumber
				&=\mE\bigl[\bigl(\Div(D_{\bu}\bG[\bu^n]\bv^n\widehat{\Delta}W_n), \Div\bv^{n+1}\bigr)\bigr] \\
				&\quad- \frac{1}{\tau}\mE\bigl[\bigl(\Div(D_{\bu}\bG[\bu^n]\bv^n\widehat{\Delta}W_n),\Div\bG[\bu^n]\widehat{\Delta}W_n\bigr)\bigr]\\
				&\quad+\mE\bigl[\bigl(\veps(D_{\bu}\bG[\bu^n] \bv^n\widehat{\Delta}W_n), \veps(\bv^{n+1})\bigr)\bigr] \\
				&\quad- \frac{1}{\tau}\mE\bigl[\bigl(\veps(D_{\bu}\bG[\bu^n]\bv^n\widehat{\Delta}W_n),\veps(\bG[\bu^n]\widehat{\Delta}W_n)\bigr)\bigr]\\
				&\leq C\tau\mE\bigl[\|\Div\bv^{n+1}\|^2_{\bL^2}\bigr]+C\tau\mE\bigl[\|\veps(\bv^{n+1})\|^2_{\bL^2}\bigr]\\
				&\qquad+ CC_A\tau\mE\bigl[\| \bv^n\|^2_{\bL^2}\bigr] + C\mE\bigl[\|\Div\bG[\bu^n]\widehat{\Delta}W_n\|^2_{\bL^2}\bigr]\\
				&\leq C\tau\mE\bigl[\widetilde{\mathcal{J}}(\bu^{n+1},\bv^{n+1})\bigr] + CC_B\tau\mE\bigl[\mathcal{J}(\bu^n,\bv^n)\bigr] + CC_A\tau,
			\end{align*}
			where Lemma \ref{lemma3.1.2} is used to get the last inequality.
			
			Finally, we bound $\mathcal{X}_7$ and $\mathcal{X}_8$ as follows.  Using \eqref{eq20230703_30} with test function  $\bphi = \bF[\bu^{n+1}] + \bF[\bu^{n-1}]$ and applying \eqref{assump:GradientF}, \eqref{assump:Fuu} and Lemma \ref{lemma3.1} yield
			\begin{align*}
				\mE\bigl[\mathcal{X}_7\bigr]+\mE\bigl[\mathcal{X}_8\bigr] &= \frac{\tau}{2}\mE\bigl[\bigl(\Div(\bF[\bu^{n+1}] + \bF[\bu^{n-1}]), \Div(\bv^{n+1} + \bv^n)\bigr)\bigr] \\\nonumber
				&\qquad+\frac{\tau}{2}\mE\bigl[\bigl(\veps(\bF[\bu^{n+1}] + \bF[\bu^{n-1}]), \veps(\bv^{n+1} + \bv^n)\bigr)\bigr] \\\nonumber
				&\leq CC_A\tau\mE\bigl[\|\Div\bv^{n+1}\|^2_{\bL^2} + \|\Div\bv^n\|^2_{\bL^2}\bigr]+C\tau\mE\bigl[\|\veps(\bv^{n+1})\|^2_{\bL^2} + \|\veps(\bv^n)\|^2_{\bL^2}\bigr] \\\nonumber
				&\qquad+ CC_A\tau\mE\bigl[\|\Div(\bu^{n+1})\|^2_{\bL^2} + \|\Div(\bu^{n-1})\|^2_{\bL^2}\bigr]\\\nonumber
				&\qquad+ CC_A\tau\mE\bigl[\|\veps(\bu^{n+1})\|^2_{\bL^2} + \|\veps(\bu^{n-1})\|^2_{\bL^2}\bigr]\\\nonumber
				&\leq C\tau\mE\bigl[\widetilde{\mathcal{J}}(\bu^{n+1},\bv^{n+1}) + \widetilde{\mathcal{J}}(\bu^{n},\bv^{n})+ \widetilde{\mathcal{J}}(\bu^{n-1},\bv^{n-1})\bigr] + CC_A\tau.
			\end{align*} 
			
			Combining all the estimates for $\mathcal{X}_1, \cdots, \mathcal{X}_8$ and substituting them into \eqref{20241218_3} followed  by taking summation from $n = 1$ to $n = \ell$ for some $1 \leq \ell <N$, we get
			\begin{align}\label{eq3.26}
				\mE\bigl[\widetilde{\mathcal{J}}(\bu^{\ell+1},\bv^{\ell+1})\bigr]  
				&	\leq C\tau\sum_{n=1}^{\ell}\mE\bigl[\widetilde{\mathcal{J}}(\bu^{n+1},\bv^{n+1})\\\nonumber
				&\qquad+\mathcal{J}(\bu^{n},\bv^{n})\bigr] + \mE[\widetilde{\mathcal{J}}(\bu^0,\bv^0) + \widetilde{\mathcal{J}}(\bu^1,\bv^1)]+ CC_AT.
			\end{align}
			The proof is complete after applying the discrete Gronwall's inequality to \eqref{eq3.26}.
		\end{proof}
		
		\medskip
		
		The next lemma establishes some high moment stability estimate, which plays an important role in obtaining our error estimates in Theorem \ref{theo:fem-form-order-l2}.
		
		Define 
		\begin{align}\label{Q1}
			\widetilde{\mathcal{Q}}(\bu^{n+1},\overline{\bu}^{n+1}):=  \frac{\lambda}{2}\|\Div {\bu}^{n+1}\|^2_{\bL^2} + \frac{\mu}{2}\|\varepsilon(\bu^{n+1})\|^2_{\bL^2} + \frac{\tau^2}{4}\|\mathcal{L}\overline{\bu}^{n+1}\|^2_{\bL^2}
		\end{align}

		\begin{lemma}\label{lemma3.4}
			Let $\{(\bu^n, \bv^n)\}$ denote the numerical solutions generated by Algorithm 1. Suppose that $(\bu_{0},\bv_{0}) \in L^p(\Omega; (\bH^2\cap\bH_0^1)\times\bH_0^1)$. For any $1\leq p <\infty$, under the assumptions \eqref{assump:F(0)}--\eqref{assump:Fuu}, there exists a constant $\overline{C}_{2,p}>0$ such that
			\begin{align*}
				\max_{1 \leq n \leq N}	\mE\left[(\widetilde{\mathcal{Q}}(\bu^n,\overline{\bu}^{n}))^p\right] \leq \overline{C}_{2,p},
			\end{align*}
			where 
			\begin{align*}
				\overline{C}_{2,p} &= \biggl\{C_p\mE\Bigl[ \bigl\{\widetilde{\mathcal{Q}}(\bu^{1},\overline{\bu}^1) 
				+\frac{\tau^2}{4}  \|\mathcal{L}\overline{\bu}^{0}\|^2_{\bL^2}\bigr\}^p\Bigr]\\\nonumber
				&\quad + C_p \lambda^p \mE\left[\|\Div(\bv^0 - \frac{\tau}{2}\mathcal{L}\bu^0 + \bG[\bu^0]\overline{\Delta}W_0 + \tau\mathcal{L}\bu^{\frac12})\|^{2p}_{\bL^2}\right]\\\nonumber
				&\quad+ C_p \mu^p \mE\left[\|\varepsilon(\bv^0 - \frac{\tau}{2}\mathcal{L}\bu^0 + \bG[\bu^0]\overline{\Delta}W_0 + \tau\mathcal{L}\bu^{\frac12})\|^{2p}_{\bL^2}\right] +C_pC_A^p T^2 \biggr\}e^{(C_p+C_pC_A^{2p})T}.
			\end{align*}
		\end{lemma}
		
		\begin{proof} 
			Without loss of generality, we let $\bF \equiv 0$ in this proof to avoid some technicalities and save space. Taking $\bphi = -\mathcal{L}\bu^{n+\frac12}$ in \eqref{eq3.244} and using the integration by parts, we obtain
			\begin{align*}
				&	\frac{\lambda}{2}\left[\|\Div \bu^{n+1}\|^2_{\bL^2} - \|\Div \bu^n\|^2_{\bL^2}\right] + \frac{\mu}{2}\left[\|\varepsilon(\bu^{n+1})\|^2_{\bL^2} - \|\varepsilon(\bu^n)\|^2_{\bL^2}\right] \\\nonumber
				&\qquad+ \frac{\tau^2}{4} \left[\|\mathcal{L}\overline{\bu}^{n+1}\|^2_{\bL^2} - \|\mathcal{L}\overline{\bu}^{n-1}\|^2_{\bL^2}\right]\\\nonumber
				&= -\bigl(\bu^1 - \bu^0 + \bG[\bu^0]\widehat{\Delta} W_0 + \tau^2\opL \bu^{\frac12},\mathcal{L}\bu^{n+\frac12}\bigr)  + \bigl(\bG[\bu^n]\widehat{\Delta} W_n, \mathcal{L}\bu^{n+\frac12}\bigr)\\\nonumber 
				&\qquad - \tau\sum_{m=1}^n \bigl(\bG[\bu^m],\mathcal{L}\bu^{n+\frac12}\bigr)\overline{\Delta} W_m -\tau\sum_{m=1}^n\bigl(D_{\bu}\bG[\bu^m]\bv^m\widehat{\Delta}W_m, \mathcal{L}\bu^{n+\frac12}\bigr),
			\end{align*}
			which, using the $\widetilde{\mathcal{Q}}$ notation from \eqref{Q1}, can be rewritten as 
			\begin{align*}
				\widetilde{\mathcal{Q}}(\bu^{n+1},\overline{\bu}^{n+1}) - \widetilde{\mathcal{Q}}(\bu^{n},\overline{\bu}^n) 
				&= -\frac{\tau^2}{4} \left[\|\mathcal{L}\overline{\bu}^{n}\|^2_{\bL^2} - \|\mathcal{L}\overline{\bu}^{n-1}\|^2_{\bL^2}\right]  \\\nonumber
				&\qquad-\bigl(\bu^1 - \bu^0 + \bG[\bu^0]\widehat{\Delta} W_0 + \tau^2\opL \bu^{\frac12},\mathcal{L}\bu^{n+\frac12}\bigr) \\\nonumber
				&\qquad + \bigl(\bG[\bu^n]\widehat{\Delta} W_n, \mathcal{L}\bu^{n+\frac12}\bigr)\\\nonumber 
				&\qquad - \tau\sum_{m=1}^n \bigl(\bG[\bu^m],\mathcal{L}\bu^{n+\frac12}\bigr)\overline{\Delta} W_m \\\nonumber
				&\qquad-\tau\sum_{m=1}^n\bigl(D_{\bu}\bG[\bu^m]\bv^m\widehat{\Delta}W_m, \mathcal{L}\bu^{n+\frac12}\bigr).
			\end{align*}
			Applying the summation operator $\sum_{n=1}^{\ell}$ for any $1\leq \ell < N$, we get
			\begin{align}\label{eq_3.29}
				\widetilde{\mathcal{Q}}(\bu^{\ell+1},\overline{\bu}^{\ell+1}) &= \widetilde{\mathcal{Q}}(\bu^{1},\overline{\bu}^1) 
				-\frac{\tau^2}{4} \left[\|\mathcal{L}\overline{\bu}^{{\ell}}\|^2_{\bL^2} - \|\mathcal{L}\overline{\bu}^{0}\|^2_{\bL^2}\right]  \\\nonumber
				&\qquad-\bigl(\bu^1 - \bu^0 + \bG[\bu^0]\widehat{\Delta} W_0 + \tau^2\opL \bu^{\frac12},\sum_{n=1}^{\ell}\mathcal{L}\bu^{n+\frac12}\bigr) \\\nonumber
				&\qquad + \sum_{n=1}^{\ell}\bigl(\bG[\bu^n]\widehat{\Delta} W_n, \mathcal{L}\bu^{n+\frac12}\bigr)\\\nonumber 
				&\qquad - \tau\sum_{n=1}^{\ell}\left(\sum_{m=1}^{n}\bG[\bu^m]\overline{\Delta}W_m,\mathcal{L}\bu^{n+\frac12}\right) \\\nonumber
				&\qquad-\tau\sum_{n=1}^{\ell}\left(\sum_{m=1}^{n}D_{\bu}\bG[\bu^m]\bv^m\widehat{\Delta}W_m,\mathcal{L}\bu^{n+\frac12}\right)\\\nonumber
				&\leq \widetilde{\mathcal{Q}}(\bu^{1},\overline{\bu}^1) 
				+\frac{\tau^2}{4}  \|\mathcal{L}\overline{\bu}^{0}\|^2_{\bL^2}  \\\nonumber
				&\qquad+\left|\bigl(\bu^1 - \bu^0 + \bG[\bu^0]\widehat{\Delta} W_0 + \tau^2\opL \bu^{\frac12},\sum_{n=1}^{\ell}\mathcal{L}\bu^{n+\frac12}\bigr)\right| \\\nonumber 
				&\qquad - \tau\sum_{n=1}^{\ell}\left|\left(\sum_{m=1}^{n}\bG[\bu^m]\overline{\Delta}W_m,\mathcal{L}\bu^{n+\frac12}\right)\right| \\\nonumber
				&\qquad-\tau\sum_{n=1}^{\ell}\left|\left(\sum_{m=1}^{n}D_{\bu}\bG[\bu^m]\bv^m\widehat{\Delta}W_m,\mathcal{L}\bu^{n+\frac12}\right)\right| \\\nonumber
				&\qquad + \sum_{n=1}^{\ell}|\bigl(\bG[\bu^n]\widehat{\Delta} W_n, \mathcal{L}\bu^{n+\frac12}\bigr)|.
			\end{align}
			
			Taking the $p$th-power to \eqref{eq_3.29} for any $p\geq 1$ and use the inequality $(a+b) \leq 2^p(a^p + b^p)$ we obtain 
			\begin{align}
				\widetilde{\mathcal{Q}}^p(\bu^{\ell+1},\overline{\bu}^{\ell+1}) & \leq C_p \left\{\widetilde{\mathcal{Q}}(\bu^{1},\overline{\bu}^1) 
				+\frac{\tau^2}{4}  \|\mathcal{L}\overline{\bu}^{0}\|^2_{\bL^2}\right\}^p\\\nonumber
				&\qquad+ C_p\left|\bigl(\bu^1 - \bu^0 + \bG[\bu^0]\widehat{\Delta} W_0 + \tau^2\opL \bu^{\frac12},\sum_{n=1}^{\ell}\mathcal{L}\bu^{n+\frac12}\bigr) \right|^p\\\nonumber
				&\qquad+ C_p\left\{\tau\sum_{n=1}^{\ell}\left|\left(\sum_{m=1}^{n}\bG[\bu^m]\overline{\Delta}W_m,\mathcal{L}\bu^{n+\frac12}\right)\right|\right\}^p\\\nonumber
				&\qquad + C_p\left\{\tau\sum_{n=1}^{\ell}\left|\left(\sum_{m=1}^{n}D_{\bu}\bG[\bu^m]\bv^m\widehat{\Delta}W_m,\mathcal{L}\bu^{n+\frac12}\right)\right|\right\}^p\\\nonumber
				&\qquad + C_p\left\{\sum_{n=1}^{\ell}|\bigl(\bG[\bu^n]\widehat{\Delta} W_n, \mathcal{L}\bu^{n+\frac12}\bigr)|\right\}^{p}\\\nonumber
				&:= \widetilde{K}_1 + \cdots + \widetilde{K}_5.
			\end{align}
			
			We now employ the similar techniques for estimating $K_1,\cdots,K_5$ in the proof of Lemma \ref{lemma3.1.2} to bound $\widetilde{K}_1,\cdots,\widetilde{K}_5$ below. 
			On noting  that having the given  initial values $(\bu_0,\bv_0)$ and defining $\bu^1$ by using \eqref{eq3.3} help to control $\mE[\widetilde{K}_1]$.  Then,  using \eqref{eq3.3} and integration by parts, we obtain
			\begin{align*}
				\mE[\widetilde{K}_2] &\leq C_p\lambda^p\tau\mE\left[\left|\left(\Div\left(\bv^0 - \frac{\tau}{2}\mathcal{L}\bu^0 + \bG[\bu^0]\overline{\Delta}W_0 + \tau\mathcal{L}\bu^{\frac12}\right), \sum_{n=1}^{\ell}\Div\bu^{n+\frac12}\right)\right|^p\right]\\\nonumber
				&\qquad+C_p\mu^p\tau\mE\left[\left|\left(\varepsilon\left(\bv^0 - \frac{\tau}{2}\mathcal{L}\bu^0 + \bG[\bu^0]\overline{\Delta}W_0 + \tau\mathcal{L}\bu^{\frac12}\right), \sum_{n=1}^{\ell}\varepsilon(\bu^{n+\frac12})\right)\right|^p\right]\\\nonumber
				&\leq C_p \lambda^p \mE\left[\|\Div(\bv^0 - \frac{\tau}{2}\mathcal{L}\bu^0 + \bG[\bu^0]\overline{\Delta}W_0 + \tau\mathcal{L}\bu^{\frac12})\|^{2p}_{\bL^2}\right]\\\nonumber
				&\qquad+ C_p \mu^p \mE\left[\|\varepsilon(\bv^0 - \frac{\tau}{2}\mathcal{L}\bu^0 + \bG[\bu^0]\overline{\Delta}W_0 + \tau\mathcal{L}\bu^{\frac12})\|^{2p}_{\bL^2}\right]\\\nonumber
				&\qquad+ C_p \tau\sum_{n}^{\ell}\mE[\widetilde{\mathcal{Q}}^p(\bu^{n+1},\overline{\bu}^{n+1}) + \widetilde{\mathcal{Q}}^p(\bu^n,\overline{\bu}^n)].
			\end{align*}
			
			To bound   $\widetilde{K}_3$,  using the integration by parts, the inequality from \cite[Lemma 2.3]{feng2024NS}, the assumption \eqref{assump:GradientF}, and the Korn's inequality, we obtain
			\begin{align*}
				\mE[\widetilde{K}_3] &\leq C_p\lambda^p \mE\left[\tau\sum_{n=1}^{\ell} \left|\left(\sum_{m=1}^{n}\Div(\bG[\bu^m])\overline{\Delta}W_m, \Div\bu^{n+\frac12}\right)\right|^p\right]\\\nonumber
				&\qquad + C_p\mu^p \mE\left[\tau\sum_{n=1}^{\ell} \left|\left(\sum_{m=1}^{n}\varepsilon(\bG[\bu^m])\overline{\Delta}W_m, \varepsilon(\bu^{n+\frac12})\right)\right|^p\right]\\\nonumber
				&\leq C_p \lambda^p \tau\sum_{n=1}^{\ell}\tau\sum_{m=1}^{n}\mE\left[\|\Div\bG[\bu^m]\|^{2p}_{\bL^2}\right] + C_p \mu^p \tau\sum_{n=1}^{\ell}\tau\sum_{m=1}^{n}\mE\left[\|\varepsilon(\bG[\bu^m])\|^{2p}_{\bL^2}\right]\\\nonumber
				&\qquad+ C_p \tau\sum_{n}^{\ell}\mE[\widetilde{\mathcal{Q}}^p(\bu^{n+1},\overline{\bu}^{n+1}) + \widetilde{\mathcal{Q}}^p(\bu^n,\overline{\bu}^n)]
			\end{align*}
			\begin{align*}
				&\leq C_pC_A^{2p} \tau\sum_{n=1}^{\ell}\tau\sum_{m=1}^{n}\mE\left[\widetilde{\mathcal{Q}}^p(\bu^m,\overline{\bu}^m)\right] \\\nonumber
				&\qquad+ C_p \tau\sum_{n}^{\ell}\mE[\widetilde{\mathcal{Q}}^p(\bu^{n+1},\overline{\bu}^{n+1}) + \widetilde{\mathcal{Q}}^p(\bu^n,\overline{\bu}^n)].
			\end{align*}
			
			To estimate $\widetilde{K}_4$, we proceed following the same lines as bounding $K_4$ in the proof of Lemma \ref{lemma3.1.2}.  It follows from using the integration by parts  that 
			\begin{align}
				\mE[\widetilde{K}_4] &\leq C_p\lambda^p \mE\left[\tau\sum_{n=1}^{\ell} \left|\left(\sum_{m=1}^{n}\Div(D_{\bu}\bG[\bu^m]\bv^m)\widehat{\Delta}W_m, \Div\bu^{n+\frac12}\right)\right|^p\right]\\\nonumber
				&\qquad + C_p\mu^p \mE\left[\tau\sum_{n=1}^{\ell} \left|\left(\sum_{m=1}^{n}\varepsilon(D_{\bu}\bG[\bu^m]\bv^m)\widehat{\Delta}W_m, \varepsilon(\bu^{n+\frac12})\right)\right|^p\right]\\\nonumber
				&\leq C_p \lambda^p\tau\sum_{n=1}^{\ell} \mE\left[\left\|\sum_{m=1}^{n}\Div(D_{\bu}\bG[\bu^m]\bv^m)\widehat{\Delta}W_m\right\|^{2p}_{\bL^2}\right] \\\nonumber
				&\qquad+ C_p \mu^p\tau\sum_{n=1}^{\ell} \mE\left[\left\|\sum_{m=1}^{n}\varepsilon(D_{\bu}\bG[\bu^m]\bv^m)\widehat{\Delta}W_m\right\|^{2p}_{\bL^2}\right]\\\nonumber
				&\qquad+ C_p \tau\sum_{n}^{\ell}\mE[\widetilde{\mathcal{Q}}^p(\bu^{n+1},\overline{\bu}^{n+1}) + \widetilde{\mathcal{Q}}^p(\bu^n,\overline{\bu}^n)]\\\
				&	:=\widetilde{K}_{4,1} + \widetilde{K}_{4,2} + C_p \tau\sum_{n}^{\ell}\mE[\widetilde{\mathcal{Q}}^p(\bu^{n+1},\overline{\bu}^{n+1}) + \widetilde{\mathcal{Q}}^p(\bu^n,\overline{\bu}^n)]. \nonumber
			\end{align}
			
			It suffices to estimate $\widetilde{K}_{4,1}$ because  $\widetilde{K}_{4,2}$  can be bounded in the same way.  Using the inequality from \cite[Lemma 2.3]{feng2024NS} and the definition of $\widehat{\Delta}W_m$, we obtain
			\begin{align}
				\mE[\widetilde{K}_{4,1}]&\leq C_p \lambda^p\tau\sum_{n=1}^{\ell} \mE\left[\tau\sum_{m=1}^{n}\|\Div(D_{\bu}\bG[\bu^m]\tau\bv^m)\|^{2p}_{\bL^2}\right] \\\nonumber
				&\leq C_p \lambda^p\tau\sum_{n=1}^{\ell} \mE\left[\tau\sum_{m=1}^{n}\|\Div(D_{\bu}\bG[\bu^m](\bu^m-\bu^{m-1}))\|^{2p}_{\bL^2}\right] \\\nonumber&\qquad+C_p \lambda^p\tau\sum_{n=1}^{\ell} \mE\left[\tau\sum_{m=1}^{n}\|\Div(D_{\bu}\bG[\bu^m]\bG[\bu^m]\widehat{\Delta}W_m)\|^{2p}_{\bL^2}\right].
			\end{align}
			Using the formula $\Div({\bf A}\bu) = \Div({\bf A}) \bu + \text{trace}({\bf A}\nabla \bu)$, where $\mathbf{A}$ and $\bu$ are matrix- and vector-valued, respectively, and the assumptions \eqref{assump:GradientF} and \eqref{assump:Fuu}, we get 
			\begin{align}
				\mE[\widetilde{K}_{4,1}]
				&\leq C_p C_A^{2p}\lambda^p\tau\sum_{n=1}^{\ell} \mE\left[\tau\sum_{m=1}^{n}\|\bu^m-\bu^{m-1}\|^{2p}_{\bL^2}\right] \\\nonumber
				&\qquad+ C_p C_A^{2p}\lambda^p\tau\sum_{n=1}^{\ell} \mE\left[\tau\sum_{m=1}^{n}\|\nabla(\bu^m-\bu^{m-1})\|^{2p}_{\bL^2}\right] \\\nonumber&\qquad+C_pC_A^{2p} \lambda^p\tau\sum_{n=1}^{\ell} \mE\left[\tau^{1+3p}\sum_{m=1}^{n}\|\bG[\bu^m]\|^{2p}_{\bL^2}\right] \\\nonumber
				&\qquad+C_pC_A^{2p} \lambda^p\tau\sum_{n=1}^{\ell} \mE\left[\tau^{1+3p}\sum_{m=1}^{n}\|\nabla\bG[\bu^m]\|^{2p}_{\bL^2}\right].
			\end{align}
		It then follows from Lemma \ref{lemma3.1.2} and the Korn's inequality that 
			\begin{align*}
				\mE[\widetilde{K}_{4,1}] &\leq C_pC_A^{2p}T^2 + C_pC_A^{2p}\tau\sum_{n=1}^{\ell}\tau\sum_{m=1}^{n}\mE[\widetilde{\mathcal{Q}}^p(\bu^m,\overline{\bu}^m) + \widetilde{\mathcal{Q}}^p(\bu^{m-1},\overline{\bu}^{m-1})].
			\end{align*}
			Therefore,  $\widetilde{K}_4$ can be controlled by
			\begin{align*}
				\mE[\widetilde{K}_4] &\leq C_pC_A^{2p}T^2 + C_pC_A^{2p}\tau\sum_{n=1}^{\ell}\tau\sum_{m=1}^{n}\mE[\widetilde{\mathcal{Q}}^p(\bu^m,\overline{\bu}^m) + \widetilde{\mathcal{Q}}^p(\bu^{m-1},\overline{\bu}^{m-1})] \\\nonumber
				&\qquad+ C_p \tau\sum_{n}^{\ell}\mE[\widetilde{\mathcal{Q}}^p(\bu^{n+1},\overline{\bu}^{n+1}) + \widetilde{\mathcal{Q}}^p(\bu^n,\overline{\bu}^n)].
			\end{align*}
			
			To bound $\widetilde{K}_5$,  following the same lines as  bounding $\widetilde{K}_5$ in the proof of  Lemma \ref{lemma3.1.2}  we get 
			\begin{align}
				\mE[\widetilde{K}_5] &\leq C_p \lambda^p\mE\left[\tau^{1-p}\sum_{n=1}^{\ell}\left|\left(\Div(\bG[\bu^n]\widehat{\Delta}W_n), \Div(\bu^{n+\frac12})\right)\right|^p\right] \\\nonumber
				&\qquad+ C_p \mu^p\mE\left[\tau^{1-p}\sum_{n=1}^{\ell}\left|\left(\varepsilon(\bG[\bu^n]\widehat{\Delta}W_n), \varepsilon(\bu^{n+\frac12})\right)\right|^p\right]\\\nonumber
				&\leq C_pC_A^{2p}\tau^{1+\frac{p}{2}}\sum_{n=1}^{\ell} \mE\left[\widetilde{\mathcal{Q}}^p(\bu^n,\overline{\bu}^n)\right] \\\nonumber
				&\qquad+ C_p\tau^{1+p/2}\sum_{n=1}^{\ell} \mE\left[\widetilde{\mathcal{Q}}^p(\bu^{n+1},\overline{\bu}^{n+1}) + \widetilde{\mathcal{Q}}^p(\bu^{n-1}, \overline{\bu}^{n-1})\right],
			\end{align}
			where the last inequality is obtained by using  \eqref{assump:GradientF} and the Korn's inequality.
			
			Finally, the proof is complete after combining all the estimates for $\widetilde{K}_1,\cdots, \widetilde{K}_5$ and using the discrete Gronwall's inequality.
		\end{proof}

		\subsection{Temporal error estimates}\label{sec-3.2} 
		
		In this subsection, our goal is to establish an optimal $\mathcal{O}(\tau^{\frac32})$ order of error estimates for  $\{(\bu^n,\bv^n)\}$ generated by Algorithm 1. 
		
		\begin{theorem}\label{theorem3.3} 
			Let $(\bu,\bv)$ be the variational solution of \eqref{eq2.4} and  $\{(\bu^n,\bv^n)\}$ be the solution generated by Algorithm 1.  Suppose that $\{\bu^0, \bv^0, \bu^1\}$ satisfy \eqref{eq3.3} and \eqref{eq3.4}. 
			Then, under the assumptions \eqref{assump:F(0)}--\eqref{assump:Fuu} and $(\bu_0,\bv_0) \in L^4(\Omega; (\bH^2\cap\bH_0^1)\times \bH_0^1)$, there holds the following estimate:
			\begin{align}\label{optimal_rate}
				\max_{1 \leq n \leq N}\E\bigl[\|\bu(t_n)-\bu^n\|^2_{\bL^2}\bigr]\leq \tilde{C}_1\tau^{3},
			\end{align}
			where $\tilde{C}_1 = C(C_A + C_A\tilde{K}_{1,4}C_{s2} + C_AC_{s5} + C_{s4} + C_{s7})$.
		\end{theorem}
		
		\begin{proof} 
			We proceed following the same approach as in the proof of \cite[Theorem 4.4]{feng2024optimal} but also highlighting the differences. For the sake of completeness, we present the full proof below. Since the proof is long and for clarity, we divide it into {\em eight} steps. 
			
			\medskip 
			{\em Step 1:}
			First, substituting \eqref{eq20230703_30} into \eqref{eq20230703_31}, we obtain
			\begin{align}\label{eq3.23}
				&\bigl((\bu^{m+1} - \bu^m) - (\bu^m - \bu^{m-1}), \bphi\bigr) + \bigl(\bG[\bu^m] \widehat{\Delta} W_m - \bG[\bu^{m-1}]\widehat{\Delta} W_{m-1}, \bphi\bigr)\\\nonumber
				&\qquad \qquad +\tau^2 \bigl(\Div \bu^{m,\frac12},\Div \bphi\bigr) +\tau^2 \bigl(\veps (\bu^{m,\frac12}),\veps(\bphi)\bigr) \\\nonumber
				&\quad = \tau^2 \bigl(\bF^{m,\frac12},\bphi\bigr) + \tau \bigl(\bG[\bu^m]\overline{\Delta} W_m, \bphi\bigr) + \tau\bigl(D_{\bu}\bG[\bu^m]\bv^m\widehat{\Delta}W_m, \bphi\bigr).
			\end{align}
			Applying the summation operator $\sum_{m=1}^n$ ($1\leq n <N$) to \eqref{eq3.23} and denoting $\overline{\bu}^{n+1} := \sum_{m=1}^n \bu^{m+1}$ and $\bu^{\frac12} = \bu^{0+\frac12} = \frac12\bigl(\bu^1 + \bu^0\bigr)$, we get 
			\begin{align}\label{eq3.24}
				&	\bigl(\bu^{n+1} - \bu^n, \bphi\bigr) + \bigl(\bG[\bu^n]\widehat{\Delta} W_n, \bphi\bigr) + \tau^2\bigl(\Div\overline{\bu}^{n,\frac12},\Div\bphi\bigr) + \tau^2\bigl(\veps(\overline{\bu}^{n,\frac12}),\veps(\bphi)\bigr)\\\nonumber
				&= \bigl(\bu^1 - \bu^0,\bphi\bigr) + \bigl(\bG[\bu^0]\widehat{\Delta} W_0 + \tau^2\opL \bu^{\frac12},\bphi\bigr) + \tau^2\sum_{m=1}^n \bigl(\bF^{m,\frac12}, \bphi\bigr)\\\nonumber 
				&\qquad \qquad + \tau\sum_{m=1}^n \bigl(\bG[\bu^m],\bphi\bigr)\overline{\Delta} W_m +\tau\sum_{m=1}^n\bigl(D_{\bu}\bG[\bu^m]\bv^m\widehat{\Delta}W_m, \bphi\bigr).
			\end{align}
			
			Next, integrating \eqref{eqn:weakform-2} in time and using the fact that $\bv = 
			\bu_t$ we obtain for any $0 \leq \alpha \leq \beta  \leq T$ and $\bphi \in \bH^1_0(D)$
			\begin{align}\label{eq3.25}
				&	\bigl(\bu(\beta) - \bu(\alpha),\bphi\bigr) + \int_{\alpha}^{\beta} \int_0^s \lambda\bigl(\Div \bu(\xi), \Div \bphi\bigr)\, d\xi\,ds
				\\ \nonumber&\qquad\qquad\qquad\qquad
				+ \int_{\alpha}^{\beta} \int_0^s \mu\bigl(\veps(\bu(\xi)), \veps(\bphi)\bigr)\, d\xi\,ds \\\nonumber
				&\qquad\quad = (\beta - \alpha) \bigl(\bv_0,\bphi\bigr) + \int_{\alpha}^{\beta} \int_0^s \bigl(\bF[\bu(\xi)], \bphi\bigr)\,d\xi\,ds \\\nonumber
				&\qquad\qquad\qquad\qquad + \int_{\alpha}^{\beta} \int_0^s \bigl(\bG[\bu(\xi)],\bphi\bigr)\, dW(\xi)\, ds. 
			\end{align}
			
			Letting $\beu^n := \bu(t_n)-\bu^n$ and $\bev^n := \bv(t_n) - \bv^n$, and setting $\beta = t_{n+1}$ and $\alpha = t_n$ in \eqref{eq3.25}, and subtracting \eqref{eq3.25} from \eqref{eq3.24} yield the following error equation:
			
			\begin{align}\label{error_equation}
				&\bigl(\beu^{n+1} - \beu^n,\bphi\bigr) + \tau^2 \lambda\bigl(\Div\bbeu^{n,\frac12},\Div\bphi\bigr) + \tau^2 \mu\bigl(\veps(\bbeu^{n,\frac12}),\veps(\bphi)\bigr)\\\nonumber
				&\quad = \bigl(\tau \bv_0 - (\bu^1 - \bu^0) - \bG[\bu^0]\widehat{\Delta} W_0 + \frac{\tau^2}{2}\opL \beu^1,\bphi\bigr) \\\nonumber 
				&\qquad - \lambda\Bigl(\int_{t_n}^{t_{n+1}}\int_0^s \Div \bu(\xi)\, d\xi\, ds - \tau^2\sum_{m=1}^n \frac{\Div\bigl(\bu(t_{m+1}) + \bu(t_{m-1})\bigr)}{2},\Div\bphi\Bigr) \\\nonumber
				&\qquad - \mu\Bigl(\int_{t_n}^{t_{n+1}}\int_0^s \veps(\bu(\xi))\, d\xi\, ds - \tau^2\sum_{m=1}^n \frac{\veps\bigl(\bu(t_{m+1}) + \bu(t_{m-1})\bigr)}{2},\veps(\bphi)\Bigr) \\\nonumber
				&\qquad +\Biggl\{ \int_{t_n}^{t_{n+1}} \int_0^s \bigl(\bF[\bu(\xi)], \bphi\bigr)\, d\xi\, ds - \tau^2 \sum_{m=1}^n \bigl(\bF^{m,\frac12} , \bphi\bigr)   \Biggr\}\\\nonumber
				&\qquad + \Biggl\{\int_{t_n}^{t_{n+1}} \int_0^s \bigl(\bG[\bu(\xi)],\bphi\bigr)\, dW(\xi)\, ds - \tau\sum_{m=1}^n \bigl(\bG[\bu^m],\bphi\bigr)\overline{\Delta} W_m \\\nonumber
				&\qquad -\tau\sum_{m=1}^n\bigl(D_u\bG[\bu^m]\bv^m\widehat{\Delta}W_m, \bphi\bigr)+ \bigl(\bG[\bu^n]\widehat{\Delta} W_n, \bphi\bigr)  \Biggr\}.
			\end{align}
			
			Taking $\bphi = \beu^{n+\frac12} := \frac{1}{2}(\beu^{n+1} + \beu^n)$ in \eqref{error_equation},  we obtain 
			\begin{align}\label{error_27}
				&\bigl(\beu^{n+1} - \beu^n,\beu^{n+\frac12}\bigr) + \tau^2 \lambda\bigl(\Div\bbeu^{n,\frac12},\Div\beu^{n+\frac12}\bigr) + \tau^2 \mu\bigl(\veps(\bbeu^{n,\frac12}),\veps(\beu^{n+\frac12})\bigr)\\\nonumber
				&= \bigl(\tau \bv_0 - (\bu^1 - \bu^0) - \bG[\bu^0]\widehat{\Delta} W_0 + \frac{\tau^2}{2}\opL \beu^1,\beu^{n+\frac12}\bigr) \\\nonumber 
				& - \lambda\Bigl(\int_{t_n}^{t_{n+1}}\int_0^s \Div \bu(\xi)\, d\xi\, ds - \tau^2\sum_{m=1}^n \frac{\Div\bigl(\bu(t_{m+1}) + \bu(t_{m-1})\bigr)}{2},\Div\beu^{n+\frac12}\Bigr) \\\nonumber
				&\qquad - \mu\Bigl(\int_{t_n}^{t_{n+1}}\int_0^s \veps(\bu(\xi))\, d\xi\, ds - \tau^2\sum_{m=1}^n \frac{\veps\bigl(\bu(t_{m+1}) + \bu(t_{m-1})\bigr)}{2},\veps(\beu^{n+\frac12})\Bigr) \\\nonumber
				&\qquad +\Biggl\{ \int_{t_n}^{t_{n+1}} \int_0^s \bigl(\bF[\bu(\xi)], \beu^{n+\frac12}\bigr)\, d\xi\, ds - \tau^2 \sum_{m=1}^n \bigl(\bF^{m,\frac12} , \beu^{n+\frac12}\bigr)   \Biggr\}\\\nonumber
				&\qquad + \Biggl\{\int_{t_n}^{t_{n+1}} \int_0^s \bigl(\bG[\bu(\xi)],\beu^{n+\frac12}\bigr)\, dW(\xi)\, ds - \tau\sum_{m=1}^n \bigl(\bG[\bu^m],\beu^{n+\frac12}\bigr)\overline{\Delta} W_m \\\nonumber
				&\qquad -\tau\sum_{m=1}^n\bigl(D_u\bG[\bu^m]\bv^m\widehat{\Delta}W_m, \beu^{n+\frac12}\bigr)+ \bigl(\bG[\bu^n]\widehat{\Delta} W_n, \beu^{n+\frac12}\bigr)  \Biggr\}\\\nonumber
				&\quad =: {\tt I + II + III + IV + V}.
			\end{align}
			
			\medskip
			{\em Step 2:} 
			The two terms on the left-hand side of \eqref{error_27} can be rewritten as follows: 
			\begin{align}\label{lhs_first}
				\bigl(\beu^{n+1} - \beu^n,\beu^{n+\frac12}\bigr) &= \frac12\bigl[\|\beu^{n+1}\|^2_{\bL^2} - \|\beu^{n}\|^2_{\bL^2}\bigr], \\  \label{lhs_second}
				\tau^2 \lambda\bigl(\Div \bbeu^{n,\frac12},\Div \beu^{n+\frac12}\bigr) &= \frac{\tau^2}{4}\lambda \bigl(\Div\bbeu^{n+1} + \Div\bbeu^{n-1}, \Div \bbeu^{n+1} - \Div\bbeu^{n-1}\bigr)\\\nonumber
				&= \frac{\tau^2}{4}\lambda\bigl[\|\Div\bbeu^{n+1}\|^2_{\bL^2}  - \|\Div\bbeu^{n-1}\|^2_{\bL^2}\bigr],\\\nonumber
				\tau^2 \mu\bigl(\veps(\bbeu^{n,\frac12}),\veps(\beu^{n+\frac12})\bigr) &= \frac{\tau^2}{4}\mu \bigl(\veps(\bbeu^{n+1}) + \veps(\bbeu^{n-1}), \veps(\bbeu^{n+1}) - \veps(\bbeu^{n-1})\bigr)\\\nonumber
				&= \frac{\tau^2}{4}\mu\bigl[\|\veps(\bbeu^{n+1})\|^2_{\bL^2}  - \|\veps(\bbeu^{n-1})\|^2_{\bL^2}\bigr].
			\end{align}
			
			We now bound each term on the right-hand side of \eqref{error_equation} from above. To control the noise term ${\tt V}$, we need to use the following fact:  
			\begin{align*}
				\int_0^s \bG[\bu(\xi)]\, dW(\xi) = \biggl(\sum_{m=0}^{n-1} \int_{t_m}^{t_{m+1}} + \int_{t_n}^s\biggr) \bG[\bu(\xi)]\, dW(\xi), \qquad t_n \leq s \leq t_{n+1}. 
			\end{align*}
			Therefore, the first two terms of {\tt V} can be written as
			\begin{align}\label{eq3.29}
				&\Bigl(\int_{t_n}^{t_{n+1}}\int_0^s \bG[\bu(\xi)]\, dW(\xi)\, ds - \tau \sum_{m=1}^n \bG[\bu^m]\overline{\Delta} W_m, \beu^{n+\frac12}\Bigr) \\\nonumber
				&\quad =  \Bigl(\int_{t_n}^{t_{n+1}} \sum_{m=1}^{n-1} \int_{t_m}^{t_{m+1}} \bigl[\bG[\bu(\xi)] - \bG[\bu(t_m)]\bigr]\, dW(\xi)\, ds, \beu^{n+\frac12}\Bigr)  \\\nonumber
				&\qquad\quad +\Bigl(\int_{t_n}^{t_{n+1}} \sum_{m=1}^{n-1} \int_{t_m}^{t_{m+1}} \bigl[\bG[\bu(t_m)] - \bG[\bu^m]\bigr]\, dW(\xi)\, ds, \beu^{n+\frac12}\Bigr)  \\\nonumber
				&\qquad\quad +\Bigl(\int_{t_n}^{t_{n+1}} \int_{t_n}^s \bigl[\bG[\bu(\xi)] - \bG[\bu(t_n)]\bigr]\, dW(\xi)\, ds, \beu^{n+\frac12}\Bigr)\\\nonumber
				&\qquad\quad +\Bigl(\bG[\bu(t_n)]\int_{t_n}^{t_{n+1}} \bigl(W(s) - W(t_n)\bigr)\,ds, \beu^{n+\frac12}\Bigr) \\\nonumber
				&\qquad\quad  +\tau \Bigl(\int_{t_0}^{t_1}\bG[\bu(\xi)]\, dW(\xi), \beu^{n+\frac12}\Bigr) -  \tau\Bigl( \bG[\bu^n]\overline{\Delta} W_n, \beu^{n+\frac12} \Bigr).
			\end{align}
			Substituting \eqref{eq3.29} for the first two terms of {\tt V},  we obtain
			\begin{align}\label{eq3.30}
				{\tt V} &= \Bigl(\int_{t_n}^{t_{n+1}} \sum_{m=1}^{n-1} \int_{t_m}^{t_{m+1}} \bigl[\bG[\bu(\xi)] - \bG[\bu(t_m)]\bigr]\, dW(\xi)\, ds, \beu^{n+\frac12}\Bigr)  \\\nonumber
				&\qquad+\Bigl(\int_{t_n}^{t_{n+1}} \sum_{m=1}^{n-1} \int_{t_m}^{t_{m+1}} \bigl[\bG[\bu(t_m)] - \bG[\bu^m]\bigr]\, dW(\xi)\, ds, \beu^{n+\frac12}\Bigr)  \\\nonumber
				&\qquad+\Bigl(\int_{t_n}^{t_{n+1}} \int_{t_n}^s \bigl[\bG[\bu(\xi)] - \bG[\bu(t_n)]\bigr]\, dW(\xi)\, ds ,\beu^{n+\frac12}\Bigr)\\\nonumber
				&\qquad+\Bigl(\bG[\bu(t_n)]\int_{t_n}^{t_{n+1}} \bigl(W(s) - W(t_n)\bigr)\,ds, \beu^{n+\frac12}\Bigr)\\\nonumber 
				&\qquad - \tau \Bigl(\sum_{m=1}^n D_{\bu} \bG[\bu^m] \bv^m \widehat{\Delta} W_m, \beu^{n+\frac12}\Bigr) + \Bigl(\bG[\bu^n] \widehat{\Delta} W_n, \beu^{n+\frac12}\Bigr)  \\\nonumber
				&\qquad+\tau \Bigl(\int_{t_0}^{t_1}\bG[\bu(\xi)]\, dW(\xi), \beu^{n+\frac12}\Bigr) -  \tau\Bigl( \bG[\bu^n]\overline{\Delta} W_n, \beu^{n+\frac12} \Bigr)\\\nonumber
				&=: {\tt V_1 + V_2 + \cdots + V_8}.
			\end{align}
			
			\medskip
			{\em Step 3:} 
			Adding and subtracting $D_{\bu} \bG[\bu(t_n)] \bv(t_n) (\xi -t_n)$  to the term ${\tt V_1}$ and combining with the term ${\tt V_5}$ yield
			\begin{align*}
				&{\tt V_1 + V_5}\\\nonumber
				&= \tau\Bigl(\sum_{m=1}^{n-1} \int_{t_m}^{t_{m+1}} \bigl[\bG[\bu(\xi) - \bG[\bu(t_m)] - D_{\bu}\bG[\bu(t_m)]\bv(t_m) (\xi - t_m) \bigr]\, dW(\xi), \beu^{n+\frac12}\Bigr)  \\\nonumber
				&\qquad+  \tau\Bigl(\sum_{m=1}^{n-1} \int_{t_m}^{t_{m+1}} \bigl(D_{\bu}\bG[\bu(t_m)] \bv(t_m) - D_{\bu} \bG[\bu^m] \bv^m\bigr)(\xi - t_m)\, dW(\xi), \beu^{n+\frac12}\Bigr) \\\nonumber 
				&\qquad+ \tau \Bigl(\sum_{m=1}^{n-1} D_{\bu}\bG[\bu^m] \bv^m \bigl(\widetilde{\Delta} W_m - \widehat{\Delta} W_{m}\bigr),  \beu^{n+\frac12}\Bigr)  - \tau \bigl(D_{\bu}\bG[\bu^n]\bv^n \widehat{\Delta} W_n, \beu^{n+\frac12}\bigr)\\\nonumber
				&=: {\tt V_{1,5}^1 + V_{1,5}^2 + V_{1,5}^3 + V_{1,5}^4}.
			\end{align*}
			By the Mean Value Theorem, we obtain
			\begin{align*}
				\bG[\bu(\xi)] - \bG[\bu(t_m)] &- D_{\bu}\bG[\bu(t_m)]\bv(t_m) (\xi - t_m) \\\nonumber
				&= \Bigl(D_{\bu}\bG[\bu(\zeta)] \bv(\zeta) - D_{\bu}\bG[\bu(t_m)]\bv(t_m) \Bigr) (\xi - t_m),
			\end{align*}
			where $\zeta \in (t_m, \xi)$. 
			
			Next, using the It\^o's isometry, the above identity, the assumptions \eqref{assump:GradientF} and \eqref{assump:Fuu}, Lemma \ref{lem:holder-cont-u}, Lemma \ref{lem:holder-cont-v}, the Ladyzhenskaya's inequality and Lemma \ref{lem:highmomentstability} we obtain
			\begin{align*}
				\mE\bigl[{\tt V_{1,5}^1}\bigr] 
				&\,\,= \tau\mE\Bigl[ \Bigl(\sum_{m=1}^{n-1} \int_{t_m}^{t_{m+1}} \Bigl(D_{\bu}\bG[\bu(\zeta)] \bv(\zeta) - D_{\bu}\bG[\bu(t_m)]\bv(t_m) \Bigr) (\xi - t_m)\, dW(\xi),  \beu^{n+\frac12}\Bigr) \Bigr] \\\nonumber
				&\,\,\leq C\tau \mE\Bigl[\sum_{m=1}^{n-1}\int_{t_m}^{t_{m+1}} \|D_{\bu}\bG[\bu(\zeta)] \bv(\zeta) - D_{\bu}\bG[\bu(t_m)]\bv(\zeta)\|^2_{\bL^2} (\xi - t_m)^2\, d\xi\Bigr]\\\nonumber
				&\qquad+C\tau \mE\Bigl[\sum_{m=1}^{n-1}\int_{t_m}^{t_{m+1}} \|D_{\bu}\bG[\bu(t_m)] \bv(\zeta) - D_{\bu}\bG[\bu(t_m)]\bv(t_m)\|^2_{\bL^2} (\xi - t_m)^2\, d\xi\Bigr]\\\nonumber
				&\qquad + C\tau \mE\bigl[\|\beu^{n+1}\|^2_{\bL^2} + \|\beu^{n}\|^2_{\bL^2}\bigr]\\\nonumber
				&\,\,\leq CC_A\tau \mE\Bigl[\sum_{m=1}^{n-1}\int_{t_m}^{t_{m+1}} \|(\bu(\zeta) - \bu(t_m))\bv(\zeta)\|^2_{\bL^2} (\xi - t_m)^2\, d\xi\Bigr] \\\nonumber
				&\qquad+CC_A\tau \mE\Bigl[\sum_{m=1}^{n-1}\int_{t_m}^{t_{m+1}} \|\bv(\zeta) - \bv(t_m)\|^2_{\bL^2} (\xi - t_m)^2\, d\xi\Bigr]  + C\tau \mE\bigl[\|\beu^{n+1}\|^2_{\bL^2} + \|\beu^{n}\|^2_{\bL^2}\bigr]\\\nonumber
				&\,\, \leq CC_A\tau^3 \sum_{m=1}^{n-1}\int_{t_m}^{t_{m+1}}\mE[\|\bu(\zeta) - \bu(t_{m})\|^2_{\bH^1}\|\bv(\zeta)\|^2_{\bH^1}]\\\nonumber
				&\qquad+ CC_AC_{s5}\tau^4+C\tau \mE\bigl[\|\beu^{n+1}\|^2_{\bL^2} + \|\beu^{n}\|^2_{\bL^2}\bigr] \\\nonumber
				&\,\, \leq CC_A\tau^3 \sum_{m=1}^{n-1}\int_{t_m}^{t_{m+1}}\left(\mE[\|\bu(\zeta) - \bu(t_{m})\|^4_{\bH^1}]\right)^{1/2}\left(\mE[\|\bv(\zeta)\|^4_{\bH^1}]\right)^{1/2}\\\nonumber
				&\qquad+ CC_AC_{s5}\tau^4+C\tau \mE\bigl[\|\beu^{n+1}\|^2_{\bL^2} + \|\beu^{n}\|^2_{\bL^2}\bigr]\\\nonumber
				&\,\, \leq CC_AC_{s2}\tilde{K}_{1,2}\tau^5+ CC_AC_{s5}\tau^4+C\tau \mE\bigl[\|\beu^{n+1}\|^2_{\bL^2} + \|\beu^{n}\|^2_{\bL^2}\bigr].
			\end{align*}
			Similarly, we can show
			\begin{align*}
				\mE\bigl[{\tt V_{1,5}^2}\bigr] \leq C\tau^4 \sum_{m=1}^{n-1} \mE\bigl[\|\beu^m\|^2_{\bL^2}\bigr] + C\tau \mE\bigl[\|\beu^{n+1}\|^2_{\bL^2} + \|\beu^{n}\|^2_{\bL^2}\bigr].
			\end{align*}
			
			To bound ${\tt V_{1,5}^3}$,  using Remark \ref{increment}, Lemma $\ref{lemma3.1}$, the assumption \eqref{assump:GradientF}, and the independence property of the  increments $\{\Delta W_n\}$,  we obtain
			\begin{align*}
				\mE\bigl[{\tt V_{1,5}^3}\bigr] &\leq \tau\mE\biggl[\Bigl\|\sum_{m=1}^{n-1} D_{\bu}\bG[\bu^m] \bv^m \bigl(\widetilde{\Delta}W_m - \widehat{\Delta} W_m\bigr)\Bigr\|_{\bL^2}\|\beu^{n+\frac12}\|_{\bL^2}\biggr]\\\nonumber
				&\leq C\tau\mE\biggl[\Bigl\|\sum_{m=1}^{n-1} D_{\bu}\bG[\bu^m] \bv^m \bigl(\widetilde{\Delta}W_m - \widehat{\Delta} W_m\bigr)\Bigr\|^2_{\bL^2}\biggr] + C\tau \mE\bigl[\|\beu^{n+\frac12}\|^2_{\bL^2}\bigr]\\\nonumber
				&\leq C\tau \sum_{m=1}^{n-1} \mE\bigl[\|D_{\bu}\bG[\bu^m] \bv^m\bigl( \widetilde{\Delta} W_m - \widehat{\Delta} W_m\bigr) \|^2_{\bL^2}\bigr] +C \tau \mE\bigl[\|\beu^{n+\frac12}\|^2_{\bL^2}\bigr] \\\nonumber
				&\leq CC_A\tau^6\sum_{m=1}^{n-1} \mE\bigl[\|\bv^m\|^2_{\bL^2}\bigr] + C\tau \mE\bigl[\|\beu^{n+\frac12}\|^2_{\bL^2}\bigr] \leq CC_AC_1 \tau^5 +C \tau \mE\bigl[\|\beu^{n+\frac12}\|^2_{\bL^2}\bigr].
			\end{align*}
			Similarly, we can show 
			\begin{align*}
				\mE\bigl[{\tt V_{1,5}^4}\bigr] &\leq C\tau\mE\bigl[\|D_{\bu}\bG[\bu^n] \bv^n \widehat{\Delta} W_n\|^2_{\bL^2}\bigr] + C\tau\mE\bigl[\|\beu^{n+\frac12}\|^2_{\bL^2}\bigr] \leq CC_AC_1\tau^4 + C\tau \mE\bigl[\|\beu^{n+\frac12}\|^2_{\bL^2}\bigr].
			\end{align*}
			This then completes bounding ${\tt V_1 + V_5}$. 
			
			\medskip 
			{\em Step 4:}		
			Using Cauchy-Schwarz inequality, the assumption \eqref{assump:GradientF}, It\^o's isometry, and Lemma \ref{lem:holder-cont-u},  we get 
			\begin{align*}
				\mE\bigl[{\tt V_2}\bigr] &= \mE\biggl[\Bigl(\int_{t_n}^{t_{n+1}} \sum_{m=1}^{n-1} \int_{t_m}^{t_{m+1}} \bigl[\bG[\bu(t_m)] - \bG[\bu^m]\bigr]\, dW(\xi)\, ds, \beu^{n+\frac12} \Bigr) \biggr]\\\nonumber
				&\leq CC_A\tau^2\sum_{m=1}^{n-1}\mE\bigl[\|\beu^m\|^2_{\bL^2}\bigr] + C\tau\mE\bigl[\|\beu^{n+\frac12}\|^2_{\bL^2}\bigr].\\\nonumber
				\mE\bigl[{\tt V_3}\bigr] &= \mE\biggl[\Bigl(\int_{t_n}^{t_{n+1}} \int_{t_n}^s \bigl[\bG[\bu(\xi)] - \bG[\bu(t_n)]\bigr]\, dW(\xi)\, ds , \beu^{n+\frac12}\Bigr)\biggr]\\\nonumber
				&\leq C\tau^4 + C\tau\mE\bigl[\|\beu^{n+\frac12}\|^2_{\bL^2}\bigr].
			\end{align*}
			
			Now we bound ${\tt V_4 + V_6 + V_8}$ by using the following rewriting:
			\begin{align*}
				{\tt V_4 + V_6 + V_8}&= \Bigl(\bG[\bu(t_n)]\int_{t_n}^{t_{n+1}} \bigl(W(s) - W(t_n)\bigr)\,ds, \beu^{n+\frac12}\Bigr) \\\nonumber
				&\qquad+ \Bigl(\bG[\bu^n] \widehat{\Delta} W_n, \beu^{n+\frac12}\Bigr) - \tau\Bigl(\bG[\bu^n] \overline{\Delta} W_n, \beu^{n+\frac12}\Bigr)\\\nonumber
				&= \Bigl(\bG[\bu(t_n)]\int_{t_n}^{t_{n+1}} \bigl(W(s) - W(t_{n+1})\bigr)\,ds, \beu^{n+\frac12}\Bigr) \\\nonumber
				&\qquad+\Bigl(\bG[\bu(t_n)]\int_{t_n}^{t_{n+1}} \bigl(W(t_{n+1}) - W(t_n)\bigr)\,ds, \beu^{n+\frac12}\Bigr) \\\nonumber
				&\qquad+ \Bigl(\bG[\bu^n] \widehat{\Delta} W_n, e_u^{n+\frac12}\Bigr) - \tau\Bigl(\bG[\bu^n] \overline{\Delta} W_n, \beu^{n+\frac12}\Bigr)\\\nonumber
				&=-\Bigl(\bG[\bu(t_n)] \widetilde{\Delta} W_n, \beu^{n+\frac12}\Bigr) + \tau\Bigl(\bG[\bu(t_n)] \overline{\Delta} W_n, \beu^{n+\frac12}\Bigr) \\\nonumber
				&\qquad+ \Bigl(\bG[\bu^n] \widehat{\Delta} W_n, \beu^{n+\frac12}\Bigr) - \tau\Bigl(\bG[\bu^n] \overline{\Delta} W_n,\beu^{n+\frac12}\Bigr)\\\nonumber
				&= \tau \Bigl(\bigl(\bG[\bu(t_n)] - \bG[\bu^n]\bigr)\, \overline{\Delta} W_n, \beu^{n+\frac12} \Bigr) \\\nonumber
				&\qquad- \Bigl(\bigl(\bG[\bu(t_n)] - \bG[\bu^n]\bigr)\, \widetilde{\Delta} W_n, \beu^{n+\frac12}\Bigr) + \Bigl(\bG[\bu^n] \bigl(\widehat{\Delta} W_n - \widetilde{\Delta} W_n\bigr), \beu^{n+\frac12}\Bigr).
			\end{align*}
			Taking the expectation on both sides and using the assumption \eqref{assump:GradientF} and Lemma  \ref{increment}, we obtain
			\begin{align*}
				\mE\bigl[{\tt V_4 + V_6 + V_8}\bigr] &= \mE\biggl[\tau \Bigl(\bigl(\bG[\bu(t_n)] - \bG[\bu^n]\bigr)\, \overline{\Delta} W_n, \beu^{n+\frac12} \Bigr) \\\nonumber
				&\qquad- \Bigl(\bigl(\bG[\bu(t_n)] - \bG[\bu^n]\bigr)\, \widetilde{\Delta} W_n, \beu^{n+\frac12}\Bigr)\\\nonumber
				&\qquad + \Bigl(\bG[\bu^n] \bigl(\widehat{\Delta} W_n - \widetilde{\Delta} W_n\bigr), \beu^{n+\frac12}\Bigr)\biggr]\\\nonumber
				&\leq CC_A\tau^2 \mE\bigl[\|\beu^n\|^2_{\bL^2}\bigr] + C\tau\mE\bigl[\|\beu^{n+\frac12}\|^2_{\bL^2}\bigr]  \\\nonumber
				&\qquad + \frac{C}{\tau}\mE\bigl[\|\bigl(\bG[\bu(t_n)] - \bG[\bu^n]\bigr)\, \widetilde{\Delta} W_n\|^2_{\bL^2}\bigr] + C\tau\mE\bigl[\|\beu^{n+\frac12}\|^2_{\bL^2}\bigr] \\\nonumber
				&\qquad + \frac{C}{\tau} \mE\bigl[\|\bG[\bu^n] \bigl(\widehat{\Delta} W_n - \widetilde{\Delta} W_n\bigr)\|^2_{\bL^2}\bigr] + C\tau\mE\bigl[\|\beu^{n+\frac12}\|^2_{\bL^2}\bigr] \\\nonumber
				&\leq CC_A\tau^2 \mE\bigl[\|\beu^n\|^2_{\bL^2}\bigr] + C\tau\mE\bigl[\|\beu^{n+\frac12}\|^2_{\bL^2}\bigr] + CC_A\tau^4\mE\bigl[\|\bu^n\|^2_{\bL^2} + 1\bigr] \\\nonumber
				&\leq CC_A\tau^2 \mE\bigl[\|\beu^n\|^2_{\bL^2}\bigr] + C\tau\mE\bigl[\|\beu^{n+\frac12}\|^2_{\bL^2}\bigr] + CC_A\tau^4.
			\end{align*}
			
			Finally, using \eqref{assump:GradientF} and Lemma \ref{lem:holder-cont-u}, we get
			\begin{align*}
				\mE	\bigl[{\tt V_7}\bigr]  &=  \mE\Bigl[\tau \Bigl(\int_{t_0}^{t_1}\bG[\bu(\xi)]\, dW(\xi), \beu^{n+\frac12}\Bigr) \Bigr]\\\nonumber
				&=  \mE\Bigl[\tau \Bigl(\int_{t_0}^{t_1}[\bG[\bu(\xi)] - \bG[\bu(0)]]\, dW(\xi), \beu^{n+\frac12}\Bigr) \Bigr] + \tau\mE\bigl[\bigl(\bG[\bu(0)]\, \overline{\Delta} W_0, \beu^{n+\frac12}\bigr)\bigr] \\\nonumber
				&\leq CC_A\tau^4 + C\tau\mE\bigl[\|\beu^{n+\frac12}\|^2_{\bL^2}\bigr] + \tau\mE\bigl[\bigl(\bG[\bu(0)]\, \overline{\Delta} W_0, \beu^{n+\frac12}\bigr)\bigr].
			\end{align*}
			We note that the last term 
			will be combined with the term {\tt I} in \eqref{error_27} and handled together later.
			
			In summary, we have obtained the following upper bound for $\mE\bigl[{\tt IV}\bigr]$:
			\begin{align*}
				\mE\bigl[{\tt V}\bigr] &\leq CC_A\tau^2 \mE\bigl[\|\beu^n\|^2_{\bL^2}\bigr] + C\tau\mE\bigl[\|\beu^{n+\frac12}\|^2_{\bL^2}\bigr] + C\tau^2\sum_{m=1}^{n-1}\mE\bigl[\|\beu^m\|^2_{\bL^2}\bigr]\\\nonumber
				&\qquad+ CC_A\tau^4 + \tau\mE\bigl[\bigl(\bG[\bu(0)]\, \overline{\Delta} W_0, \beu^{n+\frac12}\bigr)\bigr].
			\end{align*}
			
			\medskip 
			{\em Step 5:}
			Our attention now shifts to bounding the term {\tt II + III} in \eqref{error_27}. To the end, we borrow an idea of \cite{feng2022higher,feng2024optimal} by using Lemma \ref{approx_integral} to do the job. First, we state the following identity:
			\begin{align}\label{eq_3.33}
				\int_0^s f(\xi)\,d\xi &:= \biggl(\frac12\sum_{m=1}^n \int_{t_{m-1}}^{t_{m+1}} + \frac12 \int_{t_0}^{t_1} + \frac12 \int_{t_n}^{t_{n+1}} - \int_s^{t_{n+1}}\biggr) f(\xi)\, d\xi\\\nonumber
				&=\frac12\biggl(\sum_{m=1}^n \int_{t_{m-1}}^{t_{m+1}} +  \int_{t_0}^{t_1} +  \int_{t_n}^{s} - \int_s^{t_{n+1}}\biggr) f(\xi)\, d\xi.
			\end{align}
			Using this identity, we obtain
			\begin{align*}
				{\tt II} &= \frac{\lambda}{2}\int_{t_n}^{t_{n+1}}\sum_{m=1}^n\biggl(  \int_{t_{m-1}}^{t_{m+1}} \Div \bu(\xi)\, d\xi -  \frac{2\tau\Div\bigl[\bu(t_{m+1}) + \bu(t_{m-1})\bigr]}{2}, \Div \beu^{n+\frac12}\biggr) \, ds\\\nonumber
				& \qquad+ \frac{\lambda}{2} \int_{t_n}^{t_{n+1}}\biggl(\Bigl(\int_{t_0}^{t_{1}} + \int_{t_n}^s - \int_s^{t_{n+1}}\Bigr) \Div \bu(\xi)\, d\xi, \Div \beu^{n+\frac12}\biggr)\, ds,
			\end{align*}
			and
			\begin{align*}
				{\tt III} &= \frac{\mu}{2}\int_{t_n}^{t_{n+1}}\sum_{m=1}^n\biggl(  \int_{t_{m-1}}^{t_{m+1}} \veps(\bu(\xi))\, d\xi -  \frac{2\tau\veps\bigl(\bu(t_{m+1}) + \bu(t_{m-1})\bigr)}{2}, \veps (\beu^{n+\frac12})\biggr) \, ds\\\nonumber
				& \qquad+ \frac{\mu}{2} \int_{t_n}^{t_{n+1}}\biggl(\Bigl(\int_{t_0}^{t_{1}} + \int_{t_n}^s - \int_s^{t_{n+1}}\Bigr) \veps(\bu(\xi))\, d\xi, \veps (\beu^{n+\frac12})\biggr)\, ds.
			\end{align*}
			
			Thus,   we obtain
			\begin{align*}
				{\tt II} &+ {\tt III} \\
				&= \biggl\{\frac{\lambda}{2}\int_{t_n}^{t_{n+1}}\sum_{m=1}^n\biggl(  \int_{t_{m-1}}^{t_{m+1}} \Div \bu(\xi)\, d\xi -  \frac{2\tau\Div\bigl[\bu(t_{m+1}) + \bu(t_{m-1})\bigr]}{2}, \Div \beu^{n+\frac12}\biggr) \, ds\\ 
				&\qquad +\frac{\mu}{2}\int_{t_n}^{t_{n+1}}\sum_{m=1}^n\biggl(  \int_{t_{m-1}}^{t_{m+1}} \veps(\bu(\xi))\, d\xi -  \frac{2\tau\veps\bigl(\bu(t_{m+1}) + \bu(t_{m-1})\bigr)}{2}, \veps (\beu^{n+\frac12})\biggr) \, ds\biggr\}\\\nonumber
				& \qquad+ \biggl\{\frac{\lambda}{2} \int_{t_n}^{t_{n+1}}\biggl(\Bigl(\int_{t_0}^{t_{1}} + \int_{t_n}^s - \int_s^{t_{n+1}}\Bigr) \Div \bu(\xi)\, d\xi, \Div \beu^{n+\frac12}\biggr)\, ds\\\nonumber & \qquad+ \frac{\mu}{2} \int_{t_n}^{t_{n+1}}\biggl(\Bigl(\int_{t_0}^{t_{1}} + \int_{t_n}^s - \int_s^{t_{n+1}}\Bigr) \veps(\bu(\xi))\, d\xi, \veps (\beu^{n+\frac12})\biggr)\, ds\biggr\}\\\nonumber
				& := X_1 + X_2.
			\end{align*}
			
			Now, let $\Phi(\xi) := \mE\bigl[\bigl(\opL\bu(\xi), \beu^{n+\frac12}\bigr)\bigr]$, then $\Phi'(\xi) =  \mE\bigl[\bigl(\opL \bv(\xi), \beu^{n+\frac12}\bigr)\bigr]$. By Lemma \ref{lem:holder-cont-v}, we have 
			\begin{align*}
				|\Phi'(t) - \Phi'(s)| \leq C_{s7}\bigl(\mE\bigl[\|\beu^{n+\frac12}\|^2_{\bL^2}\bigr]\bigr)^{\frac12} \, |t-s|^{\frac12}.
			\end{align*}
			Therefore, performing integration by parts followed by applying Lemma \ref{approx_integral}, we obtain 
			\begin{align*}
				\mE\bigl[X_1\bigr] 
				&=-\tau^2\sum_{m=1}^n \mE\biggl[\biggl(\frac{1}{2\tau}\int_{t_{m-1}}^{t_{m+1}} \opL \bu(\xi)\, d\xi -  \frac{\opL\bigl[\bu(t_{m+1}) + \bu(t_{m-1})\bigr]}{2},  \beu^{n+\frac12}\biggr)\biggr]\\\nonumber
				&\leq C_{s7}\tau^{\frac52}\bigl(\mE\bigl[\|\beu^{n+\frac12}\|^2_{\bL^2}\bigr]\bigr)^{\frac12} \leq CC_{s7}\tau^4 + \tau\mE\bigl[\|\beu^{n+\frac12}\|^2_{\bL^2}\bigr].
			\end{align*}
			
			Next, rewriting  $\mE\bigl[X_2\bigr]$ as 
			\begin{align*}
				\mE\bigl[{ X_2}\bigr] &= \biggl\{\frac{\tau\lambda}{2}\int_0^{t_1} \mE\bigl[\bigl(\Div \bu(\xi), \Div \beu^{n+\frac12}\bigr)\bigr]\, d\xi  + \frac{\tau\mu}{2}\int_0^{t_1} \mE\bigl[\bigl(\veps(\bu(\xi)), \veps(\beu^{n+\frac12})\bigr)\bigr]\, d\xi\biggr\} \\\nonumber
				&\qquad+ \biggl\{\frac{\lambda}{2} \int_{t_n}^{t_{n+1}}\mE\biggl[\biggl(\int_{t_n}^s \Div \bu(\xi)\, d\xi - \int_s^{t_{n+1}} \Div \bu(\xi)\, d\xi,\Div \beu^{n+\frac12}\biggr)\biggr]\,ds\\\nonumber
				&\qquad+ \frac{\mu}{2} \int_{t_n}^{t_{n+1}}\mE\biggl[\biggl(\int_{t_n}^s \veps(\bu(\xi))\, d\xi - \int_s^{t_{n+1}} \veps(\bu(\xi))\, d\xi,\veps(\beu^{n+\frac12})\biggr)\biggr]\,ds\biggr\}\\\nonumber
				&:={X_2^1 + X_2^2}.
			\end{align*}
			Performing integration by parts,  then adding and subtracting the term $ \opL\bu(0)$, and using Lemma \ref{lem:holder-cont-u}, we obtain
			\begin{align*}
				{X_2^1} &=  -\frac{\tau}{2}\int_0^{t_1} \mE\bigl[\bigl(\opL (\bu(\xi)- \bu(0)),  \beu^{n+\frac12}\bigr)\bigr]\, d\xi  - \frac{\tau^2}{2}\mE\bigl[\big(\opL \bu(0), \beu^{n+\frac12}\big)\bigr]\\\nonumber
				&\leq CC_{s4}\tau^5 + \tau \mE\bigl[\|\beu^{n+\frac12}\|^2_{\bL^2}\bigr] - \frac{\tau^2}{2}\mE\bigl[\big(\opL \bu(0), \beu^{n+\frac12}\big)\bigr].
			\end{align*}
			We note that the last term 
			will be combined with the term ${\tt I}$ in \eqref{error_27}.
			
			To estimate ${X_2^2}$, we proceed as follows:
			\begin{align*}
				{X_2^2} &= \frac{\lambda}{2} \int_{t_n}^{t_{n+1}}\mE\biggl[\biggl(\int_{t_n}^s \Div\bigl( \bu(\xi) - \bu(t_n)\bigr)\, d\xi,\Div \beu^{n+\frac12}\biggr)\biggr]\,ds \\\nonumber
				&\qquad+\frac{\mu}{2} \int_{t_n}^{t_{n+1}}\mE\biggl[\biggl(\int_{t_n}^s \veps\bigl( \bu(\xi) - \bu(t_n)\bigr)\, d\xi,\veps(\beu^{n+\frac12})\biggr)\biggr]\,ds \\\nonumber
				&\qquad- \frac{\lambda}{2} \int_{t_n}^{t_{n+1}}\mE\biggl[\biggl( \int_s^{t_{n+1}} \Div \bigl(\bu(\xi)-\bu(t_n)\bigr)\, d\xi,\Div \beu^{n+\frac12}\biggr)\biggr]\,ds\\\nonumber
				&\qquad- \frac{\mu}{2} \int_{t_n}^{t_{n+1}}\mE\biggl[\biggl( \int_s^{t_{n+1}} \veps \bigl(\bu(\xi)-\bu(t_n)\bigr)\, d\xi,\veps(\beu^{n+\frac12})\biggr)\biggr]\,ds\\\nonumber
				&\qquad + \frac{\lambda}{2} \int_{t_n}^{t_{n+1}}\mE\biggl[\biggl(\int_{t_n}^s \Div \bu(t_n)\, d\xi - \int_s^{t_{n+1}} \Div \bu(t_n)\, d\xi,\Div\beu^{n+\frac12}\biggr)\biggr]\,ds\\\nonumber
				&\qquad + \frac{\mu}{2} \int_{t_n}^{t_{n+1}}\mE\biggl[\biggl(\int_{t_n}^s \veps(\bu(t_n))\, d\xi - \int_s^{t_{n+1}} \veps(\bu(t_n))\, d\xi,\veps(\beu^{n+\frac12})\biggr)\biggr]\,ds.
			\end{align*}
			Performing integration by parts and using Lemma \ref{lem:holder-cont-u}, we obtain
			\begin{align*}
				{\tt X_2^2} &\leq C\int_{t_n}^{t_{n+1}} \mE\biggl[\Bigl\|\int_{t_n}^s \opL(\bu(\xi) - \bu(t_n))\, d\xi\Bigr\|^2_{\bL^2}\, ds\biggr] \\\nonumber
				&\qquad+ C\int_{t_n}^{t_{n+1}} \mE\biggl[\Bigl\|\int_{s}^{t_{n+1}} \opL(\bu(\xi) - \bu(t_n))\, d\xi\Bigr\|^2_{\bL^2}\, ds\biggr] + \tau\mE\bigl[\|\beu^{n+\frac12}\|^2_{\bL^2}\bigr]\\\nonumber
				&\qquad- \frac12\mE\biggl[\biggl(\opL \bu(t_n) \int_{t_n}^{t_{n+1}} [(s - t_n) - (t_{n+1} - s)]\, ds,  \beu^{n+\frac12}\biggr)\biggr]\\\nonumber
				&\leq CC_{s4}\tau^5 + \tau\mE\bigl[\|\beu^{n+\frac12}\|^2_{\bL^2}\bigr],
			\end{align*}
			where we have used the fact that $\int_{t_n}^{t_{n+1}} [(s-t_n) - (t_{n+1} - s)]\, ds  =0$.
			
			In summary,  we have proved that 
			\begin{align*}
				\mE[{\tt II + III}] \leq C(C_{s4} + C_{s7})\tau^4 + \tau\mE\bigl[\|\beu^{n+\frac12}\|^2_{\bL^2}\bigr] - \frac{\tau^2}{2}\mE\big[\bigl(\opL \bu_0, \beu^{n+\frac12}\bigr)\big].
			\end{align*}
			
			\medskip
			{\em Step 6:}
			Similar to the approach for bounding ${\tt II+III}$,  we rewrite ${\tt IV}$ as
			\begin{align*}
				{\tt IV} &= \int_{t_n}^{t_{n+1}} \int_0^s \bigl(\bF[\bu(\xi)], \beu^{n+\frac12}\bigr)\, d\xi\, ds - \frac{\tau^2}2 \sum_{m=1}^n \Bigl( \bF[\bu(t_{m+1})]  + \bF[\bu(t_{m-1})], \beu^{n+\frac12}\Bigr) \\\nonumber
				&\qquad+\tau^2\biggl(\sum_{m=1}^n \frac{\bF[\bu(t_{m+1})] - \bF[\bu^{m+1}] + \bF[\bu(t_{m-1})] - \bF[\bu^{m-1}]}{2}, \beu^{n+\frac12}\biggr)\\\nonumber
				&=: {\tt IV_1 + IV_2}.
			\end{align*}
			It follows from the assumption \eqref{assump:GradientF} that
			\begin{align*}
				\mE\bigl[{\tt IV_2}\bigr] \leq CC_A\tau^2\sum_{m=1}^n \bigl[\|\beu^{m+1}\|^2_{\bL^2} + \|\beu^{m-1}\|^2_{\bL^2}\bigr] + \tau \mE\bigl[\|\beu^{n+\frac12}\|^2_{\bL^2}\bigr].
			\end{align*}
			
			We now bound ${\tt IV_1}$ similarly as we did for ${ X_1}$ using the identity \eqref{eq_3.33}. 
			\begin{align*}
				{\tt IV_1} &= \frac12\int_{t_n}^{t_{n+1}}\sum_{m=1}^n\biggl(  \int_{t_{m-1}}^{t_{m+1}}  \bF[\bu(\xi)]\, d\xi -  \frac{2\tau\bigl[\bF[\bu(t_{m+1})] + \bF[\bu(t_{m-1})]\bigr]}{2},  \beu^{n+\frac12}\biggr) \, ds\\\nonumber
				& \qquad+ \frac{1}{2} \int_{t_n}^{t_{n+1}}\biggl(\Bigl(\int_{t_0}^{t_{1}} + \int_{t_n}^s - \int_s^{t_{n+1}}\Bigr)  \bF[\bu(\xi)]\, d\xi, \beu^{n+\frac12}\biggr)\, ds\\\nonumber
				&=: {\tt IV_1^1 + IV_1^2}.
			\end{align*}
			
			Using the Mean Value Theorem, the assumptions \eqref{assump:Fuu} and \eqref{assump:GradientF}, and then Poincare's inequality and the Korn's inequality, we obtain
			\begin{align*}
				& \left|\left( \partial_t\bF[\bu(t)]-  \partial_s\bF[\bu(s)], \beu^{n+\frac12}\right)\right| \\
				& \qquad =\left|\left(D_{\bu} \bF[\bu(t)] \bv(t)-D_{\bu} \bF[\bu(s)] \bv(s), \beu^{n+\frac12}\right)\right| \\
				&\qquad  =\left|\left(\left(D_{\bu} \bF[\bu(t)]-D_{\bu} \bF[\bu(s)]\right) \bv(t)+D_{\bu} \bF[\bu(s)](\bv(t)-\bv(s)), \beu^{n+\frac12}\right)\right| \\
				&\qquad  =\left|\left(\left({D_{\bu}^2 \bF}[\tilde{\bu}(\rho)](\bu(t)-\bu(s))\right)(\bv(t))+D_{\bu} \bF[\bu(s)](\bv(t)-\bv(s)), \beu^{n+\frac12}\right)\right| \\
				&\qquad  \leq {CC_A}\|\nabla(\bu(t)-\bu(s))\|_{{\bL}^2}\|\nabla\bv(t)\|_{{\bL}^2} \|\beu^{n+\frac12} \|_{{\bL}^2}+ {C_A}\|\bv(t)-\bv(s)\|_{{\bL}^2}\|\beu^{n+\frac12}\|_{{\bL}^2}\\
				&\qquad  \leq {CC_A}\|\veps(\bu(t)-\bu(s))\|_{{\bL}^2}\|\bv(t)\|_{{\bH}^1} \|\beu^{n+\frac12} \|_{{\bL}^2}+ {C_A}\|\bv(t)-\bv(s)\|_{{\bL}^2}\|\beu^{n+\frac12}\|_{{\bL}^2},
			\end{align*}
			where $\tilde{\bu}(\rho) := \rho \bu(t) + (1-\rho) \bu(s)$ for some $\rho \in [0, 1]$.
			
			Now, using Lemma \ref{lem:highmomentstability}, Lemma \ref{lem:holder-cont-u}, and Lemma \ref{lem:holder-cont-v} we get 
			\begin{align*}
				\mE\Bigl[\left|\left( \partial_t\bF[\bu(t)]-  \partial_s\bF[\bu(s)], \beu^{n+\frac12}\right)\right| \Bigr] \leq  C(C_A\tilde{K}_{1,4}C_{s2} + C_AC_{s5})\bigl(\mE\bigl[\|\beu^{n+\frac12}\|^2_{\bL^2}\bigr]\bigr)^{\frac12} |t-s|^{\frac12}.
			\end{align*}
			
			It then follows from Lemma \ref{approx_integral} that
			\begin{align*}
				\mE\bigl[{\tt IV_1^1}\bigr] \leq C(C_A\tilde{K}_{1,4}C_{s2} + C_AC_{s5})\tau^4 + \tau\mE\bigl[\|\beu^{n+\frac12}\|^2_{\bL^2}\bigr].
			\end{align*}
			
			Since ${\tt IV_1^2}$ can be bounded following the same lines as we did for ${ X_2}$, we omit the derivation. 
			In summary, we obtain
			\begin{align*}
				\mE\bigl[{\tt IV}\bigr] &\leq C(C_AC_{s2} + C_AC_{s5})\tau^4 + C\tau\mE\bigl[\|\beu^{n+\frac12}\|^2_{\bL^2}\bigr] \\\nonumber
				&\qquad+ CC_A\tau^2\sum_{m=1}^n\bigl[\|\beu^{m+1}\|^2_{\bL^2} + \|\beu^{m-1}\|^2_{\bL^2}\bigr] - \frac{\tau^2}{2}\mE\bigl[\bigl(\bF[\bu_0], \beu^{n+\frac12}\bigr)\bigr].
			\end{align*}
			
			\medskip
			{\em Step 7:}
			Collecting all the terms related to the initial data $(u^0, u^1)$ and combining them with the term {\tt I}, and then using \eqref{eq3.3} and \eqref{eq3.4}, we get 
			\begin{align*}
				\mE[{\tt \tilde{I}}] &:= \mE[{\tt I}] + \tau\mE\bigl[\bigl(\bG[\bu_0]\overline{\Delta}W_0,\beu^{n+\frac12}\bigr)\bigr] - \frac{\tau^2}{2} \mE\bigl[\bigl(\opL \bu_0, \beu^{n+\frac12}\bigr)\bigr] - \frac{\tau^2}{2} \mE\bigl[\bigl(\bF[\bu_0], \beu^{n+\frac12}\bigr)\bigr]\\\nonumber
				&= \mE\bigl[\bigl(\tau \bv_0 - (\bu^1 - \bu^0) - \bG[\bu^0]\widehat{\Delta} W_0, \beu^{n+\frac12}\bigr) \bigr]  + \frac{\tau^2}{2}\mE\bigl[\bigl(\opL \beu^1,\beu^{n+\frac12}\bigr)\bigr]\\\nonumber
				&\qquad+ \tau\mE\bigl[\bigl(\bG[\bu_0]\overline{\Delta}W_0,\beu^{n+\frac12}\bigr)\bigr] - \frac{\tau^2}{2} \mE\bigl[\bigl(\opL \bu_0, \beu^{n+\frac12}\bigr)\bigr] - \frac{\tau^2}{2} \mE\bigl[\bigl(\bF[\bu_0], \beu^{n+\frac12}\bigr)\bigr]\\\nonumber
				&=  \frac{\tau^2}{2}\mE\bigl[\bigl(\opL \beu^1, \beu^{n+\frac12}\bigr)\bigr],
			\end{align*}
			which implies that
			\begin{align*}
				\mE[{\tt \tilde{I}}] &\leq C\tau^4 + \tau\mE\bigl[\|\beu^{n+\frac12}\|^2_{\bL^2}\bigr].
			\end{align*}
			
			\medskip
			{\em Step 8:}
			Finally, substituting all the estimates for the terms {\tt I - V} into the right-hand side of \eqref{error_27} and the identities \eqref{lhs_first} and \eqref{lhs_second} into the left-hand side of \eqref{error_27}, we get
			\begin{align}\label{last_form}
				&\frac12\mE\bigl[\|\beu^{n+1}\|^2_{\bL^2} - \|\beu^{n}\|^2_{\bL^2}\bigr] + \frac{\tau^2\lambda}{4}\mE\bigl[\|\Div\bbeu^{n+1}\|^2_{\bL^2}  - \|\Div\bbeu^{n-1}\|^2_{\bL^2}\bigr] \\  \nonumber
				&\qquad\qquad\qquad\qquad\qquad+ \frac{\tau^2\mu}{4}\mE\bigl[\|\veps(\bbeu^{n+1})\|^2_{\bL^2}  - \|\veps(\bbeu^{n-1})\|^2_{\bL^2}\bigr] \\  \nonumber
				& \leq \tilde{C}_1\tau^4 + C\tau\mE\bigl[\|\beu^{n + 1}\|^2_{\bL^2} + \|\beu^{n}\|^2_{\bL^2}\bigr] + C\tau^2\sum_{m=1}^n\mE\bigl[\|\beu^{m+1}\|^2_{\bL^2} + \|\beu^{m-1}\|^2_{\bL^2}\bigr],
			\end{align}
			where $\tilde{C}_1 = C(C_A + C_A\tilde{K}_{1,4}C_{s2} + C_AC_{s5} + C_{s4} + C_{s7} )$.
			Applying the summation operator $\sum_{n=1}^{\ell}$ ($1 \leq \ell <N$) on both sides and using the discrete Gronwall's inequality yield the desired estimate \eqref{optimal_rate}.
		\end{proof}

		\medskip
		
		We finish this section by stating the following error estimates for $\{(\bu^n,\bv^n)\}$ in $\bH^1\times \bL^2$ norm. Its proof follows the same lines as the proof of \cite[Theorem 4.3]{feng2024optimal}. So, we omit it to save space.
		
		\begin{theorem}\label{theorem3.4} 
			Let $(\bu,\bv)$ be the variational solution of \eqref{eq2.4} and  $\{(\bu^n,\bv^n)\}$ be the solution generated by Algorithm 1.  Suppose 
			the initial approximations $\{\bu^0, \bv^0, \bu^1\}$ satisfy \eqref{eq3.3} and \eqref{eq3.4}.
			Then, under the assumptions \eqref{assump:F(0)}--\eqref{assump:Fuu} and $(\bu_0,\bv_0) \in L^4(\Omega; (\bH^2\cap\bH_0^1)\times \bH_0^1)$, there holds the following estimate:
			\begin{align}
				\max_{1 \leq n \leq N}\Bigl\{\mathbb{E}\left[\|\bu(t_n) - \bu^n\|_{\bH^1}^2\right]
				+\mathbb{E}\left[\|\bv(t_n) - \bv^n\|_{\bL^2}^2\right] \Bigr\}	\le \tilde{C}_2\tau^2,
			\end{align}
			where $\tilde{C}_2 = C(C_A,C_B,\tilde{C}_1)$.
		\end{theorem}
		
		\section{Fully discrete finite element method}\label{sec-4}
		In this section, we discretize the time semi-discrete scheme in space using the finite element methods and give detailed error analysis for the fully discrete scheme.
		\subsection{Formulation of finite element method}\label{ssec:fem-form}
		Let $\cal{T}_h$ be a quasi-uniform triangulation of $\cal{D}$ with mesh size $0<h<<1$. We consider the finite element spaces
		\begin{align*}
			\bUh^{r_1} = \left\{\bu_h\in\bH^1: \bu_h|_{K}\in\mbf{P}_{r_1}({K})\quad \forall K\in\cal{T}_h\right\},\\
			\bVh^{r_2} = \left\{\bv_h\in\bH^1: \bv_h|_{K}\in\mbf{P}_{r_2}({K})\quad \forall K\in\cal{T}_h\right\},
		\end{align*}
		where $\mbf{P}_r({K})$ denotes the space of polynomials with degree not exceeding a given integer $r$ on $K\in\cal{T}_h$. Next, we define two types of projection as follows.
		\begin{definition}\label{def:projection}\
			\begin{itemize}
			\item[{\rm (1)}]     We define $\bL^2$-projection $\Ph:\bL^2\to\bUh^{r_1}$ as follows.  For any $\bw\in\bL^2$,  we seek $\Ph\bw\in\bUh^{r_1}$ such that 
			\begin{align*}
				(\Ph\bw, \bv_h) = (\bw, \bv_h)  \qquad \forall \bv_h\in\bUh^{r_1}.
			\end{align*}

			\item[{\rm (2)}]The $\bH^1$-projection $\Rh:\bH^1_0\to\bVh^{r_2}$ is defined similarly. For any $\bw\in\bH^1_0$, we seek $\Rh\bw\in\bVh^{r_2}$ such that 
			\begin{align*}
				&\innerprd{\Rh\bw}{\bv_h}  \\
				&\hskip 0.8in  = \innerprd{\bw}{\bv_h} \qquad \forall \bv_h\in\bVh^{r_2}. 
			\end{align*}
			\end{itemize} 
		\end{definition}
		
		The following error estimate results are well-known \cite{baker1976error,cohen2018numerical}.
			\begin{subequations}
		\begin{align}\label{eqn:prj-error-estimate-ph}
			\left\lVert{\bw-\Ph\bw}\right\rVert_{\bL^2}\leq& C_{\Ph}h^{\min{(r_1+1, s)}}\left\lVert\bw\right\rVert_{\bH^s},\\
			\left\lVert{\bw-\Rh\bw}\right\rVert_{\bL^2}\leq& C_{\Rh}h^{\min{(r_2+1, s)}}\left\lVert\bw\right\rVert_{\bH^s},\label{20220608_2}\\
			\Bigl(\lambda\left\lVert{\Div(\bw-\Rh\bw)}\right\rVert^2_{\bL^2}+\mu\left\lVert{\veps(\bw-\Rh\bw)}\right\rVert^2_{\bL^2}\Bigr)^{\half}\leq& C_{\Rh}h^{\min{(r_2, s-1)}}\left\lVert\bw\right\rVert_{\bH^s}.\label{20220608_3}
		\end{align}
			\end{subequations}
		
		\smallskip
		
		{\bf Algorithm 2.}\, Let $N>>1$ be a positive integer and $\tau:=T/N$. Let $\{t_n\}_{n=0}^N$ be a uniform partition of  the interval $[0,T]$ with mesh size $ \tau$. Find $\mathcal{F}_{t_n}$ adapted process $\{(\bu_h^n,\bv_h^n)\}_{n=0}^N\subset \bUh^{r_1}\times\bVh^{r_2}$ such that there hold $\mathbb{P}$-almost surely
		\begin{subequations} \label{theta_scheme_FEM}
			\begin{align}
				(\bu_h^{n+1}-\bu_h^n,\bphi_h) &= \tau (\bv_h^{n+1},\bphi_h) -\bigl(\bG[\bu_h^n]\widehat{\Delta}W_n,\bphi_h \bigr) \qquad\qquad\quad  \forall \bphi_h \in \bUh^{r_1},\label{FEM1}\\
				(\bv_h^{n+1} - \bv_h^n, \bpsi_h) &+ \tau\lambda \bigl(\Div \bu_h^{n,\frac12}, \Div \bpsi_h \bigr) + \tau\mu \bigl(\veps(\bu_h^{n,\frac12}),\veps(\bpsi_h) \bigr) \label{FEM2} \\
				&= \tau\bigl(\bF_h^{n,\frac12}, \bpsi_h\bigr) 
				+ \bigl(\bG[\bu_h^n]\overline{\Delta} W_n, \bpsi_h \bigr)  \nonumber\\
				&\quad +  \bigl(D_{\bu}\bG[\bu_h^n]\bv^n_h\widehat{\Delta}W_n, \bpsi_h \bigr) \qquad\qquad\qquad\quad \forall \, \bpsi_h \in \bVh^{r_2}, \nonumber
			\end{align}
		\end{subequations}
		where 
		\begin{align*}
			\bu_h^{n,\frac12} :=  \frac12(\bu_h^{n+1} + \bu_h^{n-1}),\qquad  \bF_h^{n,\frac12} := \frac12( \bF[u_h^{n+1}] +  \bF[u_h^{n-1}]).
		\end{align*}
		At each time step, the above scheme is a nonlinear random algebraic system for $(\bu^{n+1}_h, \bv^{n+1}_h)$ whose well-posedness can be proved by a standard fixed point argument based on the stability estimates of the next lemma, which states that the fully discrete solution $\bu_h^n$ defined by Algorithm 2 also satisfies similar high moment stability estimates.   
		
		\begin{lemma}\label{lemma:fem-form-stab}
			Suppose $(\bu^0_h, \bv^0_h)=(\Rh\bu^0, \Ph\bv^0)$ and $(\bu_0,\bv_0)\in L^p(\Omega; \bH_0^1\times\bL^2)$. Let $\{(\bu_h^n,\bv_h^n)\}$ be the solution generated by Algorithm 2. For any $1\leq p <\infty$, under the assumptions \eqref{assump:F(0)}--\eqref{assump:Fuu}, there holds
			\begin{align}\label{eqn:fem-form-stab}
				\max_{1\leq n\leq N}\mE[(\mathcal{Q}(\bu_h^n, \overline{\bu}_h^n))^p]&\leq C_{3,p},
			\end{align}
			where $C_{3,p} = C(\bu_0,\bv_0,p,C_A)>0$.
		\end{lemma}
		
		\begin{lemma}\label{lemma:fem-form-stab_4.3}
			Suppose $(\bu^0_h, \bv^0_h)=(\Rh\bu^0, \Ph\bv^0)$ and $(\bu_0,\bv_0)\in L^p(\Omega; (\bH^2\cap\bH_0^1)\times\bH_0^1)$. Let $\{(\bu_h^n,\bv_h^n)\}$ be the solution generated by Algorithm 2. For any $1\leq p <\infty$, under the assumptions \eqref{assump:F(0)}--\eqref{assump:Fuu}, there holds
			\begin{align}\label{eqn:fem-form-stab_4.6}
				\max_{1\leq n\leq N}\mE[(\widetilde{\mathcal{Q}}(\bu_h^n, \overline{\bu}_h^n))^p]&\leq \overline{C}_{3,p},
			\end{align}
			where $\overline{C}_{3,p} = C(\bu_0,\bv_0,p,C_A)>0$.
		\end{lemma}
		
		 The proofs of Lemmas \ref{lemma:fem-form-stab} and \ref{lemma:fem-form-stab_4.3} follow the same lines as the proofs of Lemmas \ref{lemma3.1.2} and \ref{lemma3.4}, respectively.  We omit them to save space. Moreover, it is worth noting that the high moment stability estimates play an important role in our error analysis presented in Theorems \ref{theo:fem-form-order-h1} and \ref{theo:fem-form-order-l2}.
		
		\subsection{Error estimates for the finite element fully discrete scheme}\label{ssec:error-estimate-fem}
		
		The linear finite elements are used in this section,  that is,  $r_1 = r_2 = 1$ are chosen in the finite element space definitions. We start with the discretization for the operator $\opL$.
		
		\begin{definition}
			We define the discrete operator $\opLh:\bUh^1\to\bVh^1$ as follows. For any $\bw_h\in \bUh^1$, we seek $\opLh\bw_h\in \bVh^1$ such that 
			\begin{align}
				-\left(\opLh\bw_h, \bv_h\right)=\innerprd{\bw_h}{\bv_h}    \qquad \forall\bv_h\in\bVh^1.
			\end{align}
		\end{definition}
		
		Define $\bfErr{u}{n} := \bu^n - \bu^n_h$ and $\bfErr{v}{n} := \bv^n - \bv^n_h$, we now want to derive the estimates for $\bfErr{u}{n}$ and $\bfErr{v}{n}$.
		
		\begin{theorem}\label{theo:fem-form-order-h1}
			Let $\{(\bu^n, \bv^n)\}$ and $\{(\bu^n_h, \bv^n_h)\}$ be the solutions generated by Algorithm 1 and Algorithm 2, respectively. Let $(\bu^0_h, \bv^0_h)=(\Rh\bu^0, \Ph\bv^0)$, and $(\bu^1_h, \bv^1_h)=(\Rh\bu^1, \Ph\bv^1)$. Suppose that $(\bu^0, \bv^0)\in L^4(\Omega; (\bH^2\cap \bH_0^1)\times\bH^1_0)$. Under the assumptions \eqref{assump:F(0)}--\eqref{assump:LipsG}, we have
			\begin{equation}\label{eqn:fem-form-order-h1}
				\max\limits_{1\leq n \leq N}\exc{\normbfh{\bu^n-\bu^n_h}+\normbf{\bv^n-\bv^n_h}}\leq \tilde{C}_3\left(\tau^4 + h^2\right),
			\end{equation}
			where $\tilde{C}_3 = \max\{CC_BC_{2,2}(C_{\Rh}+C_{\Ph}), CC_A(C_{2,4} + C_{3,4})C_{2,8}\}e^{C(C_A + C_B)T}$.
		\end{theorem}
		
		\begin{proof} 
			Subtracting \eqref{theta_scheme_FEM} from \eqref{theta_scheme} leads to
			\begin{align}\label{eqn:prf-fem-form-order-h1-1}
				\big(\bfErr{u}{n+1}-\bfErr{u}{n}, \bphi_h\big) = &\tau\big(\bfErr{v}{n+1}, \bphi_h\big) \\
				&- \bigl((\bG[\bu^n]-\bG[\bu^n_h])\widehat{\Delta}W_n,\bphi_h\bigr)\qquad \forall \bphi_h \in \bUh^{r_1},\notag\\
				\big(\bfErr{v}{n+1}-\bfErr{v}{n}, \bpsi_h\big) = &\tau\big(\opL\bfErr{u}{n,\half}, \bpsi_h\big)+\big((\bG[\bu^n]-\bG[\bu^n_h])\overline{\Delta}W_n, \bpsi_h\big)\label{20220608_1}\\\nonumber
				&\qquad+\bigl((D_{\bu}\bG[\bu^n]\bv^n - D_{\bu}\bG[\bu_h^n]\bv^n_h)\, \widehat{\Delta}W_n,\bpsi_h\bigr)\\\nonumber
				&\qquad+\tau\big(\bF^{n,\half}-\bF_h^{n,\half}, \bpsi_h\big)\qquad \forall \, \bpsi_h \in \bVh^{r_2}.\notag
			\end{align}
			Rewriting \eqref{eqn:prf-fem-form-order-h1-1} as 
			\begin{align}\label{eq4.9}
				\big(\bfErr{u}{n+1}-\bfErr{u}{n-1}, \bphi_h\big) &= \bigl(\bfErr{u}{n} - \bfErr{u}{n-1},\bphi_h\bigr)+\tau\big(\bfErr{v}{n+1}, \bphi_h\big) \\\nonumber
				&\qquad- \bigl((\bG[\bu^n] - \bG[\bu^n_h])\widehat{\Delta}W_n,\bphi_h\bigr).
			\end{align}
			Choosing $\bphi_h=-\opLh\Rh\bfErr{u}{n,\half}$ in \eqref{eq4.9} we obtain
			\begin{align}\label{eqn:prf-fem-form-order-h1-2}
				\big(\bfErr{u}{n+1}-\bfErr{u}{n-1}, \opLh\Rh\bfErr{u}{n,\half}\big) &= \tau\big(\bfErr{v}{n+1}, \opLh\Rh\bfErr{u}{n,\half}\big) \\\nonumber
				&\qquad+\bigl(\bfErr{u}{n} - \bfErr{u}{n-1},\opLh\Rh\bfErr{u}{n,\half}\bigr) \\\nonumber&\qquad- \bigl((\bG[\bu^n]-\bG[\bu^n_h])\widehat{\Delta}W_n, \opLh\Rh\bfErr{u}{n,\half}\bigr).
			\end{align}
			Choosing $\bpsi_h=\Ph\bfErr{v}{n+1}$ in \eqref{20220608_1},  we get
			\begin{align}\label{eqn:prf-fem-form-order-h1-3}
				\big(\bfErr{v}{n+1}-\bfErr{v}{n}, \Ph\bfErr{v}{n+1}\big) = &\tau\big(\opL\bfErr{u}{n,\half}, \Ph\bfErr{v}{n+1}\big)\\\nonumber
				&\quad+\tau\big(\bF^{n,\half}  -\bF_h^{n,\half}, \Ph\bfErr{v}{n+1}\big)\\\nonumber
				&\quad+\big((\bG[\bu^n] -\bG[\bu^n_h])\overline{\Delta}W_n, \Ph\bfErr{v}{n+1}\big)\\\nonumber
				&\quad+\bigl((D_{\bu}\bG[\bu^n]\bv^n - D_{\bu}\bG[\bu_h^n]\bv^n_h)\, \widehat{\Delta}W_n,\Ph\bfErr{v}{n+1}\bigr).
			\end{align}
			Using the definitions of $\Ph$, $\opL$, and $\opLh$, we have
			\begin{align}\label{20220613_1}
				\big(\opL\bfErr{u}{n,\half}, \Ph\bfErr{v}{n+1}\big)&=\big(\opL\Rh\bfErr{u}{n,\half}, \Ph\bfErr{v}{n+1}\big)\\
				&=\big(\opLh\Rh\bfErr{u}{n,\half}, \Ph\bfErr{v}{n+1}\big)\notag\\
				&=\big(\opLh\Rh\bfErr{u}{n,\half}, \bfErr{v}{n+1}\big)\notag.
			\end{align}
			
			Next, the left-hand side of \eqref{eqn:prf-fem-form-order-h1-3} can be written as
			\begin{align}\label{20220613_2}
				&\big(\bfErr{v}{n+1}-\bfErr{v}{n}, \Ph\bfErr{v}{n+1}\big)\\
				&\qquad =\big(\Ph\bfErr{v}{n+1}-\Ph\bfErr{v}{n}, \Ph\bfErr{v}{n+1}\big)\notag\\
				&\qquad = \frac12\|\Ph\bfErr{v}{n+1}\|_{{\bf L}^2}^2-\frac12\|\Ph\bfErr{v}{n}\|_{{\bf L}^2}^2+\frac12\|\Ph\bfErr{v}{n+1}-\Ph\bfErr{v}{n}\|_{{\bf L}^2}^2\notag.
			\end{align}
			Combining \eqref{eqn:prf-fem-form-order-h1-2}-\eqref{20220613_2}, we have
			\begin{align}\label{eqn:prf-fem-form-order-h1-4}
				&\frac12\bigl[\|\Ph\bfErr{v}{n+1}\|_{{\bf L}^2}^2-\|\Ph\bfErr{v}{n}\|_{{\bf L}^2}^2\bigr]+\half\exc{\normbf{\Ph\bfErr{v}{n+1}-\Ph\bfErr{v}{n}}}\\
				&\notag\qquad+ \frac12\bigl[\normdel{\Rh\bfErr{u}{n+1}}\bigr] \\\nonumber&\qquad- \frac12\bigl[\normdel{\Rh\bfErr{u}{n-1}}\bigr]\\\nonumber
				&\quad =\tau\big(\bF^{n,\half} -\bF_h^{n,\half}, \Ph\bfErr{v}{n+1}\big)+\big((\bG[\bu^n]-\bG[\bu^n_h])\overline{\Delta}W_n, \Ph\bfErr{v}{n+1}\big)\\\nonumber
				&\qquad+\bigl((D_{\bu}\bG[\bu^n]\bv^n - D_{\bu}\bG[\bu_h^n]\bv^n_h)\, \widehat{\Delta}W_n,\Ph\bfErr{v}{n+1}\bigr)\\\nonumber
				&\qquad-\bigl(\bfErr{u}{n} -\bfErr{u}{n-1},\opLh\Rh\bfErr{u}{n,\half}\bigr)+ \bigl((\bG[\bu^n]-\bG[\bu^n_h])\widehat{\Delta}W_n, \opLh\Rh\bfErr{u}{n,\half}\bigr)\\\nonumber
				&\quad :=I_1 + I_2 + I_3 + I_4 + I_5.
			\end{align}
			
			Using \eqref{20220608_2} and \eqref{assump:LipsF}, the first term $I_1$ could be bounded by
			\begin{align*}
				\mE[I_1]& \leq \tau\exc{\normbf{\bF^{n,\half}-\bF_h^{n,\half}}}+\frac{\tau}{4}\exc{\normbf{\Ph\bfErr{v}{n+1}}}\notag\\
				&\quad \leq \tau C_B\exc{\normdel{\bfErr{u}{n,\half}}+\normbf{\bfErr{u}{n,\half}}}+\frac{\tau}{4}\exc{\normbf{\Ph\bfErr{v}{n+1}}}\notag\\
				&\quad \leq \tau CC_B\exc{\normdel{\Rh\bfErr{u}{n,\half}}+\normbf{\Rh\bfErr{u}{n,\half}}}\notag\\
				&\qquad +\tau CC_B\mathbb{E}\Bigl[\normdel{\bu^{n,\half}-\Rh\bu^{n,\half}}\notag\\
				&\qquad +\normbf{\bu^{n,\half}-\Rh\bu^{n,\half}}\Bigr]+\frac{\tau}{4}\exc{\normbf{\Ph\bfErr{v}{n+1}}}\notag\\
				&\quad \leq \tau CC_B\exc{\normdel{\Rh\bfErr{u}{n,\half}}}\notag\\
				&\qquad + \tau C_BC_{\Rh}h^2\exc{\left\lVert\bu^{n,\half}\right\rVert^2_{\bH^2}}+\frac{\tau}{4}\exc{\normbf{\Ph\bfErr{v}{n+1}}}.\notag
			\end{align*}	
			
			To  bound $I_2$,  using the martingale property of the It\^o's integral yields 
			\begin{align*}
				\mE[I_2]&\leq C\tau\exc{\normbf{\bG[\bu^n]-\bG[\bu^n_h]}}+\frac{1}{16}\exc{\normbf{\Ph\bfErr{v}{n+1}-\Ph\bfErr{v}{n}}}\notag\\
				&\leq CC_B\tau\exc{\normdel{\Rh\bfErr{u}{n}}}\notag\\
				&\qquad\quad+ C_BC_{\Ph}\tau h^4\exc{\left\lVert\bu^n\right\rVert^2_{\bH^2}}+\frac{1}{16}\exc{\normbf{\Ph\bfErr{v}{n+1}-\Ph\bfErr{v}{n}}}.\notag
			\end{align*}
			By using the same approach  for bounding  $I_2$ with help of  the Ladyzhenskaya's inequality, and \eqref{assump:GradientF}, and \eqref{assump:Fuu}, we obtain
			\begin{align}\label{eq_4.16}
				\mE[I_3] &= \mE\left[\bigl((D_{\bu}\bG[\bu^n]\bv^n - D_{\bu}\bG[\bu_h^n]\bv^n_h)\, \widehat{\Delta}W_n,\Ph\bfErr{v}{n+1} - \Ph\bfErr{v}{n}\bigr)\right]\\\nonumber
				&= \mE\left[\bigl((D_{\bu}\bG[\bu^n] - D_{\bu}\bG[\bu_h^n])\bv^n\, \widehat{\Delta}W_n,\Ph\bfErr{v}{n+1} - \Ph\bfErr{v}{n}\bigr)\right]\\\nonumber
				&\qquad+\mE\left[\bigl(D_{\bu}\bG[\bu_h^n](\bv^n - \bv^n_h)\, \widehat{\Delta}W_n,\Ph\bfErr{v}{n+1} - \Ph\bfErr{v}{n}\bigr)\right]\\\nonumber
				&:= I_{3,1} + I_{3,2}.
			\end{align}
			
			To bound  $I_{3,2}$, we use  the assumption \eqref{assump:GradientF}, \eqref{eqn:prj-error-estimate-ph} and Lemma \ref{lemma:time-semiform-stab-uh2-vh1} to get 
			\begin{align}
				I_{3,2} &\leq CC_A\tau^3 \mE[\|\bfErr{\bv}{n}\|^2_{\bL^2}] + \frac{1}{16}\mE[\|\Ph\bfErr{v}{n+1} - \Ph\bfErr{v}{n}\|^2_{\bL^2}]\\\nonumber
				&\leq CC_A\tau^3 \mE[\|\Ph\bfErr{\bv}{n}\|^2_{\bL^2}] + CC_A\tau^3h^2\mE[\|\bv^n\|^2_{\bH^1}]\\\nonumber
				&\qquad+ \frac{1}{16}\mE[\|\Ph\bfErr{v}{n+1} - \Ph\bfErr{v}{n}\|^2_{\bL^2}]\\\nonumber
				&\leq CC_A\tau^3 \mE[\|\Ph\bfErr{\bv}{n}\|^2_{\bL^2}] + CC_AC_3C_{\Ph}\tau^3h^2+ \frac{1}{16}\mE[\|\Ph\bfErr{v}{n+1} - \Ph\bfErr{v}{n}\|^2_{\bL^2}]. \nonumber
			\end{align}
			To estimate $I_{3,1}$,  adding and subtracting $\bv(t_n)$ we obtain
			\begin{align}
				I_{3,1} &=  \mE\left[\bigl((D_{\bu}\bG[\bu^n] - D_{\bu}\bG[\bu_h^n])(\bv^n - \bv(t_n))\, \widehat{\Delta}W_n,\Ph\bfErr{v}{n+1} - \Ph\bfErr{v}{n}\bigr)\right]\\\nonumber
				&\qquad +  \mE\left[\bigl((D_{\bu}\bG[\bu^n] - D_{\bu}\bG[\bu_h^n])\bv(t_n)\, \widehat{\Delta}W_n,\Ph\bfErr{v}{n+1} - \Ph\bfErr{v}{n}\bigr)\right]\\\nonumber
				&:= {I}_{3,1}^1 + {I}_{3,1}^2.
			\end{align}
			
			It follows from the assumption \eqref{assump:GradientF} and Theorem \ref{theorem3.4} that
			\begin{align}
				I_{3,1}^1 &\leq \frac{1}{16}\mE[\|\Ph\bfErr{v}{n+1} - \Ph\bfErr{v}{n}\|^2_{\bL^2}] + CC_A\tau^3\mE[\|\bv(t_n) - \bv^n\|^2_{\bL^2}]\\\nonumber
				&\leq \frac{1}{16}\mE[\|\Ph\bfErr{v}{n+1} - \Ph\bfErr{v}{n}\|^2_{\bL^2}] + CC_A\tilde{C}_2\tau^5.
			\end{align}
			Using the assumption \eqref{assump:Fuu}, the Ladyzhenskaya's inequality, Lemmas \ref{lem:highmomentstability},  \ref{lemma3.1.2}  and   \ref{lemma:fem-form-stab}, we obtain
			\begin{align}
				I_{3,1}^2 &\leq \frac{1}{16}\mE[\|\Ph\bfErr{v}{n+1} - \Ph\bfErr{v}{n}\|^2_{\bL^2}] + CC_A\tau^3\mE[\| \bfErr{u}{n}\|^2_{\bL^4}\|\bv(t_n)\|^2_{\bL^4}]\\\nonumber
				&\leq \frac{1}{16}\mE[\|\Ph\bfErr{v}{n+1} - \Ph\bfErr{v}{n}\|^2_{\bL^2}] + CC_A\tau^3\mE[\|\bfErr{u}{n}\|_{\bL^2}\|\nabla \bfErr{u}{n}\|_{\bL^2}\|\bv(t_n)\|^2_{\bH^1}]\\\nonumber
				&\leq \frac{1}{16}\mE[\|\Ph\bfErr{v}{n+1} - \Ph\bfErr{v}{n}\|^2_{\bL^2}] \\\nonumber
				&\qquad+ CC_A\tau^5\mE[\|\bfErr{u}{n}\|^2_{\bL^2}\|\bv(t_n)\|^4_{\bH^1}] + CC_A\tau\mE[\|\nabla \bfErr{u}{n}\|^2_{\bL^2}]\\\nonumber
				&\leq \frac{1}{16}\mE[\|\Ph\bfErr{v}{n+1} - \Ph\bfErr{v}{n}\|^2_{\bL^2}] \\\nonumber
				&\qquad+ CC_A\tau^5\mE[\|\bfErr{u}{n}\|^4_{\bL^2} + \|\bv(t_n)\|^8_{\bH^1}] + CC_A\tau\mE[\|\nabla \bfErr{u}{n}\|^2_{\bL^2}]\\\nonumber
				&\leq \frac{1}{16}\mE[\|\Ph\bfErr{v}{n+1} - \Ph\bfErr{v}{n}\|^2_{\bL^2}] \\\nonumber
				&\qquad+ CC_A(C_{2,2} + {C}_{3,2} + C_{s,8})\tau^5+ CC_A\tau\mE[\|\nabla \bfErr{u}{n}\|^2_{\bL^2}] \\\nonumber
				&\leq \frac{1}{16}\mE[\|\Ph\bfErr{v}{n+1} - \Ph\bfErr{v}{n}\|^2_{\bL^2}] \\\nonumber
				&\qquad+ CC_A\tau^5\mE[\|\bfErr{u}{n}\|^4_{\bL^2} + \|\bv(t_n)\|^8_{\bH^1}] + CC_A\tau\mE[\|\nabla \bfErr{u}{n}\|^2_{\bL^2}]\\\nonumber
				&\leq \frac{1}{16}\mE[\|\Ph\bfErr{v}{n+1} - \Ph\bfErr{v}{n}\|^2_{\bL^2}] \\\nonumber
				&\qquad+ CC_A(C_{2,2} + {C}_{3,2} + \tilde{K}_{1,4})\tau^5+ CC_A\tau\mE[\|\varepsilon (\bfErr{u}{n})\|^2_{\bL^2}]\\\nonumber
				&\leq \frac{1}{16}\mE[\|\Ph\bfErr{v}{n+1} - \Ph\bfErr{v}{n}\|^2_{\bL^2}] + CC_A(C_{2,2} + {C}_{3,2} + \tilde{K}_{1,4})\tau^5\\\nonumber
				&\qquad+ CC_A\tau\mE[\|\varepsilon (\Rh\bfErr{u}{n})\|^2_{\bL^2}] + CC_AC_{\Rh}\tau h^2\mE[\|\bu^n\|^2_{\bH^2}],
			\end{align}
			where the second to last inequality  is obtained by using the Korn's inequality, while the last inequality holds because of \eqref{20220608_3}.
			
			Using \eqref{eqn:prf-fem-form-order-h1-1} with $\bphi_h=\opLh\Rh\bfErr{u}{n,\half}$ and \eqref{assump:LipsG},  we obtain
			\begin{align*}
				\mE[I_4] &= -\tau\exc{\bigl(\bfErr{v}{n+1}, \opLh\Rh\bfErr{u}{n,\half}\bigr)} + \exc{\bigl((\bG[\bu^n] - \bG(\bu^n_h))\widehat{\Delta}W_n,\opLh\Rh\bfErr{u}{n,\half}\bigr)}\\\nonumber
				&\leq C\tau\exc{\normdel{\bfErr{v}{n+1}}} \\\nonumber
				&\qquad+ C\tau \exc{\normdel{\Rh\bfErr{u}{n,\half}}}  \\\nonumber
				&\qquad + CC_B\tau\exc{\normdel{\Rh\bfErr{u}{n}}}.
			\end{align*}
			Similarly,  
			\begin{align*}
				\exc{I_5} &\leq CC_B\tau \exc{\normdel{\Rh\bfErr{u}{n}}} \\\nonumber
				&\qquad+ C\tau \exc{\normdel{\Rh\bfErr{u}{n,\half}}}.
			\end{align*}
			
			Combining all the estimates from $I_1,\cdots, I_5$ into the right-hand side of \eqref{eqn:prf-fem-form-order-h1-4}, we get
			\begin{align}\label{eqn:prf-fem-form-order-h1-8}
				&\frac12\exc{\|\Ph\bfErr{v}{n+1}\|_{{\bf L}^2}^2-\|\Ph\bfErr{v}{n}\|_{{\bf L}^2}^2+ \normdel{\Rh\bfErr{u}{n+1}}} \\\nonumber&
				\qquad -\frac{1}{2}\exc{\normdel{\Rh\bfErr{u}{n-1}}}\\\nonumber&\qquad+\frac14\exc{\normbf{\Ph\bfErr{v}{n+1}-\Ph\bfErr{v}{n}}}\\\nonumber
				&\leq C(C_B+1)\tau \exc{J(\Rh\bfErr{u}{n,\half},\Ph\bfErr{v}{n+1})} + C(C_B+C_A)\tau\exc{J(\Rh\bfErr{u}{n},\Ph\bfErr{v}{n})}\\\nonumber
				&\qquad +  \tau C_BC_{\Rh}h^2\exc{\left\lVert\bu^{n,\half}\right\rVert^2_{\bH^2}} + C_BC_{\Ph}\tau h^4\exc{\left\lVert\bu^n\right\rVert^2_{\bH^2}} \\\nonumber
				&\qquad+CC_A(\tilde{C}_2+C_{2,2} + C_{3,2} + \tilde{K}_{1,4})\tau^5 \\\nonumber
				&\qquad+ CC_AC_3C_{\Ph}\tau^3h^2 + CC_AC_{\Rh}\tau h^2\mE[\|\bu^n\|^2_{\bH^2}].
			\end{align}
			
			Applying the summation operator $\sum_{n=1}^{\ell}$ for any $1\leq \ell <N$ to \eqref{eqn:prf-fem-form-order-h1-8} and then using Lemma \ref{lemma:time-semiform-stab-uh2-vh1}, we obtain
			\begin{align}\label{eq4.16}
				&\frac12\exc{J(\Rh\bfErr{u}{\ell+1},\Ph\bfErr{v}{\ell+1})} \\\nonumber
				&\leq \frac12\exc{J(\Rh\bfErr{u}{1}, \Ph\bfErr{v}{1})} + \frac12\exc{J(\Rh\bfErr{u}{0}, \Ph\bfErr{v}{0})}\\\nonumber
				&\qquad+ C(C_BC_{\Rh} + C_AC_{\Ph} + C_AC_{\Rh} + h^2C_BC_{\Ph})C_3h^2 \\\nonumber
				&\qquad+CC_A(\tilde{C}_2+C_{2,2} + C_{3,2} + \tilde{K}_{1,4})\tau^4 \\\nonumber
				&\qquad + C(C_A+C_B+1)\tau\sum_{n=1}^{\ell}\mE\left[J(\Rh\bfErr{\bu}{n,\half},\Ph\bfErr{\bv}{n+1}) + J(\Rh\bfErr{\bu}{n},\Ph\bfErr{\bv}{n})\right].
			\end{align}
			
			It follows from  the discrete Gronwall's inequality that 
			\begin{align}\label{eqn:prf-fem-form-order-h1-9}
				\exc{J(\Rh\bfErr{u}{\ell+1}, \Ph\bfErr{v}{\ell+1})} \leq \tilde{C}_3(h^2 + \tau^4),
			\end{align}
			where 
			\begin{align*}
				\tilde{C}_3 &= \max\{(C_BC_{\Rh} + C_AC_{\Ph} + C_AC_{\Rh} + h^2C_BC_{\Ph})C_3, \\
				&\hskip 1.2in C_A(\tilde{C}_2+C_{2,2} + C_{3,2} + \tilde{K}_{1,4})\}e^{C(C_A+C_B+1)T}.
			\end{align*}
			
			The proof is complete after combining \eqref{eqn:prf-fem-form-order-h1-9}, \eqref{eqn:prj-error-estimate-ph}, and \eqref{20220608_3}.
		\end{proof}
		
		Theorem \ref{theo:fem-form-order-h1} shows that the linear finite element method converges with rate $O(h)$ in the $\bH^1$-norm.  The next theorem shows that the convergence rate becomes $O(h^2)$ in the $\bL^2$-norm .
		
		\begin{theorem}\label{theo:fem-form-order-l2}
			Let $\{(\bu^n, \bv^n)\}$ and $\{(\bu^n_h, \bv^n_h)\}$ be the solutions generated by Algorithm 1 and Algorithm 2, respectively. Let $(\bu^0_h, \bv^0_h)=(\bu^0, \bv^0)$, and $(\bu^1_h, \bv^1_h)=(\bu^1, \bv^1)$. Suppose that $(\bu^0, \bv^0)\in L^4(\Omega; (\bH^2\cap \bH_0^1)\times\bH^1_0)$. Suppose that the assumptions \eqref{assump:F(0)}--\eqref{assump:LipsG} satisfy. Then, there holds
			\begin{equation}\label{eqn:fem-form-order-h2}
				\max_{1 \leq n \leq N}\exc{\normbf{\bu^n - \bu^n_h}}\leq \tilde{C}_4\,(\tau^3+h^4),
			\end{equation}
			where $\tilde{C}_4 = e^{CC_AT}\max\left\{CC_AC_{\Ph}C_3, CC_A[(\overline{C}_{2,2} + \overline{C}_{3,2})\tilde{K}_{1,4} + \widetilde{C}_2]\tau\right\}$.
		\end{theorem}
		
		\begin{proof}\label{prf:fem-form-order-l2}
			Let $d_t\bu^l_h:=\tau^{-1}(\bu^l_h-\bu^{l-1}_h)$. Then, from \eqref{FEM1}, we have
			\begin{align}\label{eq4.20}
				\bigl(\bv_h^{n+1},\bphi_h\bigr) = \bigl(d_t\bu^{l+1}_h,\bphi_h\bigr) + \frac{1}{\tau}\bigl(\bG[\bu_h^l]\widehat{\Delta}W_l,\bphi_h\bigr).
			\end{align}
			Using \eqref{eq4.20} and \eqref{FEM2} we get
			\begin{align}\label{eqn:prf-fem-form-order-l2-1}
				&\big(d_t\bu^{l+1}_h-d_t\bu^l_h, \bphi_h\big)+\tau\lambda\big(\Div(\bu^{l,\half}_h), \Div(\bphi_h)\big) + \tau\mu\big(\varepsilon(\bu^{l,\half}_h), \varepsilon(\bphi_h)\big)\\\nonumber
				&\qquad = \tau\bigl(\bF_h^{n,\frac12}, \bphi_h\bigr) 
				+ \bigl(\bG[\bu_h^n]\overline{\Delta} W_n, \bphi_h \bigr)  +  \bigl(D_{\bu}\bG[\bu_h^n]\bv^n_h\widehat{\Delta}W_n, \bphi_h \bigr)\\\nonumber
				&\qquad \qquad-\frac{1}{\tau}\bigl(\bG[\bu^l_h]\widehat{\Delta}W_l,\bphi_h\bigr)-\frac{1}{\tau}\bigl(\bG[\bu^{l-1}_h]\widehat{\Delta}W_{l-1},\bphi_h\bigr).
			\end{align}
			Setting $\bar{\bu}^{n}_h=\sum\limits_{l = 1}^{n}\bu^l_h$ and taking the summation over $l$ in  \eqref{eqn:prf-fem-form-order-l2-1} yield
			\begin{align}\label{eqn:prf-fem-form-order-l2-2}
				&	\bigl(\bu_h^{n+1} - \bu_h^n, \bphi_h\bigr) + \bigl(\bG[\bu^n_h]\widehat{\Delta} W_n, \bphi\bigr) - \tau^2\bigl(\opL\overline{\bu}_h^{n,\frac12},\bphi_h\bigr)\\\nonumber
				&\qquad  = \left(\bu_h^1 - \bu_h^0,\bphi_h\right) + \bigl(\bG[\bu_h^0]\widehat{\Delta} W_0 + \tau^2\opL \bu_h^{\frac12},\bphi_h\bigr) + \tau^2\sum_{l=1}^n \bigl(\bF_h^{l,\frac12}, \bphi_h\bigr)\\\nonumber 
				&\qquad \qquad + \tau\sum_{l=1}^n \bigl(\bG[\bu_h^l],\bphi_h\bigr)\overline{\Delta} W_l +\tau\sum_{l=1}^n\bigl(D_{\bu}\bG[\bu_h^l]\bv_h^l\widehat{\Delta}W_l, \bphi_h\bigr).
			\end{align}
			
			Subtracting \eqref{eqn:prf-fem-form-order-l2-2} from \eqref{eq3.24} leads to
			\begin{align}\label{eqn:prf-fem-form-order-l2-3}
				&\big(\bfErr{u}{n+1}-\bfErr{u}{n}, \bphi_h\big)+\tau^2\left[\innerprd{\bar{\bf E}_{\bf u}^{n+\half}}{\bphi_h}\right]\\\nonumber
				&\qquad = \tau\left(\sum\limits_{l=1}^{n}(\bG[\bu^l]-\bG[\bu^l_h])\overline{\Delta} W_l, \bphi_h\right) + \tau^2\left(\sum\limits_{l = 1}^{n}(\bF^{l,\half}-\bF_h^{l,\half}), \bphi_h\right)\\\nonumber
				&\qquad \quad+\tau\sum_{l=1}^n\bigl((D_{\bu}\bG[\bu^l]\bv^l - D_{\bu}\bG[\bu_h^l]\bv_h^l)\widehat{\Delta}W_l, \bphi_h\bigr) \\\nonumber
				&\qquad \quad- \bigl((\bG[\bu^n] - \bG[\bu^n_h])\widehat{\Delta} W_n, \bphi_h\bigr).
			\end{align}
			
			Choosing $\bphi_h=\Ph\bfErr{u}{n+\half}$ in \eqref{eqn:prf-fem-form-order-l2-3}, the first term on the left-hand side of \eqref{eqn:prf-fem-form-order-l2-3} can be written as
			\begin{align}\label{20220613_4}
				\big(\bfErr{u}{n+1}-\bfErr{u}{n}, \Ph\bfErr{u}{n+\half}\big)&=\big(\Ph\bfErr{u}{n+1}-\Ph\bfErr{u}{n}, \Ph\bfErr{u}{n+\half}\big)\\\nonumber
				&=\frac12\|\Ph\bfErr{u}{n+1}\|_{{\bf L}^2}^2-\frac12\|\Ph\bfErr{u}{n}\|_{{\bf L}^2}^2.
			\end{align}
			Moreover, 
			\begin{align}\label{20220615_1}
				&\innerprd{\bar{\bf E}_{\bf u}^{n,\half}}{\Ph\bfErr{u}{n+\half}}\\
				&\qquad =\innerprd{\Rh\bar{\bf E}_{\bf u}^{n,\half}}{\Ph\bfErr{u}{n+\half}}\notag\\
				&\qquad =-\big(\opLh\Rh\bar{\bf E}_{\bu}^{n,\half}, \Ph\bfErr{u}{n+\half}\big)\notag\\
				&\qquad=-\big(\opLh\Rh\bar{\bf E}_{\bu}^{n,\half}, \bfErr{u}{n+\half}\big)\notag\\
				&\qquad =-\half\big(\opLh\Rh\bar{\bf E}_{\bu}^{n,\half}, \Rh(\bar{\bf E}_{\bu}^{n+1}-\bar{\bf E}_{\bu}^{n-1})\big)\notag\\\nonumber
				&\qquad = \frac{\lambda}{4}\left(\normbf{\Div \Rh\bar{\bf E}_{\bu}^{n+1}} - \normbf{\Div \Rh\bar{\bf E}_{\bu}^{n-1}}\right) \\\nonumber
				&\qquad\qquad + \frac{\mu}{4}\left(\normbf{\veps(\Rh\bar{\bf E}_{\bu}^{n+1})} - \normbf{\veps(\Rh\bar{\bf E}_{\bu}^{n-1})}\right).
			\end{align}
			
			Define  the following energy functional $\cal{E}(\bfErr{u}{n})$:
			\begin{equation}
				\cal{E}(\bfErr{u}{n}):=\half\normbf{\Ph\bfErr{u}{n}}+\frac{\tau^2}{4}\left(\normdel{\Rh\bar{\bf E}_{u}^{n}}\right).
			\end{equation}
			Choosing $\bphi_h=\Ph\bfErr{u}{n+\half}$ in \eqref{eqn:prf-fem-form-order-l2-3} and using \eqref{20220613_4}--\eqref{20220615_1} yield
			\begin{align}\label{eqn:prf-fem-form-order-l2-4}
				&\cal{E}(\bfErr{u}{n+1})-\cal{E}({\bfErr{u}{n}})+  \frac{\tau^2\lambda}{4}\left(\normbf{\Div \Rh\bar{\bf E}_{\bu}^{n}} - \normbf{\Div \Rh\bar{\bf E}_{\bu}^{n-1}}\right) \\\nonumber
				&\qquad+ \frac{\tau^2\mu}{4}\left(\normbf{\veps(\Rh\bar{\bf E}_{\bu}^{n})} - \normbf{\veps(\Rh\bar{\bf E}_{\bu}^{n-1})}\right) \\\nonumber
				&= \tau^2\left(\sum\limits_{l = 1}^{n}(\bF^{l,\half}-\bF_h^{l,\half}), \Ph\bfErr{u}{n+\half}\right) \\\nonumber&\qquad+\tau\left(\sum\limits_{l=1}^{n}(\bG[\bu^l]-\bG[\bu^l_h])\overline{\Delta} W_l, \Ph\bfErr{u}{n+\half}\right) \\\nonumber
				&\qquad+\tau\sum_{l=1}^n\bigl((D_{\bu}\bG[\bu^l]\bv^l - D_{\bu}\bG[\bu_h^l]\bv_h^l)\widehat{\Delta}W_l, \Ph\bfErr{u}{n+\half}\bigr) \\\nonumber
				&\qquad- \bigl((\bG[\bu^n] - \bG[\bu^n_h])\widehat{\Delta} W_n, \Ph \bfErr{u}{n+\half}\bigr)\\\nonumber
				&:= I_1 + I_2 + I_3 + I_4.
			\end{align}
			
		To bound $I_1$, by the assumption \eqref{assump:GradientF} for $\bF$,  we have
			\begin{align}\label{eqn:prf-fem-form-order-l2-6}
				\exc{I_1} &=\exc{\tau^2\bigl(\sum_{l = 1}^{n}(\bF^{l,\half}-\bF_h^{l,\half}), \Ph\bfErr{u}{n+\half}\bigr)}\\\nonumber
				&\leq C\tau^3\exc{\normbf{\sum_{l = 1}^{n}(\bF^{l,\half}-\bF_h^{l,\half})}} + \tau\exc{\normbf{\Ph\bfErr{u}{n+\half}}}\\\nonumber
				&\leq C\tau^2\exc{\sum_{l = 1}^{n}\normbf{\bF^{l,\half}-\bF_h^{l,\half}}} + \tau\exc{\normbf{\Ph\bfErr{u}{n+\half}}}\\\nonumber
				&\leq CC_A\tau^2\exc{\sum_{l = 1}^{n}\left(\normbf{\bfErr{u}{l+1}} + \normbf{\bfErr{u}{l-1}}\right)} + \tau\exc{\normbf{\Ph\bfErr{u}{n+\half}}}\\\nonumber
				&\leq CC_A\tau^2\exc{\sum_{l = 1}^{n}\left(\normbf{\Ph\bfErr{u}{l+1}} + \normbf{\Ph\bfErr{u}{l-1}}\right)} \\\nonumber
				&\qquad+CC_AC_{\Ph}\tau h^4\max\limits_{1\le l\le N}\exc{\left\lVert{\bu^{l}}\right\rVert^2_{\bH^2}}+ \tau\exc{\normbf{\Ph\bfErr{u}{n+\frac12}}}\notag.
			\end{align}
			
			Similarly, by the assumption \eqref{assump:GradientF} and the It\^o's isometry, $I_2$ can be bounded by 
			\begin{align}\label{eqn:prf-fem-form-order-l2-7}
				\exc{I_2}&=\exc{\tau\left(\sum_{l = 1}^{n}\big(\bG[\bu^l]-\bG[\bu^l_h]\big)\commentone{\overline{\Delta} W_{l}}, \Ph\bfErr{u}{n+\half}\right)}\\\nonumber
				&\leq C\tau\exc{\normbf{\sum_{l = 1}^{n}\big(\bG[\bu^l]-\bG[\bu^l_h]\big)\commentone{\overline{\Delta} W_{l}}}} + \tau\exc{\normbf{\Ph \bfErr{u}{n+\half}}} \\\nonumber
				&= C\tau\exc{\tau\sum_{l = 1}^n\normbf{\bG[\bu^l]-\bG[\bu^l_h]}} + \tau\exc{\normbf{\Ph\bfErr{u}{n+\half}}}\notag\\
				&\leq CC_A\tau^2\sum_{l = 1}^n\exc{\normbf{\Ph\bfErr{u}{l}}}\notag\\
				&\quad+CC_{\Ph}C_A\tau h^4\max\limits_{1\le l\le n}\exc{\left\lVert{\bu^{l}}\right\rVert^2_{\bH^2}}+ \tau\exc{\normbf{\Ph\bfErr{u}{n+\half}}}\notag.
			\end{align}
			
			Next, to estimate $I_3$, inserting four terms $\pm \left(D_{\bu}\bG[\bu^l] - D_{\bu}\bG[\bu_h^l]\bv(t_l)\widehat{\Delta}W_{l}\right)$ and 
			$\pm D_{\bu}\bG[\bu_h^l]\bv^l\widehat{\Delta}W_l$, and rearranging terms, we obtain
			\begin{align}\label{eq4.31}
				\exc{I_3} &= \tau\mE\left[\left(\sum_{l=1}^{n}(D_{\bu}\bG[\bu^l] - D_{\bu}\bG[\bu^l_h])\bv(t_l)\widehat{\Delta}W_l, \Ph \bfErr{\bu}{n+\half}\right)\right]\\\nonumber
				&\qquad+  \tau\mE\left[\left(\sum_{l=1}^{n}(D_{\bu}\bG[\bu^l] - D_{\bu}\bG[\bu^l_h])\bferr{\bv}{l}\widehat{\Delta}W_l, \Ph \bfErr{\bu}{n+\half}\right)\right]\\\nonumber
				&\qquad +  \tau\mE\left[\left(\sum_{l=1}^{n}D_{\bu}\bG[\bu^l]\bfErr{\bv}{l}\widehat{\Delta}W_l, \Ph \bfErr{\bu}{n+\half}\right)\right]\\\nonumber
				&:= I_{3,1} + I_{3,2} + I_{3,3},
			\end{align}
			where $\bferr{\bv}{l} = \bv(t_l) - \bv^l$ is the semi-discrete error from Theorem \ref{theorem3.4}.
			
			Now using the It\^o's isometry, Theorem \ref{theorem3.4} and the assumption \eqref{assump:GradientF}, the second term on the right side of \eqref{eq4.31} can be controlled as 
			\begin{align*}
				I_{3,2} &\leq C\tau\mE\left[\left\|\sum_{l=1}^{n}(D_{\bu}\bG[\bu^l] - D_{\bu}[\bu^l_h])\bferr{\bv}{l}\widehat{\Delta}W_l\right\|^2_{\bL^2}\right] + \tau\mE\left[\|\Ph\bfErr{\bu}{n+\frac12}\|^2_{\bL^2}\right]\\\nonumber
				&\leq CC_A\tau\mE\left[\tau^3\sum_{l=1}^{n}\left\|\bferr{\bv}{l}\right\|^2_{\bL^2}\right] + \tau\mE\left[\|\Ph\bfErr{\bu}{n+\frac12}\|^2_{\bL^2}\right]\\\nonumber
				&\leq CC_A\widetilde{C}_2\tau^5+ \tau\mE\left[\|\Ph\bfErr{\bu}{n+\frac12}\|^2_{\bL^2}\right],
			\end{align*}
			and the third term on the right-hand side of \eqref{eq4.31} can be bounded by using \eqref{eqn:prf-fem-form-order-h1-1} with $\bPhi_h = \bfErr{v}{l}$ and \eqref{eqn:prj-error-estimate-ph} as follows. Noticing that 
			\begin{align*}
				\tau^2\exc{\normbf{\bfErr{v}{l}} }&\leq \exc{\normbf{\bfErr{u}{l+1} - \bfErr{u}{l}} + C_A\tau^3\normbf{\bfErr{u}{l}}}\\\nonumber
				&\leq \exc{\normbf{\Ph\bfErr{u}{l+1} - \Ph\bfErr{u}{l}} + C_A\tau^3\normbf{\Ph\bfErr{u}{l}}} \\\nonumber
				&\qquad+ (1+ \tau^3C_A)C_{\Ph}h^4\max_{1\leq l\leq N}\exc{\|\bu^{l}\|^2_{\bH^2}},
			\end{align*}
			which and the assumption \eqref{assump:GradientF} infer
			\begin{align*}
				I_{3,3} &\leq CC_A\mE\left[\tau^2\sum_{l=1}^{n}\tau^2\|\bfErr{\bv}{l}\|^2_{\bL^2}\right] + \tau\mE\left[\|\Ph\bfErr{\bu}{n+\frac12}\|^2_{\bL^2}\right]\\\nonumber
				&\leq CC_A\tau^2\sum_{l=1}^{n}\mE\left[\|\Ph\bfErr{\bu}{l+1} - \Ph\bfErr{\bu}{l}\|^2_{\bL^2}\right] +CC_A\tau^5\sum_{l=1}^{n}\mE\left[\| \Ph\bfErr{\bu}{l}\|^2_{\bL^2}\right] \\\nonumber
				&\qquad+ CC_A\tau(1+\tau^3C_A)C_{\Ph}h^4\max_{1\leq l\leq N}\exc{\|\bu^{l}\|^2_{\bH^2}}  + \tau\mE\left[\|\Ph\bfErr{\bu}{n+\frac12}\|^2_{\bL^2}\right].
			\end{align*}
			
			Furthermore,  we bound the first term on the right-hand side of \eqref{eq4.31} by using the assumption \eqref{assump:Fuu} and the Ladyzhenskaya's inequality as follows:
			\begin{align*}
				I_{3,1} &\leq \tau\mE\left[\|\Ph\bfErr{\bu}{n+\half}\|^2_{\bL^2}\right] + C\tau \mE\left[\left\|\sum_{l=1}^{n}(D_{\bu}\bG[\bu^l] - D_{\bu}\bG[\bu^l_h])\bv(t_l)\widehat{\Delta}W_l\right\|^2_{\bL^2}\right]\\\nonumber
				&\leq \tau\mE\left[\|\Ph\bfErr{\bu}{n+\half}\|^2_{\bL^2}\right] + C \tau^4\sum_{l=1}^{n}\mE\left[\left\|(D_{\bu}\bG[\bu^l] - D_{\bu}\bG[\bu^l_h])\bv(t_l)\right\|^2_{\bL^2}\right]\\\nonumber
				&\leq \tau\mE\left[\|\Ph\bfErr{\bu}{n+\half}\|^2_{\bL^2}\right] + CC_A \tau^4\sum_{l=1}^{n}\mE\left[\left\|\bfErr{\bu}{l}\bv(t_l)\right\|^2_{\bL^2}\right]\\\nonumber
				&\leq \tau\mE\left[\|\Ph\bfErr{\bu}{n+\half}\|^2_{\bL^2}\right] + CC_A \tau^4\sum_{l=1}^{n}\mE\left[\left\|\bfErr{\bu}{l}\|_{\bL^2}\|\nabla\bfErr{\bu}{l}\|_{\bL^2}\|\bv(t_l)\right\|^2_{\bH^1}\right].
			\end{align*}
			Using the above estimate and the high moment stability estimates from Lemmas \ref{lem:highmomentstability},  \ref{lemma3.4}, and  \ref{lemma:fem-form-stab_4.3}, we obtain
			\begin{align*}
				I_{3,1} &\leq \tau\mE\left[\|\Ph\bfErr{\bu}{n+\half}\|^2_{\bL^2}\right] + CC_A \tau^2\sum_{l=1}^{n}\mE\left[\|\bfErr{\bu}{l}\|^2_{\bL^2}\right] \\\nonumber
				&\qquad+ CC_A\tau^6\sum_{l=1}^{n}\mE\left[\|\nabla\bfErr{\bu}{l}\|^2_{\bL^2}\|\bv(t_l)\|^4_{\bH^1}\right]\\\nonumber
				&\leq \tau\mE\left[\|\Ph\bfErr{\bu}{n+\half}\|^2_{\bL^2}\right] + CC_A \tau^2\sum_{l=1}^{n}\mE\left[\|\bfErr{\bu}{l}\|^2_{\bL^2}\right] \\\nonumber
				&\qquad+ CC_A\tau^6\sum_{l=1}^{n}\left(\mE\left[\|\nabla\bfErr{\bu}{l}\|^4_{\bL^2}\right]\right)^{\half}\left(\mE\left[\|\bv(t_l)\|^8_{\bH^1}\right]\right)^{\half}\\\nonumber
				&\leq \tau\mE\left[\|\Ph\bfErr{\bu}{n+\half}\|^2_{\bL^2}\right] + CC_A \tau^2\sum_{l=1}^{n}\mE\left[\|\bfErr{\bu}{l}\|^2_{\bL^2}\right] \\\nonumber
				&\qquad+ CC_A\left(\overline{C}_{2,2}+\overline{C}_{3,2}\right)\tilde{K}_{1,4}\tau^5.
			\end{align*}
			
			Finally,  using \eqref{assump:GradientF},  $I_4$ can be controlled by 
			\begin{align*}
				\exc{I_4} &\leq CC_A\tau^2\exc{\normbf{\bfErr{u}{n}}} + \tau\exc{\normbf{\Ph\bfErr{u}{n+\half}}}\\\nonumber
				&\leq CC_A\tau^2\exc{\normbf{\Ph\bfErr{u}{n}}} +CC_AC_{\Ph}\tau^2h^4\exc{\|\bu^n\|^2_{\bH^2}} +  \tau\exc{\normbf{\Ph\bfErr{u}{n+\half}}}.
			\end{align*}
			
			Combining all the estimates of $I_1,\cdots, I_4$ into \eqref{eqn:prf-fem-form-order-l2-4} we obtain
			\begin{align}\label{eqn:prf-fem-form-order-l2-8}
				&\exc{\cal{E}(\bfErr{u}{n+1})}-\exc{\cal{E}({\bfErr{u}{n}})}+  \frac{\tau^2\lambda}{4}\exc{\normbf{\Div \Rh\bar{\bf E}_{\bu}^{n}} - \normbf{\Div \Rh\bar{\bf E}_{\bu}^{n-1}}} \\\nonumber
				&\qquad+ \frac{\tau^2\mu}{4}\exc{\normbf{\veps(\Rh\bar{\bf E}_{\bu}^{n})} - \normbf{\veps(\Rh\bar{\bf E}_{\bu}^{n-1})}} \\\nonumber
				&\leq CC_A\tau\exc{\tau\sum_{l = 1}^{n}\left(\normbf{\Ph\bfErr{u}{l+1}} +\normbf{\Ph\bfErr{u}{l}} + \normbf{\Ph\bfErr{u}{l-1}}\right)} \\\nonumber
				&\qquad+C \tau\exc{\normbf{\Ph\bfErr{u}{n+\frac12}}} + CC_A\tau^2\exc{\normbf{\Ph\bfErr{u}{n}}}\\\nonumber
				&\qquad+CC_AC_{\Ph}\tau h^4\max\limits_{1\le l\le N}\exc{\left\lVert{\bu^{l}}\right\rVert^2_{\bH^2}} + CC_A[(\overline{C}_{2,2} + \overline{C}_{3,2})\tilde{K}_{1,4} + \widetilde{C}_2]\tau^5\\\nonumber
				&\leq CC_A\tau\max_{1\leq l\leq n}\exc{\normbf{\Ph\bfErr{u}{l+1}} +\Ph\bfErr{u}{l} + \normbf{\Ph\bfErr{u}{l-1}}} \\\nonumber
				&\qquad+CC_AC_{\Ph}C_3\tau h^4 + CC_A[(\overline{C}_{2,2} + \overline{C}_{3,2})\tilde{K}_{1,4} + \widetilde{C}_2]\tau^5.
			\end{align}
		
			Applying the summation operator $\sum_{n=1}^{\ell}$ for any $1\leq \ell < N$ and the discrete Gronwall's inequality yields
			\begin{align*}\label{eqn:prf-fem-form-order-l2-9}
				\max_{1 \leq n \leq N}\exc{\cal{E}(\bfErr{u}{n})}\leq e^{CC_AT}\left[CC_AC_{\Ph}C_3 h^4 + CC_A[(\overline{C}_{2,2} + \overline{C}_{3,2})\tilde{K}_{1,4} + \widetilde{C}_2]\tau^4\right].
			\end{align*}
			The proof is complete.
		\end{proof}
		
		\section{Global error estimates}\label{sec-5}
		In this section, we combine the error estimates for Algorithms 1 and 2 to get global error estimates for the proposed fully discrete finite element methods.  They are immediate consequences     of Theorems \ref{theorem3.3}, \ref{theorem3.4}, \ref{theo:fem-form-order-h1}  \ref{theo:fem-form-order-l2},  and an application of the triangle inequality. 
		
		\begin{theorem}\label{theo:global-order-l2}
			Let $(\bu, \bv)$ be the variational solution of \eqref{eq2.4} and $\{(\bu^n_h, \bv^n_h)\}$ be the solutions generated Algorithm 2. Under the assumptions of Theorem \ref{theorem3.3} and Theorem \ref{theo:fem-form-order-l2}, there holds
			\begin{equation}
				\max_{1 \leq n \leq N}\exc{\normbf{\bu(t_n) - \bu^n_h}}\leq \tilde{C}_1\tau^3 + \tilde{C}_4\,h^4.
			\end{equation}
		\end{theorem}
		
		\begin{theorem}\label{theo:global--order-h1}
			Let $(\bu, \bv)$ be the variational solution of \eqref{eq2.4} and $\{(\bu^n_h, \bv^n_h)\}$ be the solutions generated Algorithm 2. Under the assumptions of Theorem \ref{theorem3.4} and Theorem \ref{theo:fem-form-order-h1}, there holds
			\begin{equation}
				\max\limits_{1\leq n \leq N}\exc{\normbfh{\bu(t_n)-\bu^n_h}+\normbf{\bv(t_n)-\bv^n_h}}\leq \tilde{C}_2\tau^2 + \tilde{C}_3h^2.
			\end{equation}
		\end{theorem}

		\section{Numerical tests}\label{sec-6}
		In this section, we provide two 2-D numerical tests to validate the proven rates of convergence. In our tests,  we choose the computational domain to be $\mathcal{D}=[-1,1]^2$ and use very fine spatial and temporal meshes to generate numerical exact solutions for verifying the rates of convergence. In addition, we define 
		\begin{align*}
			&\bL^{\infty}_t\bL^2_x(\bu):= \max_{1 \leq n \leq N}\bigl(\mE\bigl[\|\bu(t_n) - \bu^n\|^2_{\bL^2}\bigr]  \bigr)^{\frac12}, \\
			&\bL^{\infty}_t\bH^1_x(\bu):=\max_{1 \leq n \leq N}\bigl(\mE\bigl[\|\nabla(\bu(t_n) - \bu^n)\|^2_{\bL^2}\bigr]  \bigr)^{\frac12},\\
			&\bL^{\infty}_t\bL^2_x(\bv):=\max_{1 \leq n \leq N}\bigl(\mE\bigl[\|\bv(t_n) - \bv^n\|^2_{\bL^2}\bigr]  \bigr)^{\frac12}.
		\end{align*}
		
		{\bf Test 1.}  {\bf Linear case.}
				In this test, we consider linear drift and diffusion terms. Precisely, $\bF(\bu)=-[\bu,\ 3\bu]^T$ and $\bG[\bu]=[\bu,\ 3\bu]^T$. Table \ref{tab1_1} displays the computed errors in different norms and the computed rates of convergence. We observe that the temporal errors for $\bu^n$ and $\bv^n$ in the $\bL^2$-norm converge with rates $O(\tau^{\frac32})$ and $O(\tau)$, respectively, which matches well with our theoretical error estimates. The numerical results also suggest $O(\tau^{\frac32})$ rate of convergence for $\bu^n$ in the $\bH^1$-norm, which is not included in our theoretical results.
		
		\begin{table}[htb]
			\renewcommand{\arraystretch}{1.2}
			\centering 
			\footnotesize
			\begin{tabular}{|l||c|c||c|c||c|c|}
				\hline  
				&  $\bL^{\infty}_t\bL^2_x(\bu)$ error \quad & order &   
				$\bL^{\infty}_t\bH^1_x(\bu)$ error \quad & order &   
				$\bL^{\infty}_t\bL^2_x(\bv)$ error\quad & order \\ \hline    
				$\tau=1/8$ & $8.594\times10^{-3}$ & --- &   
				$5.404\times10^{-2}$ & --- &   
				$2.123\times10^{-2}$ & --- \\ \hline    
				$\tau=1/16$ & $2.316\times10^{-3}$ & 1.892 &   
				$1.464\times10^{-2}$ & 1.884 &   
				$1.030\times10^{-2}$ &1.004 \\ \hline    
				$\tau=1/32$ & $6.595\times10^{-4}$ & 1.812 &   
				$4.303\times10^{-3}$ & 1.766 &   
				$4.854\times10^{-3}$ & 1.085 \\ \hline    
				$\tau=1/64$ & $1.926\times10^{-4}$  & 1.776 &   
				$1.465\times10^{-3}$ & 1.555 &   
				$2.492\times10^{-3}$ & 0.962 \\ \hline       
			\end{tabular}
			\caption{Test 1: Temporal errors and convergence rates at $T = 1 $.}\label{tab1_1}
		\end{table}
		
		Table \ref{tab1_2} shows the errors and their convergence rates in the spatial variables which are  obtained using a very small time step size. We clearly observe the second order convergence for the  $\bL^{\infty}_t\bL^2_x(\bu)$-error,  the first order for the $\bL^{\infty}_t\bH^1_x(\bu)$-error,
		and the second order for $\bL^{\infty}_t\bL^2_x(\bv)$-error. 
		
		\begin{table}[!htbp]
			\centering 
			\footnotesize
			\begin{tabular}{|l||c|c||c|c||c|c|}
				\hline  
				&  $\bL^2$ error \quad & order &   
				$\bH^1$ error \quad & order &   
				$d_t\bL^2$ error\quad & order \\ \hline    
				$h=1/4$ & $6.821\times10^{-2}$ & --- &   
				$1.433\times10^{0}$ & --- &   
				$5.121\times10^{-1}$ & --- \\ \hline    
				$h=1/8$ & $1.438\times10^{-2}$ & 2.246 &   
				$7.035\times10^{-1}$ & 1.026 &   
				$1.141\times10^{-1}$ & 2.166\\ \hline    
				$h=1/16$ & $3.433\times10^{-3}$ & 2.067  &   
				$3.496\times10^{-1}$ & 1.009 &   
				$2.757\times10^{-2}$ & 2.049\\ \hline    
				$h=1/32$ & $8.481\times10^{-4}$  & 2.017 &   
				$1.745\times10^{-1}$ & 1.002 &   
				$6.834\times10^{-3}$ & 2.012 \\ \hline   
			\end{tabular}
			\caption{Test 1: Spatial errors and convergence rates when $\tau =
				1\times 10^{-3}$, $T = 0.1$.}
			\label{tab1_2}
		\end{table}
		
		\medskip
		{\bf Test 2.}  {\bf Nonlinear case.} 	
	  In this test, we consider the nonlinear drift and diffusion. Specifically, we take  $\bF[\bu]=[\cos(\bu),\ 2\cos(\bu)]^T$ and $\bG[\bu] =[\sin(\bu),\ 2\sin(\bu)]^T$. 
	  Table \ref{tab2_1} shows the computed errors in different norms and the computed convergence rates. We also observe that the temporal errors for $\bu^n$ in the $\bL^2$-norm and  the $\bH^1$-norm converge with rates $O(\tau^{\frac32})$ and $O(\tau)$, respectively. 
	  We also see that  $\bv^n$ converges in the $\bL^2$-norm with rate $O(\tau)$. All these results 
	  match our theoretical error estimates. 
		
		\begin{table}[htb]
			\renewcommand{\arraystretch}{1.2}
			\centering 
			\footnotesize
			\begin{tabular}{|l||c|c||c|c||c|c|}
				\hline  
				&  $\bL^{\infty}_t\bL^2_x(\bu)$ error \quad & order &   
				$\bL^{\infty}_t\bH^1_x(\bu)$ error \quad & order &   
				$\bL^{\infty}_t\bL^2_x(\bv)$ error\quad & order \\ \hline    
				$\tau=1/2$ & $8.703\times10^{-2}$ & --- &   
				$5.646\times10^{-1}$ & --- &   
				$9.058\times10^{-2}$ & --- \\ \hline    
				$\tau=1/4$ & $2.378\times10^{-2}$ & 1.872 &   
				$1.592\times10^{-1}$ & 1.826 &   
				$3.983\times10^{-2}$ & 1.185 \\ \hline    
				$\tau=1/8$ & $7.970\times10^{-3}$ & 1.577 &   
				$5.299\times10^{-2}$ & 1.587 &   
				$1.961\times10^{-2}$ & 1.022 \\ \hline    
				$\tau=1/16$ & $2.703\times10^{-3}$  & 1.560 &   
				$1.827\times10^{-2}$ & 1.537 &   
				$9.420\times10^{-3}$ & 1.058 \\ \hline       
			\end{tabular}
			\caption{Test 2: Temporal errors and convergence rates at $T=1 $.}\label{tab2_1}
		\end{table}
		
		Table \ref{tab2_2} shows the errors and their convergence rates in spatial variables. They are the second order convergence for the  $\bL^{\infty}_t\bL^2_x(\bu)$-error,  the first order for the $\bL^{\infty}_t\bH^1_x(\bu)$-error, and the second order for $\bL^{\infty}_t\bL^2_x(\bv)$-error.

		\begin{table}[!htbp]
			\centering 
			\footnotesize
			\begin{tabular}{|l||c|c||c|c||c|c|}
				\hline  
				&  $\bL^2$ error \quad & order &   
				$\bH^1$ error \quad & order &   
				$d_t\bL^2$ error\quad & order \\ \hline    
				$h=1/4$ & $1.540\times10^{-2}$ & --- &   
				$4.013\times10^{-1}$ & --- &   
				$3.695\times10^{-2}$ & --- \\ \hline    
				$h=1/8$ & $3.645\times10^{-3}$ & 2.079 &   
				$1.982\times10^{-1}$ & 1.018 &   
				$8.847\times10^{-3}$ & 2.062\\ \hline    
				$h=1/16$ & $8.985\times10^{-4}$ & 2.020 &   
				$9.881\times10^{-2}$ & 1.004 &   
				$2.183\times10^{-3}$ & 2.019 \\ \hline    
				$h=1/32$ & $2.238\times10^{-4}$  & 2.005 &   
				$4.937\times10^{-2}$ & 1.001 &   
				$5.440\times10^{-4}$ & 2.005 \\ \hline       
			\end{tabular}
			\caption{Test 1: Spatial errors and convergence rates when $\tau =
				1\times 10^{-3}$, $T = 0.1$.}
			\label{tab2_2}
		\end{table}
		
		
	\end{document}